\documentclass[A4paper,10pt]{article}

\usepackage[numbers]{natbib}
\usepackage{graphicx}
\usepackage{algorithm,algorithmic}
\usepackage{setspace}
\usepackage{rotating}
\usepackage{color}
\usepackage{amsmath}
\usepackage{enumerate}
\usepackage{hyperref}
\usepackage{authblk}
\usepackage{geometry}
\usepackage{indentfirst}
\usepackage{amssymb}
\usepackage{subfigure}
\usepackage{nicefrac}
\usepackage{bm}
\usepackage{multirow}
\usepackage{amsfonts,amsthm}
\usepackage[title]{appendix}
\usepackage{epstopdf}
\usepackage{mathrsfs}
\newcommand{\dt}{\Delta t}

\newcommand{\U}{\bm{U}}

\newcommand{\F}{\bm{F}}
\newcommand{\G}{\bm{G}}

\newcommand{\w}{\bm{\omega}}

\newcommand{\T}{\bm{T}}
\newcommand{\bH}{\bm{H}}
\newcommand{\mbH}{\mathcal{\bH}}
\newcommand{\z}{\bm{z}}
\newcommand{\n}{\bm{n}}

\newcommand{\mG}{\mathcal{G}}

\newcommand{\half}{\frac{1}{2}}

\newcommand{\jph}{{j+\frac{1}{2}}}

\newcommand{\kph}{{k+\frac{1}{2}}}

\newcommand{\dsum}{\displaystyle\sum}
\newcommand{\CFL}{{\rm CFL}}
\newcommand{\dvert}[1]{\left\vert #1 \right\vert}
\newcommand{\dVert}[1]{\left\Vert #1 \right\Vert}
\newcommand\chfrac[2]{\frac{\; #1\;}{\;#2\;}}

\newtheorem{remark}{Remark}
\newtheorem{lemma}{Lemma}
\newtheorem{theorem}{Theorem}
\newtheorem{example}{Example}



\title{Bound-preserving Flux Limiters for Hyperbolic Conservation Laws on Two-Dimensional Adaptive Moving Meshes}

\author[1]{Yaguang Gu
\thanks{guyaguang@scut.edu.cn}}

\author[2]{Guanghui Hu
\thanks{Corresponding author. garyhu@um.edu.mo}}

\author[3]{Mengqing Liu
\thanks{mengqingliu@m.scnu.edu.cn}}

\affil[1]{School of Mathematics, South China University of Technology, Guangzhou 510641, China}

\affil[2]{State Key Laboratory of Internet of Things for Smart City and Department of Mathematics, University of Macau, Macao SAR 999078, China}

\affil[3]{School of Mathematical Sciences, South China Normal University, Guangzhou 510631, China}

\date{}

\begin{document}

\maketitle

\begin{abstract}
In this paper, we develop a bound-preserving finite-volume method on adaptive moving meshes for two-dimensional hyperbolic conservation laws. The proposed method extends our previous one-dimensional framework [Y. Gu, G. Hu, and T. Tang, Adaptive mesh methods for hyperbolic conservation laws with bound-preserving flux limiters, arXiv:2602.12580] to two dimensions within the arbitrary Lagrangian--Eulerian framework. The mesh velocity is determined by solving a moving mesh equation, allowing the numerical solution to evolve directly between consecutive time levels. A key challenge in this extension is to ensure compatibility between the evolution of the mesh geometry and the update of the numerical solution, as geometric information is essential to the convex decomposition underlying the bound-preserving construction in two dimensions. Specifically, the high-order scheme is reformulated as a convex combination of high-order sub-cell schemes, each of which is blended with its first-order counterpart to enforce the bound-preserving property. The discrete geometric conservation law ensures the compatibility needed to maintain this convex decomposition on moving meshes. The resulting correction acts locally at individual cell interfaces and is equivalent to applying a bound-preserving flux limiter, making the approach applicable to general unstructured or irregular meshes. The associated bound-preserving CFL conditions are established for both scalar conservation laws and the nonlinear Euler equations of gas dynamics. Extensive two-dimensional numerical experiments demonstrate the accuracy, robustness, and effectiveness of the proposed method.
\end{abstract}



\section{Introduction}
Hyperbolic conservation laws arise in a wide range of scientific and engineering applications, including compressible fluid dynamics, shallow water flows, multi-material flows, etc. Numerical simulation of these problems is particularly challenging because the solutions may contain sharp gradients or even discontinuities whose locations are generally unknown in advance. High-resolution numerical methods capable of resolving such local structures therefore remain an active topic in computational physics.

Adaptive moving-mesh methods provide an effective approach to this challenge by dynamically redistributing mesh points according to the evolving solution. Compared with uniform-mesh finite-volume methods (UM-FVMs), adaptive moving-mesh finite-volume methods (AMM-FVMs) concentrate computational effort in regions of interest while maintaining a fixed number of cells and the topology of the mesh, thereby achieving substantially higher resolution than UM-FVMs. Besides, for a prescribed level of resolution, AMM-FVMs require considerably fewer degrees of freedom than UM-FVMs, which represents a stronger measure of efficiency than the former one. These advantages make adaptive moving meshes an attractive framework for high-resolution simulations of hyperbolic conservation laws.
Motivated by these advantages, we develop a third-order finite volume method on adaptive moving meshes for the two-dimensional hyperbolic conservation law
\begin{equation}\label{EQ:Govern-Eq-Scalar}
u_t + \nabla\cdot \bH(u) = 0,
\end{equation}
where $\bH(u) = (f(u),\, g(u))$ is a flux tensor.

In addition to achieving high resolution, numerical methods should preserve the physical admissibility of the solution. For scalar conservation laws, the exact solution satisfies a maximum principle, whereas for nonlinear systems one seeks to preserve suitable invariant domains, such as positivity of density and pressure or boundedness of volume fractions. Violation of these physical constraints may destroy hyperbolicity of the governing equation and lead to numerical instability or even simulation failure. Despite the advantages of AMM-FVMs, existing high-order AMM-FVMs generally do not guarantee preservation of the invariant domain. Developing high-order bound-preserving (BP) AMM-FVMs therefore remains an important yet challenging problem.

During the past decades, a variety of BP techniques have been developed. 
One representative framework is the Zhang-Shu bound-preserving limiter \cite{ZHANG20103091}, which combines reconstruction process with a simple scaling limiter, and the BP property is achieved by using a provably suitable BP CFL number. The methodology has been successfully extended to the Euler equations, balance laws, and many other hyperbolic problems; consult \cite{CUI2025114189,
FAN2022111446,QIAN2018172,Wu2021Uniformly,
ZHANG20108918,ZHANG20111238,ZHANG20122245} and references therein. It has also been adapted to moving-mesh discretizations; see e.g. \cite{CHENG2014143,
Fu2019Arbitrary,ZHANG2025113801}.
Another important class of BP methods is based on flux limiting, where high-order numerical fluxes are blended with their monotone first-order counterparts. To the best of our knowledge, this idea originates from the classical flux-corrected transport (FCT) methods and has been developed into various forms; consult \cite{BOOK1975248,BORIS1976397,BORIS197338,
kuzmin2012flux,ZALESAK1979335} and references therein. In \cite{XU2014Parametrized}, a parametrized flux limiter was proposed, where the BP property is fulfilled by using linear programming method to determine flux limiters on all edges of the cell; see also \cite{Liang2014Parametrized,
Xiong2016Parametrized} for subsequent studies. Although parametrized limiters admit relatively mild CFL restrictions, their implementation on irregular or moving meshes remains challenging.
Recently, monolithic convex limiting (MCL) BP method has provided a flexible alternative by introducing flux limiters  \cite{KUZMIN2020112804}. In this method, since the BP flux limiters are only dependent on a single cell interface, the resulting algorithms are considerably easier to implement on adaptive moving meshes.
For studies of the MCL flux limiters, we refer the readers to \cite{Abgrall2025bound,
Duan2025Active,Fu2025Bound,Gao2023High,
GU2023112190,Gu2021A,Kuzmin2023Property,
Wu2025High} and references therein. In \cite{Gu2026Adaptive}, a novel third-order BP AMM-FVM was proposed for one-dimensional hyperbolic problems. The method generalizes the uniform-mesh MCL BP schemes to adaptive moving meshes. In addition, it is first shown that the accuracy-preserving CFL condition of the third-order AMM-FVMs is $\nicefrac{1}{6}$.

Extending this framework to two-dimensional adaptive moving meshes is far from straightforward because mesh motion, geometric conservation, and high-order reconstruction become strongly coupled.
In this paper, we propose a third-order BP AMM-FVM for two-dimensional hyperbolic problems. For clarity of presentation, we consider a third-order AMM-FVM equipped with the classical Riemann-solver-free central-upwind (CU) flux \cite{Kurganov2001Semi} and the third-order strong-stability-preserving Runge-Kutta (SSP-RK3) time discretization \cite{GOTTLIEB1998Total,Shu2001Strong}. 
One main challenge is that the discrete geometric conservation law (D-GCL), which guarantees that a uniform flow remains invariant under temporal evolution, must be satisfied by the AMM-FVM at the discrete level. In addition, the D-GCL will also play an important role in the BP procedure: it allows decomposition of the AMM-FVM into a convex combination of sub-cell schemes. In other words, the D-GCL should not be violated throughout BP procedure at any time step.
Another difficulty comes from the convex decomposition mentioned above, as the mesh geometry becomes irregular in two dimensions. Computation of the geometry information and reformulation of the AMM-FVM in BP procedure become highly non-trivial.

Taking the challenges mentioned above into account, we first systematically propose a BP AMM-FVM framework for scalar problems. The main ingredients of our BP procedure are summarized as follows:
\begin{enumerate}[(i)]
\item At every SSP-RK3 stage, we deliver a convex decomposition of the scheme into a linear combination of sub-cell schemes. Subsequently, we define first-order counterparts of these sub-cell schemes. 

\item Sufficient CFL conditions are derived to guarantee both the convexity of the decomposition and the BP property of each first-order sub-cell scheme. 

\item Whenever a high-order sub-cell scheme violates the physical constraints, it is blended with its first-order counterpart through a local correction procedure. The resulting algorithm naturally yields a simple BP flux limiter involving only one cell interface, making the implementation straightforward. Meanwhile, the D-GCL is preserved after BP corrections, which implies that the BP procedure can be successfully carried out at any time step.
\end{enumerate}
The proposed framework is further extended to the compressible Euler system of gas dynamics. By utilizing rotational invariance, we establish practical BP CFL conditions that depend solely on the first-order sub-cell schemes.

The remainder of this paper is organized as follows. In \S\ref{Sec:BP-FVM-Scalar}, we propose our BP AMM-FVM for scalar conservation laws and derive the corresponding BP CFL conditions. In \S\ref{Sec:Euler},  the method is extended to the compressible Euler equations. Numerical experiments are reported in \S\ref{Sec:Num} to verify the theoretical analysis and demonstrate the accuracy, robustness, and efficiency of the proposed methods.

\section{Bound-preserving AMM-FVM for scalar hyperbolic equations}\label{Sec:BP-FVM-Scalar}

This section elaborates the BP AMM-FVM scheme for solving the scalar hyperbolic conservation laws \eqref{EQ:Govern-Eq-Scalar}. First, \S\ref{Sec:FVM} presents a classical finite-volume method on adaptive moving meshes, which evolves the solution directly from $t^n$ to $t^{n+1}$. To further enhance the numerical robustness of the AMM-FVM, a bound-preserving correction is developed in \S\ref{Sec:BP-CFL}.
The mechanism of this BP scheme is comprehensively introduced, and the corresponding BP CFL constraints are rigorously derived. Finally, \S\ref{Sec:mesh-redist} briefly reviews the adaptive mesh redistribution algorithm, which will be used to define the mesh velocity. Meanwhile, numerical implementation details are discussed.

\subsection{Finite-volume discretization}\label{Sec:FVM}
Let the physical domain $\Omega$ be partitioned by a time-dependent structured irregular grid consisting of quadrilateral finite-volume cells $C_{j,k}(t)$, $j=1,2,\dots,N_x$ and $k=1,2,\dots,N_y$. For each cell $C_{j,k}(t)$, $\z_{j\pm\half,k\pm\half}$ denote coordinates of its vertices, $\dvert{C_{j,k}(t)}$ denotes its cell area, and $\z_{j,k} = (x_{j, k}, y_{j, k})^\top$ denotes coordinates of its centroid. Although the mesh points vary temporally, each node is assumed to move along a straight trajectory within a single time step, governed by the mesh velocity $\w(\z)$, which is determined by the moving mesh equation and will be introduced in \S\ref{Sec:mesh-redist}.

Assume that the solution, realized in terms of its cell average $\overline{u}_{j,k}$ over $C_{j,k}$, is available at a certain time level $t^n$. We discretize the scalar hyperbolic equation \eqref{EQ:Govern-Eq-Scalar} via a finite-volume spatial discretization combined with the third-order strong-stability-preserving Runge-Kutta (SSP-RK3) temporal discretization, yielding the fully discrete AMM-FVM scheme as follows.
\begin{equation}\label{EQ:AMM-FVM-U}
\begin{aligned}
\overline{u}_{j,k}^{(1)} &= 
\chfrac{\dvert{C_{j,k}^n}}{\dvert{C_{j,k}^{(1)}}}\overline{u}_{j,k}^n - \chfrac{\dt^n}{\dvert{C_{j,k}^{(1)}}} \sum_{s=1}^{4} L_{j,k,s}^{n} \sum_{r=1}^{N_G}\eta_r
\mbH \left(\w_{j,k,s,r},\n_{j,k,s}^n, u_{j,k,s,r}^{n,-},u_{j,k,s,r}^{n,+}\right), \vspace{+3pt}\\
\overline{u}_{j,k}^{(2)} &= \chfrac{3}{4}\frac{\dvert{C_{j,k}^n}}{\dvert{C_{j,k}^{(2)}}}\overline{u}_{j,k}^n + \frac{1}{4}\left[
\chfrac{\dvert{C_{j,k}^{(1)}}}{\dvert{C_{j,k}^{(2)}}}\overline{u}_{j,k}^{(1)} -
\chfrac{\dt^n}{\dvert{C_{j,k}^{(2)}}} \sum_{s=1}^{4} L_{j,k,s}^{(1)} \sum_{r=1}^{N_G}\eta_r
\mbH \left(\w_{j,k,s,r}, \n_{j,k,s}^{(1)}, u_{j,k,s,r}^{(1),-},u_{j,k,s,r}^{(1),+}\right)\right],\vspace{+3pt}\\
\overline{u}_{j,k}^{n+1} &= \frac{1}{3}\chfrac{\dvert{C_{j,k}^n}}{\dvert{C_{j,k}^{n+1}}}\overline{u}_{j,k}^n + \frac{2}{3}\left[
\chfrac{\dvert{C_{j,k}^{(2)}}}{\dvert{C_{j,k}^{n+1}}}\overline{u}_{j,k}^{(1)} -
\chfrac{\dt^n}{\dvert{C_{j,k}^{n+1}}} \sum_{s=1}^{4} L_{j,k,s}^{(2)} \sum_{r=1}^{N_G}\eta_r
\mbH \left(\w_{j,k,s,r}, \n_{j,k,s}^{(2)}, u_{j,k,s,r}^{(2),-},u_{j,k,s,r}^{(2),+}\right)\right].
\end{aligned}
\end{equation}
In the above formulation, $L_{j,k,s}$ denotes the length of the $s$-th edge of $C_{j,k}$, $\n_{j,k,s}=[\cos\theta_{j,k,s},\sin\theta_{j,k,s}]^\top$ denotes the  outward unit normal vector to the edge $s$, and $N_G$ is the number of Gauss quadrature points adopted for boundary integration, with associated weights denoted by $\eta_r$, The quantities $u_{j,k,s,r}^{(\ell),-}$ and $u_{j,k,s,r}^{(\ell),+}$ for $\ell = 0, 1, 2$ (with $u_{j,k,s,r}^{(0),\pm}=u_{j,k,s,r}^{n,\pm}$) stand for the interior and exterior reconstructed states at each Gauss quadrature point on the cell interface, obtained via weighted essentially non-oscillatory (WENO) reconstruction. In this work, we implement the $k$-exact WENO reconstruction proposed in \cite{HU2016235} and adapt it to general quadrilateral meshes. For the neighboring control volume $C_{\widehat{j},\widehat{k}}$ that shares the $s$-th edge with $C_{j,k}$, the relevant variables satisfy the following relations:
\begin{equation}\label{EQ:quant-adj-vol}
L_{j,k,s} = L_{\widehat{j},\widehat{k},s},\quad
\n_{j,k,s} = -\n_{\widehat{j},\widehat{k},s}, \quad
\w_{j,k,s,r} = \w_{\widehat{j},\widehat{k},s,r},\quad
u_{j,k,s,r}^{(\ell),\pm} = u_{\widehat{j},\widehat{k},s,r}^{(\ell),\mp},
\end{equation}
The operator $\mbH$ is responsible for evaluating numerical fluxes at cell interfaces. In this paper, we adopt the Riemann-solver-free central-upwind (CU) flux proposed in \cite{Kurganov2001Semi}, which takes the form
\begin{equation*}
\mbH( \w, \n, u^-, u^+ ) = \begin{cases}
\bH(u^-)\cdot\n - \omega u^-, &\omega \leq a^-, \vspace{+3pt}\\
\dfrac{a^+ \bH(u^-)\cdot\n - a^- \bH(u^+)\cdot\n}{a^+ - a^-} + \dfrac{a^- a^+}{a^+ - a^-}(u^+ - u^-) - \omega u^*, &a^-< \omega < a^+, \vspace{+3pt}\\
\bH(u^+)\cdot\n - \omega u^+, &\omega \geq a^+,
\end{cases}
\end{equation*}
where $\omega=\w\cdot\n$, and $\bH(u)\cdot\n = f(u)\cos\theta + g(u)\sin\theta$, with $\n=(\cos\theta, \sin\theta)^\top$. The local speeds of propagation $a^\pm$ are estimated by
\begin{equation*}
a^- = \min\left\lbrace \bH'(u^-)\cdot\n,~\bH'(u^+)\cdot\n,~ -\varepsilon^{\rm CU} \right\rbrace,\quad
a^+ = \max\left\lbrace \bH'(u^-)\cdot\n,~\bH'(u^+)\cdot\n,~\varepsilon^{\rm CU} \right\rbrace.
\end{equation*}
Here, $\bH'(u) = (f'(u),g'(u))$ denotes the Jacobian of the flux vector, and a small regularization parameter $\varepsilon^{\rm CU} = 10^{-12}$ is introduced to avoid division by zero in the CU flux. In addition, $u^*$ represents the intermediate state in the CU discretization, and takes the form 
\begin{equation*}
u^* = \dfrac{a^+ u^+ - a^- u^- - \bH(u^+)\cdot\n + \bH(u^-)\cdot\n}{a^+ - a^-}.
\end{equation*}

The time step size is constrained by the standard CFL condition
\begin{equation}\label{EQ:Conventional-dt}
\dt^n \max_{j,k} \left\{ \frac{\max\left\{ \dvert{a^-_{j+\frac{1}{2},k}}, \dvert{a^+_{j-\frac{1}{2},k}}, \dvert{a^-_{j,k+\frac{1}{2}}}, \dvert{a^+_{j,k-\frac{1}{2}}} \right\}}{\mathrm{dist}(\z_{j,k}, \partial C_{j,k})} \right\} \leq \CFL,
\end{equation}
where the conventional CFL number is set to $\CFL=1$.

During the time evolution from $t^n$ to $t^{n+1}$, the time-dependent intermediate areas of moving control volumes require consistent updating, which is accomplished using the same SSP-RK3 temporal discretization. The resulting area update formula reads as
\begin{equation}\label{EQ:AMM-FVM-C}
\begin{aligned}
\dvert{C_{j,k}^{(1)}} &= 
\dvert{C_{j,k}^n} + \dt^n \sum_{s=1}^{4} L_s^n \sum_{r=1}^{N_G}\eta_r
(\w_{j,k,s,r}\cdot\n_{j,k,s}^n), \vspace{+3pt}\\
\dvert{C_{j,k}^{(2)}} &= \frac{3}{4}\dvert{C_{j,k}^n} + \frac{1}{4}\left[
\dvert{C_{j,k}^{(1)}} + 
\dt^n\sum_{s=1}^{4} L_s^{(1)} \sum_{r=1}^{N_G}\eta_r
(\w_{j,k,s,r}\cdot\n_{j,k,s}^{(1)})\right],\vspace{+3pt}\\
\dvert{C_{j,k}^{n+1}} &= \frac{1}{3}\dvert{C_{j,k}^n} + \frac{2}{3}\left[\dvert{C_{j,k}^{(2)}} + \dt^n \sum_{s=1}^{4} L_s^{(2)} \sum_{r=1}^{N_G}\eta_r
(\w_{j,k,s,r}\cdot\n_{j,k,s}^{(2)})\right].
\end{aligned}
\end{equation}
This updating strategy strictly satisfies the discrete geometric conservation law (D-GCL), which guarantees that uniform flow remains invariant under temporal evolution. Furthermore, the D-GCL serves as a fundamental prerequisite for the subsequent development of our BP scheme, enabling the convex decomposition of the fully discrete scheme \eqref{EQ:AMM-FVM-U} into a combination of sub-cell schemes.

\subsection{Bound-preserving method and CFL conditions}\label{Sec:BP-CFL}

In this section, we derive a simple bound-preserving numerical framework for two-dimensional hyperbolic conservation laws \eqref{EQ:Govern-Eq-Scalar}, and further establish the corresponding BP CFL conditions. The construction and theoretical analysis of the proposed BP scheme rely on several essential lemmas, which are presented below.

\begin{lemma}\label{Lem:Consistency-num-flux}
Let $u$ denote the conservative variable, and let $\w$ and $\n$ represent the mesh velocity and unit outward normal vector, respectively. The following relation holds:
\begin{equation*}
\mbH(\w,\n, u,u) = \bH(u)\cdot\n - (\w\cdot\n)\, u.
\end{equation*}
In particular, the relation reduces to the following form when the mesh points are static, namely $\w = \bm{0}$.
\begin{equation*}
\mbH(\bm{0},\n,u,u) = \bH(u)\cdot\n.
\end{equation*}
\end{lemma}
\begin{proof}
The above results follow directly from the consistency property of the numerical flux.
\end{proof}

\begin{lemma}\label{Lem:HmH}
Let $u$ denote the conservative variable, with $\w$ and $\n$ denoting the mesh velocity and unit outward normal vector, respectively. The identity below holds:
\begin{equation}\label{EQ:HmH}
\mbH(\w,\n, u,u) - \mbH(\bm{0},\n,u,u) = -(\w\cdot\n)\, u.
\end{equation}
\end{lemma} 
\begin{proof}
According to Lemma \ref{Lem:Consistency-num-flux}, we obtain
\begin{equation*}
\mbH(\w,\n, u,u) - \mbH(\bm{0},\n,u,u) = \left[\bH(u)\cdot\n - (\w\cdot\n)\, u\right] - \left[ \bH(u)\cdot\n \right] = -(\w\cdot\n)\, u.
\end{equation*}
This completes the proof of Lemma \ref{Lem:HmH}.
\end{proof}

\begin{lemma}\label{Lem:Inout-Flux}
Let $u$ and $v$ be the conservative variables, and let $\w$ and $\n$ denote the mesh velocity and unit outward normal vector, respectively. The flux operator satisfies the following relation:
\begin{equation}
\mbH(\w,\n, u, v) = -\mbH(\w,-\n, v, u).
\end{equation}
\end{lemma}
\begin{proof}
This property holds trivially, as the numerical flux from $C_{j,k}$ to $C_{\widehat{j},\widehat{k}}$ is exactly opposite to the flux from $C_{\widehat{j},\widehat{k}}$ to $C_{j,k}$.
\end{proof}

\begin{lemma}\label{Lem:sum-Ln}
Let $L_{j,k,s}$ ($s = 1, \dots, 4$) be the length of the $s$-th edge of the control volume $C_{j,k}$, and let $\n_{j,k,s}$ denote the corresponding unit outward normal vector. The following geometric identity holds:
\begin{equation}\label{EQ:sum-Ln}
\sum_{s = 1}^{4} L_{j,k,s} \n_{j,k,s} = \bm{0}.
\end{equation}
\end{lemma}

\begin{proof}
Define the edge vector $\bm{L}_{j,k,s}$ formed by connecting two adjacent vertices of the $s$-th edge of $C_{j,k}$ in the clockwise direction. The magnitude of $\bm{L}_{j,k,s}$ is $L_{j,k,s}$. For a quadrilateral cell, it follows immediately that $\sum_{s = 1}^{4} \bm{L}_{j,k,s} = \bm{0}$. Rotating each vector $\bm{L}_{j,k,s}$ counterclockwise by $90^\circ$ produces the vector $L_{j,k,s}\n_{j,k,s}$. Consequently, the identity \eqref{EQ:sum-Ln} stills holds, which completes the proof of Lemma \ref{Lem:sum-Ln}.
\end{proof}

\begin{lemma}\label{Lem:Sum-w}
Let $\w_{j,k,s,r}$ denote the mesh velocity evaluated at the Gauss quadrature point $\z_r$ on the $s$-th edge of $C_{j,k}(t^n)$, and let $\w_{j,k,s}$ denote the mesh velocity at the midpoint of the same edge. Then the following relation holds:
\begin{equation*}
\sum_{r=1}^{N_G} \eta_r \w_{j,k,s,r} = \w_{j,k,s}.
\end{equation*}
\end{lemma}

\begin{proof}
Denote the two endpoints of the $s$-th edge of $C_{j,k}(t)$ by $\z_0(t)$ and $\z_1(t)$, and let $\z_{\nicefrac{1}{2}}(t)$ be the corresponding midpoint. The Gauss quadrature points on each edge are symmetrically distributed about $\z_{\nicefrac{1}{2}}$, and symmetric quadrature points share identical integration weights. Thus, it suffices to verify the following identity for an arbitrary pair of symmetric quadrature points $(\z_\alpha, \z_\beta)$ with common weight $\eta$:
\begin{equation*}
\eta \w_\alpha + \eta \w_\beta =2\eta\w_\half.
\end{equation*}
To proceed, we parameterize the coordinates of $\z_\alpha$ and $\z_\beta$ as
\begin{equation}\label{EQ:z-alp-beta}
\z_\alpha = \zeta \z_0 + (1-\zeta)\z_1, \qquad
\z_\beta = (1-\zeta) \z_0 + \zeta\z_1, 
\end{equation}
for a certain parameter $\zeta\in (0,1)$. The associated mesh velocities are expressed as
\begin{equation}\label{EQ:w-alp-beta}
\w_\alpha = \frac{\z_\alpha^{n+1} - \z_\alpha^n}{\dt^n}, \qquad
\w_\beta  = \frac{\z_\beta^{n+1}  - \z_\beta^n}{\dt^n}.
\end{equation}
Substituting \eqref{EQ:z-alp-beta} and \eqref{EQ:w-alp-beta} into the summation formula yields
\begin{equation*}
\eta \w_\alpha + \eta \w_\beta = \eta \left( \frac{\z_\alpha^{n+1} - \z_\alpha^n}{\dt^n} + \frac{\z_\beta^{n+1}  - \z_\beta^n}{\dt^n} \right) = 2\eta \frac{\frac{\z_0^{n+1}+\z_1^{n+1}}{2}-\frac{\z_0^n + \z_1^n}{2}}{\dt^n} = 2\eta\w_\half.
\end{equation*}
which completes the proof of Lemma \ref{Lem:Sum-w}.
\end{proof}

With the above lemmas established, we now proceed to present our BP scheme. Consider the $\ell$-th stage of the SSP-RK3 method, which admits a unified formulation as
\begin{equation}\label{EQ:RK-stage-l}
\overline{u}_{j,k}^{(\ell)} = \xi^\ell\chfrac{\dvert{C_{j,k}^n}}{\dvert{C_{j,k}^{(\ell)}}}\overline{u}_{j,k}^n + (1-\xi^\ell)\left[
\chfrac{\dvert{C_{j,k}^{(\ell-1)}}}{\dvert{C_{j,k}^{(\ell)}}}\overline{u}_{j,k}^{(\ell-1)} -
\chfrac{\dt^n}{\dvert{C_{j,k}^{(\ell)}}} \sum_{s=1}^{4} L_{j,k,s}^{(\ell-1)} \sum_{r=1}^{N_G}\eta_r
\mbH \left(\w_{j,k,s,r}, \n_{j,k,s}^{(\ell-1)}, u_{j,k,s,r}^{(\ell-1),-},u_{j,k,s,r}^{(\ell-1),+}\right)\right],
\end{equation}
where the intermediate cell area $\dvert{C_{j,k}^{(\ell)}}$ is formulated uniformly by
\begin{equation}\label{EQ:RK-stage-Cjk-l}
\dvert{C_{j,k}^{(\ell)}} = \xi^\ell\dvert{C_{j,k}^n} + (1-\xi^\ell)\left[\dvert{C_{j,k}^{(\ell-1)}} + \dt^n \sum_{s=1}^{4} L_{j,k,s}^{(\ell-1)} \sum_{r=1}^{N_G}\eta_r
(\w_{j,k,s,r}\cdot\n_{j,k,s}^{(\ell-1)})\right].
\end{equation}
Here, $\ell = 1, 2, 3$ corresponds to the three Runge-Kutta stages, with $\overline{u}_{j,k}^{(0)}=\overline{u}_{j,k}^{n}$ and $\overline{u}_{j,k}^{(3)}=\overline{u}_{j,k}^{n+1}$. The coefficients $\xi^{(\ell)}$ are $\xi^{(1)} = 0$, 
$\xi^{(2)} = \frac{3}{4}$,
and $\xi^{(3)} = \frac{1}{3}$.
 
Based on Lemma \ref{Lem:Consistency-num-flux}, Lemma \ref{Lem:HmH}, Lemma \ref{Lem:sum-Ln}, and \eqref{EQ:RK-stage-Cjk-l}, we utilize the identity  
\begin{equation*}
\sum_{s=1}^{Ns}L_{j,k,s}^{(\ell-1)}\bH(\overline{u}_{j,k}^{(\ell-1)})\cdot\n_{j,k,s}^{(\ell-1)} = \sum_{s=1}^{4}L_{j,k,s}^{(\ell-1)}\mbH\left(\bm{0},\n_{j,k,s}^{(\ell-1)},\overline{u}_{j,k}^{(\ell-1)},\overline{u}_{j,k}^{(\ell-1)}\right) = 0
\end{equation*}
to reformulate \eqref{EQ:RK-stage-l} into an equivalent form
\begin{equation}\label{EQ:RK-stage-l-temp}
\overline{u}_{j,k}^{(\ell)} = \xi^\ell\chfrac{\dvert{C_{j,k}^n}}{\dvert{C_{j,k}^{\ell}}}\overline{u}_{j,k}^n + (1-\xi^\ell)\left[
\chfrac{\dvert{C_{j,k}^{(\ell-1)}}}{\dvert{C_{j,k}^{\ell}}}\overline{u}_{j,k}^{(\ell-1)} -
\chfrac{\dt^n}{\dvert{C_{j,k}^{\ell}}} \sum_{s=1}^{4} L_{j,k,s}^{(\ell-1)} \Delta \mbH_{j,k,s} \right],
\end{equation}
where the difference of the flux are defined and further reformulated as
\begin{equation*}
\begin{aligned}
\Delta \mbH_{j,k,s} =& \sum_{r=1}^{N_G} \eta_r\mbH \left(\w_{j,k,s,r}, \n_{j,k,s}^{(\ell-1)}, u_{j,k,s,r}^{(\ell-1),-},u_{j,k,s,r}^{(\ell-1),+}\right) - \mbH\left(\bm{0},\n_{j,k,s}^{(\ell-1)},\overline{u}_{j,k}^{(\ell-1)},\overline{u}_{j,k}^{(\ell-1)}  \right)\\
%%=& \sum_{r=1}^{N_G} \eta_r \left[\mbH \left(\w_{j,k,s,r}, \n_{j,k,s}^{(\ell-1)}, u_{j,k,s,r}^{(\ell-1),-},u_{j,k,s,r}^{(\ell-1),+}\right) - \mbH\left(\bm{0},\n_{j,k,s}^{(\ell-1)},\overline{u}_{j,k}^{(\ell-1)},\overline{u}_{j,k}^{(\ell-1)}  \right)\right]\\
=& \sum_{r=1}^{N_G} \eta_r \left[\mbH \left(\w_{j,k,s,r}, \n_{j,k,s}^{(\ell-1)}, u_{j,k,s,r}^{(\ell-1),-},u_{j,k,s,r}^{(\ell-1),+}\right) - \mbH\left(\w_{j,k,s,r},\n_{j,k,s}^{(\ell-1)},\overline{u}_{j,k}^{(\ell-1)},\overline{u}_{j,k}^{(\ell-1)}\right) \right.\\
& \left.+ \mbH\left(\w_{j,k,s,r},\n_{j,k,s}^{(\ell-1)},\overline{u}_{j,k}^{(\ell-1)},\overline{u}_{j,k}^{(\ell-1)}\right) - \mbH\left(\bm{0},\n_{j,k,s}^{(\ell-1)},\overline{u}_{j,k}^{(\ell-1)},\overline{u}_{j,k}^{(\ell-1)}\right)\right] \\
=& \sum_{r=1}^{N_G} \eta_r \left[\mbH \left(\w_{j,k,s,r}, \n_{j,k,s}^{(\ell-1)}, u_{j,k,s,r}^{(\ell-1),-},u_{j,k,s,r}^{(\ell-1),+}\right) - \mbH\left(\w_{s,r},\n_s^{(\ell-1)},\overline{u}_{j,k}^{(\ell-1)},\overline{u}_{j,k}^{(\ell-1)}\right)\right] - \sum_{r=1}^{N_G} \eta_r \omega_{j,k,s,r}^{(\ell-1)}\overline{u}_{j,k}^{(\ell-1)}.
\end{aligned}
\end{equation*}
Substituting $\Delta \mbH_{j,k,s}$ into \eqref{EQ:RK-stage-l-temp}, we obtain a decomposition of the high-order scheme \eqref{EQ:RK-stage-l} into a linear combination of sub-cell schemes
\begin{equation}\label{EQ:Convex-Decomp}
\overline{u}_{j,k}^{(\ell)} = \xi^\ell\chfrac{\dvert{C_{j,k}^n}}{\dvert{C_{j,k}^{(\ell)}}}\overline{u}_{j,k}^n + (1-\xi^\ell)\chfrac{\widehat{C}_{j,k}^{(\ell-1)}}{\dvert{C_{j,k}^{(\ell)}}}\sum_{s=1}^{4}\frac{L_{j,k,s}^{(\ell-1)}}{\sum_{q=1}^{4}L_{j,k,q}^{(\ell-1)}} u_{j,k,s}^{(\ell-1),\rm H},
\end{equation}
where $\widehat{C}_{j,k}^{(\ell-1)}$ is a rough estimation of the intermediate cell area by forward Euler, namely 
\begin{equation}\label{EQ:hat-C}
\widehat{C}_{j,k}^{(\ell-1)} 
= \dvert{C_{j,k}^{(\ell-1)}} + \dt^n \dsum_{s=1}^{4} L_{j,k,s}^{(\ell-1)} \dsum_{r=1}^{N_G}\eta_r\omega_{j,k,s,r}^{(\ell-1)}
= \dvert{C_{j,k}^{(\ell-1)}} + \dt^n \dsum_{s=1}^{4} L_{j,k,s}^{(\ell-1)} \omega_{j,k,s}^{(\ell-1)},
\end{equation} 
and $u_{j,k,s}^{(\ell-1), \rm H}$ is a sub-cell scheme using high-order numerical flux at cell interface:
\begin{equation}\label{EQ:uH}
u_{j,k,s}^{(\ell-1), \rm H} = \overline{u}_{j,k}^{(\ell-1)} - \lambda_{j,k}^{(\ell-1)}\sum_{r=1}^{N_G} \eta_r \left[\mbH \left(\w_{j,k,s,r}, \n_{j,k,s}^{(\ell-1)}, u_{j,k,s,r}^{(\ell-1),-},u_{j,k,s,r}^{(\ell-1),+}\right) - \mbH\left(\w_{j,k,s,r},\n_{j,k,s}^{(\ell-1)},\overline{u}_{j,k}^{(\ell-1)},\overline{u}_{j,k}^{(\ell-1)}\right)\right].
\end{equation} 
Here, the coefficient $\lambda_{j,k}^{(\ell-1)}$ is 
\begin{equation}\label{EQ:Lambda}
\lambda_{j,k}^{(\ell-1)} = \frac{\dt}{\widehat{C}_{j,k}^{(\ell-1)}}\dsum_{s=1}^{4}L_{j,k,s}^{(\ell-1)}.
\end{equation}

It is noteworthy that \eqref{EQ:Convex-Decomp} constitutes a valid convex combination of the cell average $\overline{u}_{j,k}^n$ and the high-order sub-cell schemes $u_{j,k,s}^{(\ell-1),\rm H}$, provided that $\widehat{C}_{j,k}^{(\ell-1)}> 0$. Combining Lemma \ref{Lem:Consistency-num-flux} and Lemma \ref{Lem:Sum-w}, we obtain 
\begin{equation*}
\begin{aligned}
\sum_{r=1}^{N_G} \eta_r \mbH\left(\w_{j,k,s,r},\n_{j,k,s}^{(\ell-1)},\overline{u}_{j,k}^{(\ell-1)},\overline{u}_{j,k}^{(\ell-1)}\right) &= \sum_{r=1}^{N_G} \eta_r \left[
\bH(\overline{u}_{j,k}^{(\ell-1)})\cdot\n_{j,k,s}^{(\ell-1)} - \left( \w_{j,k,s,r}\cdot\n_{j,k,s}^{(\ell-1)} \right)\overline{u}_{j,k}^{(\ell-1)}
\right]\\
&= \bH(\overline{u}_{j,k}^{(\ell-1)})\cdot\n_{j,k,s}^{(\ell-1)} -  \left[ \left(\sum_{r=1}^{N_G} \eta_r\w_{j,k,s,r}\right)\cdot\n_{j,k,s}^{(\ell-1)} \right]\overline{u}_{j,k}^{(\ell-1)} \\
&= \bH(\overline{u}_{j,k}^{(\ell-1)})\cdot\n_{j,k,s}^{(\ell-1)} - \left(\w_{j,k,s}\cdot\n_{j,k,s}^{(\ell-1)}\right)\overline{u}_{j,k}^{(\ell-1)} \\
&= \mbH\left(\w_{j,k,s},\n_{j,k,s}^{(\ell-1)},\overline{u}_{j,k}^{(\ell-1)},\overline{u}_{j,k}^{(\ell-1)}\right),
\end{aligned}
\end{equation*}
where $\w_{j,k,s}$ denotes the mesh velocity evaluated at the midpoint of the $s$-th edge of control volume $C_{j,k}$.
Accordingly, the high-order sub-cell scheme $u_{j,k,s}^{(\ell-1), \rm H}$ in \eqref{EQ:uH} can be further simplified as
\begin{equation}\label{EQ:uH-simplify}
u_{j,k,s}^{(\ell-1), \rm H} = \overline{u}_{j,k}^{(\ell-1)} - \lambda_{j,k}^{(\ell-1)}\left[
\sum_{r=1}^{N_G} \eta_r \mbH \left(\w_{j,k,s,r}, \n_{j,k,s}^{(\ell-1)}, u_{j,k,s,r}^{(\ell-1),-},u_{j,k,s,r}^{(\ell-1),+}\right) - \mbH\left(\w_{j,k,s},\n_{j,k,s}^{(\ell-1)},\overline{u}_{j,k}^{(\ell-1)},\overline{u}_{j,k}^{(\ell-1)}\right)\right].
\end{equation}

Analogously, the sub-cell scheme defined on the neighboring control volume $C_{\widehat{j},\widehat{k}}$ that shares the $s$-th edge takes the form
\begin{equation}\label{EQ:uH-NB-cell}
\begin{aligned}
u_{\widehat{j},\widehat{k},s}^{(\ell-1), \rm H} &= \overline{u}_{\widehat{j},\widehat{k}}^{(\ell-1)} - \lambda_{\widehat{j},\widehat{k}}^{(\ell-1)}\sum_{r=1}^{N_G} \eta_r \left[\mbH \left(\w_{\widehat{j},\widehat{k},s,r}, \n_{\widehat{j},\widehat{k},s}^{(\ell-1)}, u_{\widehat{j},\widehat{k},s,r}^{(\ell-1),-},u_{\widehat{j},\widehat{k},s,r}^{(\ell-1),+}\right) - \mbH\left(\w_{\widehat{j},\widehat{k},s,r},\n_{\widehat{j},\widehat{k},s}^{(\ell-1)},\overline{u}_{\widehat{j},\widehat{k}}^{(\ell-1)},\overline{u}_{\widehat{j},\widehat{k}}^{(\ell-1)}\right)\right]\\
&= \overline{u}_{\widehat{j},\widehat{k}}^{(\ell-1)} + \lambda_{\widehat{j},\widehat{k}}^{(\ell-1)}\sum_{r=1}^{N_G} \eta_r \left[\mbH \left(\w_{j,k,s,r}, \n_{j,k,s}^{(\ell-1)}, u_{j,k,s,r}^{(\ell-1),-},u_{j,k,s,r}^{(\ell-1),+}\right) - \mbH\left(\w_{j,k,s,r},\n_{j,k,s}^{(\ell-1)},\overline{u}_{\widehat{j},\widehat{k}}^{(\ell-1)},\overline{u}_{\widehat{j},\widehat{k}}^{(\ell-1)}\right)\right]\\
&= \overline{u}_{\widehat{j},\widehat{k}}^{(\ell-1)} + \lambda_{\widehat{j},\widehat{k}}^{(\ell-1)}\left[
\sum_{r=1}^{N_G} \eta_r \mbH \left(\w_{j,k,s,r}, \n_{j,k,s}^{(\ell-1)}, u_{j,k,s,r}^{(\ell-1),-},u_{j,k,s,r}^{(\ell-1),+}\right) - \mbH\left(\w_{j,k,s},\n_{j,k,s}^{(\ell-1)},\overline{u}_{\widehat{j},\widehat{k}}^{(\ell-1)},\overline{u}_{\widehat{j},\widehat{k}}^{(\ell-1)}\right)\right],
\end{aligned}
\end{equation}
where $\lambda_{\widehat{j},\widehat{k}}^{(\ell-1)}$ and $\widehat{C}_{\widehat{j},\widehat{k}}^{(\ell-1)}$ are, respectively, given by
\begin{equation}
\lambda_{\widehat{j},\widehat{k}}^{(\ell-1)} = \frac{\dt}{\widehat{C}_{\widehat{j},\widehat{k}}^{(\ell-1)}}\dsum_{s=1}^{4}L_{\widehat{j},\widehat{k},s}^{(\ell-1)}, \qquad
\widehat{C}_{\widehat{j},\widehat{k}}^{(\ell-1)} = \dvert{C_{\widehat{j},\widehat{k}}^{(\ell-1)}} + \dt^n \dsum_{s=1}^{4} L_{\widehat{j},\widehat{k},s}^{(\ell-1)} \omega_{\widehat{j},\widehat{k},s}^{(\ell-1)}.
\end{equation} 
In the second equality of \eqref{EQ:uH-NB-cell}, we have used \eqref{EQ:quant-adj-vol} and Lemma \ref{Lem:Inout-Flux}.

The BP property of the $\ell$-th stage of SSP-RK3 \eqref{EQ:RK-stage-l} is inherently determined by the BP property of every sub-cell scheme. Nevertheless, the high-order sub-cell scheme $u_{j,k,s}^{(\ell-1), \rm H}$ generally fails to have BP property, due to inevitable numerical errors. Following the numerical framework proposed in \cite{Gu2026Adaptive}, we enforce the BP property via a blending strategy that combines the high-order sub-cell scheme $u_{j,k,s}^{(\ell-1), \rm H}$ with its first-order counterpart $u_{j,k,s}^{(\ell-1), \rm L}$.
Different from the high-order sub-cell scheme, the first-order counterpart $u_{j,k,s}^{(\ell-1), \rm L}$ is constructed purely based on the cell averages $\overline{u}_{j,k}^{(\ell-1)}$ and $\overline{u}_{\widehat{j},\widehat{k}}^{(\ell-1)}$, rather than $u_{j,k,s,r}^{(\ell-1),\pm}$. The explicit formulation of $u_{j,k,s}^{(\ell-1), \rm L}$ is given by
\begin{equation}\label{EQ:Sub-cell-L}
u_{j,k,s}^{(\ell-1), \rm L} = \overline{u}_{j,k}^{(\ell-1)} - \lambda_{j,k}^{(\ell-1)}\left[
\sum_{r=1}^{N_G} \eta_r \mbH \left(\w_{j,k,s,r}, \n_{j,k,s}^{(\ell-1)}, \overline{u}_{j,k}^{(\ell-1)},\overline{u}_{\widehat{j},\widehat{k}}^{(\ell-1)}\right) - \mbH\left(\w_{j,k,s},\n_{j,k,s}^{(\ell-1)},\overline{u}_{j,k}^{(\ell-1)},\overline{u}_{j,k}^{(\ell-1)}\right)\right],
\end{equation} 
where $C_{\widehat{j},\widehat{k}}$ denotes the neighboring control volume that shares the $s$-th common edge with $C_{j,k}$. However, theoretical analysis of the BP property of $u_{j,k,s}^{(\ell-1), \rm L}$ is rather complicated, since the overall flux at the cell interface is a convex combination of fluxes evaluated at all Gauss quadrature points. To simplify the theoretical derivation, we adopt a simpler yet equivalent first-order sub-cell scheme with the following form:
\begin{equation}\label{EQ:Sub-cell-L-simplified}
u_{j,k,s}^{(\ell-1), \rm L} = \overline{u}_{j,k}^{(\ell-1)} - \lambda_{j,k}^{(\ell-1)} \left[\mbH \left(\w_{j,k,s}, \n_{j,k,s}^{(\ell-1)}, \overline{u}_{j,k}^{(\ell-1)},\overline{u}_{\widehat{j},\widehat{k}}^{(\ell-1)}\right) - \mbH\left(\w_{j,k,s},\n_{j,k,s}^{(\ell-1)},\overline{u}_{j,k}^{(\ell-1)},\overline{u}_{j,k}^{(\ell-1)}\right)\right].
\end{equation}
In the adjacent control volume $C_{\widehat{j},\widehat{k}}$, the simplified first-order sub-cell scheme corresponding to the $s$-th edge of $C_{j,k}$ takes the form:
\begin{equation}\label{EQ:Sub-cell-adj-L-simplified}
\begin{aligned}
u_{\widehat{j},\widehat{k},s}^{(\ell-1), \rm L} &= \overline{u}_{\widehat{j},\widehat{k}}^{(\ell-1)} - \lambda_{\widehat{j},\widehat{k}}^{(\ell-1)} \left[\mbH \left(\w_{\widehat{j},\widehat{k},s}, \n_{\widehat{j},\widehat{k},s}^{(\ell-1)}, \overline{u}_{\widehat{j},\widehat{k}}^{(\ell-1)},\overline{u}_{j,k}^{(\ell-1)}\right) - \mbH\left(\w_{\widehat{j},\widehat{k},s},\n_{\widehat{j},\widehat{k},s}^{(\ell-1)},\overline{u}_{\widehat{j},\widehat{k}}^{(\ell-1)},\overline{u}_{\widehat{j},\widehat{k}}^{(\ell-1)}\right)\right] \\
&= \overline{u}_{\widehat{j},\widehat{k}}^{(\ell-1)} +\lambda_{\widehat{j},\widehat{k}}^{(\ell-1)} \left[\mbH \left(\w_{j,k,s}, \n_{j,k,s}^{(\ell-1)}, \overline{u}_{j,k}^{(\ell-1)},\overline{u}_{\widehat{j},\widehat{k}}^{(\ell-1)}\right) - \mbH\left(\w_{j,k,s},\n_{j,k,s}^{(\ell-1)},\overline{u}_{\widehat{j},\widehat{k}}^{(\ell-1)},\overline{u}_{\widehat{j},\widehat{k}}^{(\ell-1)}\right)\right].
\end{aligned}
\end{equation}

We now rigorously investigate the BP property of the first-order sub-cell scheme \eqref{EQ:Sub-cell-L-simplified}. The main results are presented in Theorem \ref{Th:BP-Cond-l}.
\begin{theorem}\label{Th:BP-Cond-l}
Suppose $\overline{u}_{j,k}^{(\ell-1)}\in\mG$ holds for all $j$ and $k$. The sub-cell scheme $u_{j,k,s}^{(\ell-1), \rm L}$ computed by \eqref{EQ:Sub-cell-L-simplified} satisfies $u_{j,k,s}^{(\ell-1), \rm L}\in\mG$ for all $j$, $k$, and $s$, if the following two conditions are fulfilled: 
\begin{enumerate}[(1)]
\item All intermediate cell areas $\widehat{C}_{j,k}^{(\ell-1)}$ are positive, namely
\begin{equation}\label{EQ:BP-CFL-Scalar-area}
\widehat{C}_{j,k}^{(\ell-1)}= \dvert{C_{j,k}^{(\ell-1)}} + \dt^n \dsum_{s=1}^{4} L_{j,k,s}^{(\ell-1)} 
\omega_{j,k,s}^{(\ell-1)} > 0, \quad \forall\, j, \, k.
\end{equation}

\item The time step size $\dt^n$ obeys the following BP CFL constraint
\begin{equation}\label{EQ:BP-CFL-Scalar-dt}
\CFL:=\dt\max_{j,k,s}\left\lbrace \frac{1}{\widehat{C}_{j,k}^{(\ell-1)}}\left[\dsum_{q=1}^{4}L_{j,k,q}^{(\ell-1)}\right]\left(\omega_{j,k,s}^{(\ell-1)} - a_{j,k,s}^{(\ell-1),-}\right)\right\rbrace \leq 1,
\end{equation}
where local speeds of propagation are 
\begin{equation*}
\begin{aligned}
a_{j,k,s}^{(\ell-1),-} &= \min\left\lbrace\bH'(\overline{u}_{j,k}^{(\ell-1)})\cdot\n_{j,k,s}^{(\ell-1)},~\bH'(\overline{u}_{\widehat{j},\widehat{k}}^{(\ell-1)})\cdot\n_{j,k,s}^{(\ell-1)},~ -\varepsilon^{\rm CU} \right\rbrace,\\
a_{j,k,s}^{(\ell-1),+} &= \max\left\lbrace\bH'(\overline{u}_{j,k}^{(\ell-1)})\cdot\n_{j,k,s}^{(\ell-1)},~\bH'(\overline{u}_{\widehat{j},\widehat{k}}^{(\ell-1)})\cdot\n_{j,k,s}^{(\ell-1)},~\varepsilon^{\rm CU} \right\rbrace.
\end{aligned}
\end{equation*}
\end{enumerate}
\end{theorem}
\begin{proof}
The first-order sub-cell scheme $u_{j,k,s}^{(\ell-1), \rm L}$ can be rewritten in the equivalent form
\begin{equation}
u_{j,k,s}^{(\ell-1), \rm L} = \overline{u}_{j,k}^{(\ell-1)} - \lambda_{j,k}^{(\ell-1)} \left[\mbH \left(\w_{j,k,s}, \n_{j,k,s}^{(\ell-1)}, \overline{u}_{j,k}^{(\ell-1)},\overline{u}_{\widehat{j},\widehat{k}}^{(\ell-1)}\right) -
\left(\bH(\overline{u}_{j,k}^{(\ell-1)})\cdot\n_{j,k,s}^{(\ell-1)} - \omega_{j,k,s}^{(\ell-1)}\,\overline{u}_{j,k}^{(\ell-1)}
\right)\right].
\end{equation}
where the coefficient $\lambda_{j,k}^{(\ell-1)}$ is defined in \eqref{EQ:Lambda}. It is clear that $\lambda_{j,k}^{(\ell-1)}>0$ provided the cell area satisfies $\widehat{C}_{j,k}^{(\ell-1)}> 0$.

The BP property is verified via a comprehensive three-case analysis as follows.

\noindent\textbf{Case I}: $\w_{j,k,s}\cdot\n_{j,k,s}^{(\ell-1)}\leq a_{j,k,s}^{(\ell-1),-}$. In this case, the sub-cell scheme degenerates to $u_{j,k,s}^{(\ell-1), \rm L} = \overline{u}_{j,k}^{(\ell-1)}$. Therefore, the BP property $u_{j,k,s}^{(\ell-1), \rm L}\in\mG$ holds trivially from the given assumption.

\noindent\textbf{Case II}: $\w_{j,k,s}\cdot\n_{j,k,s}^{(\ell-1)}\geq a_{j,k,s}^{(\ell-1),+}$. In this case, the sub-cell scheme admits the linear decomposition $u_{j,k,s}^{(\ell-1), \rm L} = \Phi_1(\overline{u}_{j,k}^{(\ell-1)}) + \Phi_2(\overline{u}_{\widehat{j},\widehat{k}}^{(\ell-1)})$,
with the two functions given by 
\begin{equation}
\begin{aligned}
\Phi_1(\overline{u}_{j,k}^{(\ell-1)}) &= \left[ 1-\lambda_{j,k}^{(\ell-1)}(\w_{j,k,s}\cdot\n_{j,k,s}^{(\ell-1)}) \right]\overline{u}_{j,k}^{(\ell-1)} + \lambda_{j,k}^{(\ell-1)}\bH(\overline{u}_{j,k}^{(\ell-1)})\cdot\n_{j,k,s}^{(\ell-1)},\\
\Phi_2(\overline{u}_{\widehat{j},\widehat{k}}^{(\ell-1)}) &= \lambda_{j,k}^{(\ell-1)}(\w_{j,k,s}\cdot\n_{j,k,s}^{(\ell-1)}) \overline{u}_{\widehat{j},\widehat{k}}^{(\ell-1)} - \lambda_{j,k}^{(\ell-1)}\bH(\overline{u}_{\widehat{j},\widehat{k}}^{(\ell-1)})\cdot\n_{j,k,s}^{(\ell-1)}.
\end{aligned}
\end{equation}
The first-order derivatives of $\Phi_1$ and $\Phi_2$ take the form
\begin{equation}
\begin{aligned}
\Phi_1'(\overline{u}_{j,k}^{(\ell-1)}) &= 1-\lambda_{j,k}^{(\ell-1)}\left[ \w_{j,k,s}\cdot\n_{j,k,s}^{(\ell-1)} - \bH'(\overline{u}_{j,k}^{(\ell-1)})\cdot\n_{j,k,s}^{(\ell-1)} \right]\geq 1-\lambda_{j,k}^{(\ell-1)}\left[ \w_{j,k,s}\cdot\n_{j,k,s}^{(\ell-1)} - a_{j,k,s}^- \right],\\
\Phi_2'(\overline{u}_{\widehat{j},\widehat{k}}^{(\ell-1)}) &=\lambda_{j,k}^{(\ell-1)}\left[ \w_{j,k,s}\cdot\n_{j,k,s}^{(\ell-1)} - \bH'(\overline{u}_{\widehat{j},\widehat{k}}^{(\ell-1)})\cdot\n_{j,k,s}^{(\ell-1)} \right].
\end{aligned}
\end{equation}
Under the conditions \eqref{EQ:BP-CFL-Scalar-area} and \eqref{EQ:BP-CFL-Scalar-dt}, the derivatives are both non-negative, namely $\Phi_1'\geq 0$ and $\Phi_2'\geq 0$, indicating that both functions are monotonically increasing. Consequently, the BP property is guaranteed by the monotonicity of the sub-cell scheme:
\begin{equation}\label{EQ:Mono-Ineq}
u_m = \Phi_1(u_m) + \Phi_2(u_m) \leq \Phi_1(\overline{u}_{j,k}^{(\ell-1)}) + \Phi_2(\overline{u}_{\widehat{j},\widehat{k}}^{(\ell-1)}) \leq \Phi_1(u_M) + \Phi_2(u_M) = u_M,
\end{equation}
there $u_m$ and $u_M$ are the minimum and maximum values to be preserved, respectively.

\noindent\textbf{Case III}: $a_{j,k,s}^{(\ell-1),-} < \w_{j,k,s}\cdot\n_{j,k,s}^{(\ell-1)}< a_{j,k,s}^{(\ell-1),+}$. In this case, we decompose the sub-cell scheme into $u_{j,k,s}^{(\ell-1), \rm L} = \Psi_1(\overline{u}_{j,k}^{(\ell-1)}) + \Psi_2(\overline{u}_{\widehat{j},\widehat{k}}^{(\ell-1)})$, where functions $\Psi_1$ and $\Psi_2$ are
\begin{equation*}
\begin{aligned}
\Psi_1(\overline{u}_{j,k}^{(\ell-1)}) &= \left[ 1-\lambda_{j,k}^{(\ell-1)}a_{j,k,s}^{(\ell-1),+}\frac{\w_{j,k,s}\cdot\n_{j,k,s}^{(\ell-1)} - a_{j,k,s}^{(\ell-1),-}}{a_{j,k,s}^{(\ell-1),+} - a_{j,k,s}^{(\ell-1),-}} \right]\overline{u}_{j,k}^{(\ell-1)} + \lambda_{j,k}^{(\ell-1)}\frac{\w_{j,k,s}\cdot\n_{j,k,s}^{(\ell-1)} - a_{j,k,s}^{(\ell-1),-}}{a_{j,k,s}^{(\ell-1),+} - a_{j,k,s}^{(\ell-1),-}}\bH(\overline{u}_{j,k}^{(\ell-1)})\cdot\n_{j,k,s}^{(\ell-1)},\\
\Psi_2(\overline{u}_{\widehat{j},\widehat{k}}^{(\ell-1)}) &= \lambda_{j,k}^{(\ell-1)}a_{j,k,s}^{(\ell-1),+}\frac{\w_{j,k,s}\cdot\n_{j,k,s}^{(\ell-1)} - a_{j,k,s}^{(\ell-1),-}}{a_{j,k,s}^{(\ell-1),+} - a_{j,k,s}^{(\ell-1),-}} \overline{u}_{\widehat{j},\widehat{k}}^{(\ell-1)} - \lambda_{j,k}^{(\ell-1)}\frac{\w_{j,k,s}\cdot\n_{j,k,s}^{(\ell-1)} - a_{j,k,s}^{(\ell-1),-}}{a_{j,k,s}^{(\ell-1),+} - a_{j,k,s}^{(\ell-1),-}}\bH(\overline{u}_{\widehat{j},\widehat{k}}^{(\ell-1)})\cdot\n_{j,k,s}^{(\ell-1)}.
\end{aligned}
\end{equation*}
The first-order derivatives are
\begin{equation*}
\begin{aligned}
\Psi_1'(\overline{u}_{j,k}^{(\ell-1)}) &= 1-\lambda_{j,k}^{(\ell-1)}\left[a_{j,k,s}^{(\ell-1),+} - \bH'(\overline{u}_{j,k}^{(\ell-1)})\cdot\n_{j,k,s}^{(\ell-1)} \right] \frac{\w_{j,k,s}\cdot\n_{j,k,s}^{(\ell-1)} - a_{j,k,s}^{(\ell-1),-}}{a_{j,k,s}^{(\ell-1),+} - a_{j,k,s}^{(\ell-1),-}} \\
&\geq 1-\lambda_{j,k}^{(\ell-1)}\left(\w_{j,k,s}\cdot\n_{j,k,s}^{(\ell-1)} - a_{j,k,s}^{(\ell-1),-} \right),\\
\Psi_2'(\overline{u}_{\widehat{j},\widehat{k}}^{(\ell-1)}) &= \lambda_{j,k}^{(\ell-1)}\left[a_{j,k,s}^{(\ell-1),+} - \bH'(\overline{u}_{\widehat{j},\widehat{k}}^{(\ell-1)})\cdot\n_{j,k,s}^{(\ell-1)} \right]\frac{\w_{j,k,s}\cdot\n_{j,k,s}^{(\ell-1)} - a_{j,k,s}^{(\ell-1),-}}{a_{j,k,s}^{(\ell-1),+} - a_{j,k,s}^{(\ell-1),-}}. 
%\overline{u}_{\widehat{j},\widehat{k}}^{(\ell-1)}.
\end{aligned}
\end{equation*}
Subject to constraints \eqref{EQ:BP-CFL-Scalar-area} and \eqref{EQ:BP-CFL-Scalar-dt}, the non-negativity properties $\Psi_1'\geq 0$ and $\Psi_2'\geq 0$ are guaranteed. Following the same monotonicity argument as in \eqref{EQ:Mono-Ineq}, we conclude that $u_m\leq u_{j,k,s}^{(\ell-1), \rm L} \leq u_M$.

Combining the theoretical results from all three cases, we  complete the proof of Theorem \ref{Th:BP-Cond-l}.
\end{proof}

With the help of $u_{j,k,s}^{(\ell-1), \rm L}\in\mG$ and $u_{\widehat{j},\widehat{k},s}^{(\ell-1), \rm L}\in\mG$, where $u_{j,k,s}^{(\ell-1), \rm L}$ and $u_{\widehat{j},\widehat{k},s}^{(\ell-1), \rm L}$ are defined in \eqref{EQ:Sub-cell-L-simplified} and \eqref{EQ:Sub-cell-adj-L-simplified}, respectively, we perform BP modifications on the high-order sub-cell schemes $u_{j,k,s}^{(\ell-1), \rm H}$ and $u_{\widehat{j},\widehat{k},s}^{(\ell-1), \rm H}$. 
Specifically, we blend the high-order sub-cell schemes $u_{j,k,s}^{(\ell-1), \rm H}$ and $u_{\widehat{j},\widehat{k},s}^{(\ell-1), \rm H}$ with their BP first-order counterparts to construct
\begin{equation}\label{EQ:blend-RK}
u_{j,k,s}^{(\ell-1), \rm H, (\mathrm{mod})} = \theta_s^{(\ell-1)} u_{j,k,s}^{(\ell-1), \rm H} + (1-\theta_s^{(\ell-1)})u_{j,k,s}^{(\ell-1), \rm L}, \qquad
u_{\widehat{j},\widehat{k},s}^{(\ell-1), \rm H, (\mathrm{mod})} = \theta_s^{(\ell-1)} u_{\widehat{j},\widehat{k},s}^{(\ell-1), \rm H} + (1-\theta_s^{(\ell-1)})u_{\widehat{j},\widehat{k},s}^{(\ell-1), \rm L}.
\end{equation}
The original high-order sub-cell schemes $u_{j,k,s}^{(\ell-1), \rm H}$ and $u_{\widehat{j},\widehat{k},s}^{(\ell-1), \rm H}$ are then replaced by the corrected ones $u_{j,k,s}^{(\ell-1), \rm H, (\mathrm{mod})}$ and $u_{\widehat{j},\widehat{k},s}^{(\ell-1), \rm H, (\mathrm{mod})}$, respectively. Here, the common blending coefficient is defined as $\theta_s^{(\ell-1)} = \min\{ \theta_{j,k,s}^{(\ell-1)}, \theta_{\widehat{j},\widehat{k},s}^{(\ell-1)} \}$, with
\begin{equation}\label{EQ:theta}
\begin{aligned}
\theta_{j,k,s}^{(\ell-1)} &= \min\left\lbrace \Theta\left(u_{j,k,s}^{(\ell-1), \rm H}-u_m, u_{j,k,s}^{(\ell-1), \rm L}-u_m, \varepsilon^\mathrm{BP}\right),\, \Theta\left(u_M-u_{j,k,s}^{(\ell-1), \rm H}, u_M-u_{j,k,s}^{(\ell-1), \rm L}, \varepsilon^\mathrm{BP}\right)\right\rbrace, \vspace{+3pt}\\
\theta_{\widehat{j},\widehat{k},s}^{(\ell-1)} &= \min\left\lbrace \Theta\left(u_{\widehat{j},\widehat{k},s}^{(\ell-1), \rm H}-u_m, u_{\widehat{j},\widehat{k},s}^{(\ell-1), \rm L}-u_m, \varepsilon^\mathrm{BP}\right),\, \Theta\left(u_M-u_{\widehat{j},\widehat{k},s}^{(\ell-1), \rm H}, u_M-u_{\widehat{j},\widehat{k},s}^{(\ell-1), \rm L}, \varepsilon^\mathrm{BP}\right)\right\rbrace,
\end{aligned}
\end{equation}
where the limiter function $\Theta$ is defined by
\begin{equation}
\Theta\left(u^{\rm H}, u^{\rm L}, \varepsilon^\mathrm{BP}\right) = \begin{cases}
1, & \text{if } u^{\rm H} \geq \varepsilon^\mathrm{BP}, \vspace{+3pt}\\
\dfrac{u^{\rm L}-\varepsilon^\mathrm{BP}}{u^{\rm L}-u^{\rm H}}, & \text{if } u^{\rm H} < \varepsilon^\mathrm{BP} ~\&~ u^{\rm L} \geq \varepsilon^\mathrm{BP}, \vspace{+3pt}\\
0, & \text{if } u^{\rm H} < \varepsilon^\mathrm{BP} ~\&~ u^{\rm L} < \varepsilon^\mathrm{BP},
\end{cases}
\end{equation}
A small positive parameter $\varepsilon^\mathrm{BP}$ is introduced to eliminate potential round-off errors, and is set to $\varepsilon^\mathrm{BP}=10^{-13}$ in all numerical tests. Under the valid BP constraints \eqref{EQ:BP-CFL-Scalar-area} and \eqref{EQ:BP-CFL-Scalar-dt}, the proposed correction strategy ensures that the modified high-order sub-cell schemes satisfy $u_{j,k,s}^{(\ell-1), \rm H, (\mathrm{mod})}\in\mG$ for all indices $j$, $k$, and $s$.
It is worth noting that the present BP correction essentially introduces a flux-limiting procedure. This property can be clarified by reformulating the corrected sub-cell schemes $u_{j,k,s}^{(\ell-1), \rm H, (\mathrm{mod})}$ and $u_{\widehat{j},\widehat{k},s}^{(\ell-1), \rm H, (\mathrm{mod})}$ as
\begin{equation}
\begin{aligned}
u_{j,k,s}^{(\ell-1), \rm H, (\mathrm{mod})} = \theta_s^{(\ell-1)} u_{j,k,s}^{(\ell-1), \rm H} + (1-\theta_s^{(\ell-1)})u_{j,k,s}^{(\ell-1), \rm L}
&= \overline{u}_{j,k}^{(\ell-1)} - \lambda_{j,k}^{(\ell-1)}\left[ H_{j,k,s}^{(\mathrm{mod})} - \mbH\left(\w_{j,k,s},\n_{j,k,s}^{(\ell-1)},\overline{u}_{j,k}^{(\ell-1)},\overline{u}_{j,k}^{(\ell-1)}\right)\right],\\
u_{\widehat{j},\widehat{k},s}^{(\ell-1), \rm H, (\mathrm{mod})} = \theta_s^{(\ell-1)} u_{\widehat{j},\widehat{k},s}^{(\ell-1), \rm H} + (1-\theta_s^{(\ell-1)})u_{\widehat{j},\widehat{k},s}^{(\ell-1), \rm L}
&= \overline{u}_{\widehat{j},\widehat{k}}^{(\ell-1)} - \lambda_{\widehat{j},\widehat{k}}^{(\ell-1)}\left[ H_{\widehat{j},\widehat{k},s}^{(\mathrm{mod})} - \mbH\left(\w_{\widehat{j},\widehat{k},s},\n_{\widehat{j},\widehat{k},s}^{(\ell-1)},\overline{u}_{\widehat{j},\widehat{k}}^{(\ell-1)},\overline{u}_{\widehat{j},\widehat{k}}^{(\ell-1)}\right)\right]\\
&=\overline{u}_{\widehat{j},\widehat{k}}^{(\ell-1)} + \lambda_{\widehat{j},\widehat{k}}^{(\ell-1)}\left[ H_{j,k,s}^{(\mathrm{mod})} - \mbH\left(\w_{j,k,s},\n_{j,k,s}^{(\ell-1)},\overline{u}_{\widehat{j},\widehat{k}}^{(\ell-1)},\overline{u}_{\widehat{j},\widehat{k}}^{(\ell-1)}\right)\right],
\end{aligned}
\end{equation} 
where $H_{j,k,s}^{(\mathrm{mod})}$ and $H_{\widehat{j},\widehat{k},s}^{(\mathrm{mod})}$ are the modified numerical flux, satisfying
\begin{equation}
\begin{aligned}
H_{j,k,s}^{(\mathrm{mod})} &= -H_{\widehat{j},\widehat{k},s}^{(\mathrm{mod})} \\
&= \theta_s^{(\ell-1)} \sum_{r=1}^{N_G} \eta_r \mbH \left(\w_{j,k,s,r}, \n_{j,k,s}^{(\ell-1)}, u_{j,k,s,r}^{(\ell-1),-},u_{j,k,s,r}^{(\ell-1),+}\right) + (1-\theta_s^{(\ell-1)}) \mbH \left(\w_{j,k,s}, \n_{j,k,s}^{(\ell-1)}, \overline{u}_{j,k}^{(\ell-1)},\overline{u}_{\widehat{j},\widehat{k}}^{(\ell-1)}\right).
\end{aligned}
\end{equation}
This formulation demonstrates that a common modified flux is adopted at each cell interface, which is formulated as a convex combination of the original high-order flux with its first-order counterpart. Therefore, the BP correction inherently introduces a flux limiter while strictly preserving the discrete conservation of the original numerical scheme.

Next, we verify in Theorem \ref{Th:GCL} that the D-GCL is preserved by the BP corrections. 
\begin{theorem}\label{Th:GCL}
The proposed BP correction procedure preserves the discrete geometric conservation law.
\end{theorem}
\begin{proof}
Consider a uniform flow $\overline{u}_{j,k}^{(\ell-1)}=u^0=\mathrm{const}$ for all $j$ and $k$. In this case, the WENO reconstruction yields $u_{j,k,s,r}^{(\ell-1),\pm} = u^0$. Accordingly, the modified BP sub-cell scheme $u_{j,k,s}^{(\ell-1),\rm H, (\mathrm{mod})}$ simplifies to 
\begin{equation}
\begin{aligned}
u_{j,k,s}^{(\ell-1),\rm H, (\mathrm{mod})} =& \theta_s^{(\ell-1)} u_{j,k,s}^{(\ell-1), \rm H} + (1-\theta_s^{(\ell-1)})u_{j,k,s}^{(\ell-1), \rm L}\\
=&\theta_s^{(\ell-1)} \left\lbrace u^0 - \lambda_{j,k}^{(\ell-1)}\sum_{r=1}^{N_G} \eta_r \left[\mbH \left(\w_{j,k,s,r}, \n_{j,k,s}^{(\ell-1)}, u^0,u^0\right) - \mbH\left(\w_{j,k,s,r},\n_{j,k,s}^{(\ell-1)},u^0,u^0\right)\right] \right\rbrace\\
&+ (1-\theta_s^{(\ell-1)})\left\lbrace u^0 - \lambda_{j,k}^{(\ell-1)} \left[\mbH \left(\w_{j,k,s}, \n_{j,k,s}^{(\ell-1)}, u^0,u^0\right) - \mbH\left(\w_{j,k,s},\n_{j,k,s}^{(\ell-1)},u^0,u^0\right)\right]\right\rbrace\\
=& u^0.
\end{aligned}
\end{equation}
Therefore, we have
\begin{equation}
\begin{aligned}
\overline{u}_{j,k}^{(\ell)} &= \xi^\ell\chfrac{\dvert{C_{j,k}^n}}{\dvert{C_{j,k}^{\ell}}}\overline{u}_{j,k}^{n} + (1-\xi^\ell)\chfrac{\widehat{C}_{j,k}^{(\ell-1)}}{\dvert{C_{j,k}^{\ell}}}\sum_{s=1}^{4}\frac{L_{j,k,s}^{(\ell-1)}}{\sum_{q=1}^{4}L_{j,k,q}^{(\ell-1)}} \widetilde{u}_{j,k,s}^{(\ell-1),\rm H} \\
&= \chfrac{1}{\dvert{C_{j,k}^{\ell}}}\left[
\xi^\ell \dvert{C_{j,k}^n} + (1-\xi^\ell)\widehat{C}_{j,k}^{(\ell-1)}
\right] u^0.
\end{aligned}
\end{equation}
It immediately follows from \eqref{EQ:RK-stage-Cjk-l} and \eqref{EQ:hat-C} that $\overline{u}_{j,k}^{(\ell)} = u^0$. 
By mathematical induction, the evolved cell average satisfies $\overline{u}_{j,k}^{n+1} = u^0$. This indicates that the uniform flow remains invariant under the proposed numerical scheme, which verifies the preservation of the D-GCL. This completes the proof of Theorem \ref{Th:GCL}.
\end{proof}

Taking into account the three stages of the SSP-RK3 method, we finally establish the BP constraints in the following theorem. 
\begin{theorem}\label{Th:BP-Cond}
Assume that $\overline{u}_{j,k}^n\in\mG$ for all $j$ and $k$. The fully discrete numerical scheme produces $\overline{u}_{j,k}^{n+1}\in\mG$ for all grid indices, provided that the constraints \eqref{EQ:BP-CFL-Scalar-area} and \eqref{EQ:BP-CFL-Scalar-dt} hold at every SSP-RK3 stage.
\end{theorem}
\begin{proof}
For each $\ell$-th stage of the SSP-RK3, the validity of the constraints \eqref{EQ:BP-CFL-Scalar-area} and \eqref{EQ:BP-CFL-Scalar-dt} guarantees the BP property of the first-order sub-cell schemes according to Theorem \ref{Th:BP-Cond-l}. The convex blending operation \eqref{EQ:blend-RK} further produces BP sub-cell schemes, namely $u_{j,k,s}^{(\ell-1), \rm H, (\mathrm{mod})}\in\mG$. Under the condition \eqref{EQ:BP-CFL-Scalar-area}, $\overline{u}_{j,k}^{(\ell)}$ is a convex combination of BP quantities and thus remains within the admissible set $\mG$. Since the above argument holds for arbitrary SSP-RK3 stages, the evolved cell average $\overline{u}_{j,k}^{n+1}$ satisfies BP property. This completes the proof of Theorem \ref{Th:BP-Cond}.
\end{proof}

\subsection{Adaptive moving mesh redistribution}\label{Sec:mesh-redist}

The mesh velocity involved in the discrete formulations \eqref{EQ:AMM-FVM-U} and \eqref{EQ:AMM-FVM-C} is determined by the moving mesh equation. This section briefly reviews the adaptive mesh redistribution algorithm developed in \cite{Kurganov2021Adaptive,Tang2003Adaptive}, which is adopted in the present numerical framework.

Let $\xi$ and $\eta$ denote the logical coordinates defined on the computational domain $\Omega_c = [0,1]^2$, which is partitioned by a uniform square mesh with identical number of cells to the physical domain. A one-to-one continuous mapping is constructed to bridge the physical coordinates $\z_{j,k}=(x_j,y_k)$ and the logical coordinates $\widehat{\z}_{j,k}=(\xi_j,\eta_k)$. The distribution of the adaptive mesh points in the physical domain is determined by minimizing the following functional:
\begin{equation}\label{EQ:functional}
J(\z) = \half\sum_{i=1}^2 \int_{\Omega_c} (\widetilde{\nabla}\z)^\top G_i \widetilde{\nabla}\z,
\end{equation}
where $\widetilde{\nabla} = (\partial_\xi, \partial_\eta)^\top$ represents the gradient operator with respect to the logical variables, and $G_1$ and $G_2$ are two symmetric positive definite matrices. Following the standard setting in \cite{Tang2003Adaptive}, we set $G_1=G_2=\sigma\mathrm{I}$, where $\mathrm{I}$ denotes the $2\times 2$ identity matrix and $\sigma$ denotes the monitor function that governs mesh adaption. Instead of directly solving the minimization problem, we compute the adaptive mesh points by solving the corresponding Euler-Lagrange equation \cite{Tang2003Adaptive}
\begin{equation*}
(\sigma\z_\xi)_\xi + (\sigma\z_\eta)_\eta = 0.
\end{equation*}
Following \cite{Kurganov2021Adaptive}, the monitor function $\sigma$ is defined by
\begin{equation*}
\sigma = 1+ \alpha\varphi(\dVert{D\U}),
\end{equation*}
where $\dVert{D\U}_{j,k}$ measures the magnitude of the solution gradient $\dvert{\nabla\U^{(i)}}_{j,k}$ for some specific components of $\U$, $\varphi$ denotes the smoothing function, and $\alpha$ is an intensity parameter that controls concentration of mesh points. In addition, the smoothing function $\varphi$ takes the form
\begin{equation*}
\begin{aligned}
\varphi_{j,k}^{l+1} = &\dfrac{1}{4}\varphi_{j,k}^{l} + \dfrac{1}{8}\left(
\varphi_{j-1,k}^{l} + \varphi_{j+1,k}^{l} + \varphi_{j,k-1}^{l} + \varphi_{j,k+1}^{l}
\right) \\
& + \dfrac{1}{16}\left(
\varphi_{j-1,k-1}^{l} + \varphi_{j-1,k+1}^{l} + \varphi_{j+1,k-1}^{l} + \varphi_{j+1,k+1}^{l}
\right), \quad
l = 0, 1, \dots, m-1.
\end{aligned}
\end{equation*}
In this paper, eight smoothing iterations will be taken, namely $m=8$.
Equipped with $\varphi_{j,k}^{m}$, we calculate the intensity parameter $\alpha$ by
\begin{equation*}
\alpha = \left(
\dfrac{1-\beta}{\beta\dvert{\Omega}}\iint_\Omega  \varphi(\dVert{D\U})
\right)^{-1} 
=\left(
\dfrac{1-\beta}{\beta\dvert{\Omega}}\sum_{j,k} \varphi_{j,k}^m \dvert{C_{j,k}}
\right)^{-1}, 
\end{equation*}
with $\beta\in (0,1)$.

With the monitor function $\sigma_{j,k}$ precomputed, the linearized Euler-Lagrange equation is solved iteratively using Jacobi iteration, that is
\begin{equation*}
\begin{aligned}
&\dfrac{\sigma_{j+1,k+\half}(\z_{j+\frac{3}{2},k+\half}^{[\nu]}-\z_{j+\half,k+\half}^{[\nu+1]})-\sigma_{j,k+\half}(\z_{j+\half,k+\half}^{[\nu+1]}-\z_{j-\half,k+\half}^{[\nu]})}{\Delta \xi^2} \\
+& 
\dfrac{\sigma_{j+\half,k+1}(\z_{j+\half,k+\frac{3}{2}}^{[\nu]}-\z_{j+\half,k+\half}^{[\nu+1]})-\sigma_{j+\half,k}(\z_{j+\half,k+\half}^{[\nu+1]}-\z_{j+\half,k-\half}^{[\nu]})}{\Delta \eta^2}
=0, \qquad
\nu = 0, 1,\dots,s-1,
\end{aligned}
\end{equation*}
where 
\begin{equation*}
\sigma_{j,k+\half} = \half(\sigma_{j,k} + \sigma_{j,k+1}), \quad
\sigma_{j+\half,k} = \half(\sigma_{j,k} + \sigma_{j+1,k}).
\end{equation*}
During the Jacobi iteration, mesh points should be adjusted subject to the constraint
\begin{equation*}
\z_{\jph,\kph}^{\nu+1}\in [\z_{j,\kph}^\nu,\z_{j+1,\kph}^\nu]\times [\z_{\jph,k}^\nu,\z_{\jph,k+1}^\nu]\cap [\z_{j,\kph}^n,\z_{j+1,\kph}^n]\times [\z_{\jph,k}^n,\z_{\jph,k+1}^n].
\end{equation*}
This constraint strictly preserves the mesh topological structure throughout the procedure of mesh adaption. In all numerical examples, we take four Jacobi iterations, namely $s=4$.

\begin{remark}
To avoid excessively small cell sizes and improve numerical robustness, additional mesh regularization strategies can also be incorporated. We refer to \cite{Kurganov2022Well} and the references therein for more details.
\end{remark}

The updated location of mesh point $\z_{\jph,\kph}^{[s]}$ obtained from the Jacobi iteration serves as the candidate coordinates of $\z_{\jph,\kph}$ at the next time level $t^{n+1}$. Accordingly, the time-dependent mesh velocity is explicitly defined by
\begin{equation*}
\w_{\jph,\kph} = \frac{\z_{\jph,\kph}^{[s]} - \z_{\jph,\kph}^{[0]}}{\dt^n},
\end{equation*}
where $\z_{\jph,\kph}^{[0]}=\z_{\jph,\kph}^n$, and $t^n$ is computed by \eqref{EQ:Conventional-dt}.

To close this section, we provide several practical implementation details for the proposed BP AMM-FVM.
The BP constraints \eqref{EQ:BP-CFL-Scalar-area} and \eqref{EQ:BP-CFL-Scalar-dt} indicate that the admissible time step size $\dt^n$ depends implicitly on both the instantaneous numerical solution and the cell areas during the entire SSP-RK3 evolution.
In practical computations, we initialize the time step size following the conventional CFL condition \eqref{EQ:Conventional-dt} and solve the moving mesh equation to obtain the mesh velocity $\w_{\jph,\kph}$, which is frozen over the time interval $[t^n, t^{n+1}]$. If either of the following two conditions occurs,
\begin{enumerate}[(1)]
\item the intermediate cell area $\widehat{C}_{j,k}^{(\ell-1)}$ becomes non-positive for some $j$, $k$, and $\ell$,

\item The constraint $\widetilde{u}_{j,k,s}^{(\ell-1),\mathrm{H}}\in\mG$ fails to be satisfied by a convex combination of $u_{j,k,s}^{(\ell-1),\mathrm{H}}$ and $u_{j,k,s}^{(\ell-1),\mathrm{L}}$ for some $j$, $k$, and $\ell$, 
\end{enumerate}
then the current time step is repeated with a halved time step size. Accordingly, 
whenever the SSP-RK3 procedure is restarted, the updated locations of mesh points $\z_{\jph,\kph}^{n+1}$ are consistently recalculated using the frozen mesh velocity $\w_{\jph,\kph}$ and the reduced time step size.
This time-marching procedure, if triggered, terminates after a finite number of repetitions, since both BP constraints \eqref{EQ:BP-CFL-Scalar-area} and \eqref{EQ:BP-CFL-Scalar-dt} are satisfied as $\dt^n\to 0$.

\section{Application to Euler equations}\label{Sec:Euler}
This section extends the proposed third-order BP AMM-FVM to the two-dimensional Euler equations of gas dynamics, which read as
\begin{equation*}
\begin{array}{l}
\U_t + \nabla\cdot\bH(\U) = 0,\quad
\bH(\U) = (\F(\U),\,\G(\U)), \quad
\U = \begin{bmatrix}
\rho\\
\rho u\\
\rho v\\
E
\end{bmatrix},\quad
\F(\U) = \begin{bmatrix}
\rho u\\
\rho u^2 + p\\
\rho uv\\
(E+p)u
\end{bmatrix},\quad
\G(\U) = \begin{bmatrix}
\rho v\\
\rho uv\\
\rho v^2 + p\\
(E+p)v
\end{bmatrix},
\end{array}
\end{equation*}
where $\rho$ denotes the fluid density, $u$ and $v$ represent the velocity components in the $x$- and $y$-directions, respectively, $p$ is the thermodynamic pressure, and $E$ denotes the total energy. The system is closed by the ideal gas equation of state (EOS)
\begin{equation*}
p = (\gamma-1)\left[E-\dfrac{1}{2}\rho (u^2+ v^2)\right], 
\end{equation*}
where $\gamma$ is the specific heat ratio, which is fixed to $1.4$ for all numerical examples in this paper. The corresponding speed of sound is defined as $c = \sqrt{\nicefrac{\gamma p}{\rho}}$. To make the quantities in the square root positive, we define the physics-constraint invariant domain $\mG$ as
\begin{equation*}
\mG = \left\lbrace
\U\vert ~ \rho[\U]>0,\, p[\U]>0
\right\rbrace.
\end{equation*}
One can prove by Jensen's inequality that the set $\mG$ defined above is convex.

The adaptive moving mesh finite-volume method combined with the SSP-RK3 temporal discretization, developed for scalar hyperbolic equations, can be directly extended to the Euler system following the discrete formulations \eqref{EQ:AMM-FVM-U} and \eqref{EQ:AMM-FVM-C}. The major modification lies in the replacement of the scalar solution with the vector-valued conservative variables $\U_{j,k}$. Besides, the CU flux operator for the Euler system takes the form
\begin{equation*}
\mbH( \w, \n, \U^-, \U^+ ) = \begin{cases}
\bH(\U^-)\cdot\n - \omega\U^-, &\omega < a^-, \vspace{+3pt}\\
\dfrac{a^+ \bH(\U^-)\cdot\n - a^- \bH(\U^+)\cdot\n}{a^+ - a^-} + \dfrac{a^- a^+}{a^+ - a^-}(\U^+ - \U^-) - \omega \U^*, &a^-\leq \omega \leq a^+, \vspace{+3pt}\\
\bH(\U^+)\cdot\n - \omega\U^+, &\omega > a^+,
\end{cases}
\end{equation*}
where $\omega=\w\cdot\n$, and
\begin{equation*}
\begin{aligned}
a^- &= \min\left\lbrace (u^-, v^-)\cdot\n - c^-,~(u^+, v^+)\cdot\n - c^+,~ -\varepsilon^{\rm CU} \right\rbrace, \vspace{+3pt}\\
a^+ &= \max\left\lbrace (u^-, v^-)\cdot\n + c^-,~(u^+, v^+)\cdot\n + c^+,~\varepsilon^{\rm CU} \right\rbrace, \vspace{+3pt}\\
\U^* &= \dfrac{a^+ \U^+ - a^- \U^- - \bH(\U^+)\cdot\n + \bH(\U^-)\cdot\n}{a^+ - a^-}.
\end{aligned}
\end{equation*}

Consistent with the scalar framework, each SSP-RK3 stage can be reformulated as a convex combination of the cell average $\overline{\U}_{j,k}^n$ and high-order sub-cell schemes $\U_{j,k,s}^{(\ell-1),\rm H}$. The BP correction is then implemented via a blending strategy. Specifically, two successive blending operations are performed to enforce physical constraints. First, the high-order sub-cell scheme $\U_{j,k,s}^{(\ell-1),\rm H}$ is blended with its first-order BP counterpart $\U_{j,k,s}^{(\ell-1),\rm L}$ to produce an intermediate sub-cell scheme $\U_{j,k,s}^{(\ell-1),\rm H,*}$ with strictly positive density. Second, the intermediate sub-cell scheme $\U_{j,k,s}^{(\ell-1),\rm H,*}$ is further blended with the first-order sub-cell scheme $\U_{j,k,s}^{(\ell-1),\rm L}$ to yield the eventual corrected sub-cell scheme $\U_{j,k,s}^{(\ell-1),\rm H,(\mathrm{mod})}$, which admits both positive density and pressure, and is subsequently adopted to replace the original high-order sub-cell scheme $\U_{j,k,s}^{(\ell-1),\rm H}$ in the convex decomposition.

In the following theorem, we seek for conditions such that the first-order sub-cell schemes $\U_{j,k,s}^{(\ell-1),\rm L}$ are BP.
\begin{theorem}\label{Th:BP-CFL-Euler-l}
Suppose that all cell averages satisfy $\overline{\U}_{j,k}^{(\ell-1)}\in\mG$. The first-order sub-cell scheme $\U_{j,k,s}^{(\ell-1),\rm L}$ fulfills the BP property for all grid indices $j, k$, and all edge indices $s$, provided the area positivity condition \eqref{EQ:BP-CFL-Scalar-area} and the CFL constraint \eqref{EQ:BP-CFL-Scalar-dt} are satisfied.
\end{theorem}

\begin{proof}
The proof is conducted utilizing the rotational invariance of the Euler equations. Let the unit outward normal vector be $\n_{j,k,s}^{(\ell-1)} = (\cos\theta, \sin\theta)$ with $\theta\in [0,2\pi)$. As demonstrated in \cite{Toro2009Riemann}, the following rotational invariance holds for arbitrary $\U$:
\begin{equation}
\bH(\U)\cdot\n_{j,k,s}^{(\ell-1)} = \T_{j,k,s}^{-1} \F(\T_{j,k,s} \U),
\end{equation}
where $\T_{j,k,s}$ and $\T_{j,k,s}^{-1}$, which denote the rotational matrix and its inverse, respectively, are defined by
\begin{equation*}
\T_{j,k,s} = \begin{bmatrix}
1 &0          &0           &0\\
0 &\cos\theta &-\sin\theta &0\\
0 &\sin\theta & \cos\theta &0\\
0 &0          &0           &1 
\end{bmatrix}, \qquad
\T_{j,k,s}^{-1} = \begin{bmatrix}
1 &0           &0          &0\\
0 & \cos\theta &\sin\theta &0\\
0 &-\sin\theta &\cos\theta &0\\
0 &0           &0          &1 
\end{bmatrix}.
\end{equation*}

Define the rotated conservative variables $\widetilde{\U}_{j,k}^{(\ell-1)} = \T_{j,k,s}\overline{\U}_{j,k}^{(\ell-1)}$ and the rotated mesh velocity $\omega_{j,k,s}^{(\ell-1)} = \w_{j,k,s}\cdot\n_{j,k,s}^{(\ell-1)}$. The following relations naturally hold:
\begin{equation*}
\begin{array}{ll}
\rho[\widetilde{\U}_{j,k}^{(\ell-1)}] = \rho[\overline{\U}_{j,k}^{(\ell-1)}], 
&u[\widetilde{\U}_{j,k}^{(\ell-1)}] = u[\overline{\U}_{j,k}^{(\ell-1)}]\cos\theta + v[\overline{\U}_{j,k}^{(\ell-1)}]\sin\theta, \vspace{+3pt}\\
p[\widetilde{\U}_{j,k}^{(\ell-1)}] = p[\overline{\U}_{j,k}^{(\ell-1)}], 
&v[\widetilde{\U}_{j,k}^{(\ell-1)}] = -u[\overline{\U}_{j,k}^{(\ell-1)}]\sin\theta + v[\overline{\U}_{j,k}^{(\ell-1)}]\cos\theta.\end{array}
\end{equation*}
Consequently, it is straightforward to verify that both density and pressure are invariant under coordinate rotation, which further implies that $\overline{\U}_{j,k}^{(\ell-1)}\in\mG$ is equivalent to $\widetilde{\U}_{j,k}^{(\ell-1)}\in\mG$. In addition, it is apparent that
\begin{equation*}
\begin{aligned}
\widetilde{a}_{j,k,s}^{(\ell-1),-} &:= \min\left\lbrace
(u-c)[\widetilde{\U}_{j,k}^{(\ell-1)}], (u-c)[\widetilde{\U}_{\widehat{j},\widehat{k}}^{(\ell-1)}], -\varepsilon^\mathrm{CU} \right\rbrace
= a_{j,k,s}^{(\ell-1),-}, \\
\widetilde{a}_{j,k,s}^{(\ell-1),+} &:= \max\left\lbrace
(u+c)[\widetilde{\U}_{j,k}^{(\ell-1)}], (u+c)[\widetilde{\U}_{\widehat{j},\widehat{k}}^{(\ell-1)}], \varepsilon^\mathrm{CU} \right\rbrace
= a_{j,k,s}^{(\ell-1),+}.
\end{aligned}
\end{equation*}

The first-order sub-cell scheme for the Euler system can be reformulated into the rotated coordinate system as
\begin{equation*}
\begin{aligned}
\U_{j,k,s}^{(\ell-1),\rm L} &= \overline{\U}_{j,k}^{(\ell-1)} - \lambda_{j,k}^{(\ell-1)} \left[\mbH \left(\w_{j,k,s}, \n_{j,k,s}^{(\ell-1)}, \overline{u}_{j,k}^{(\ell-1)},\overline{u}_{\widehat{j},\widehat{k}}^{(\ell-1)}\right) - \mbH\left(\w_{j,k,s},\n_{j,k,s}^{(\ell-1)},\overline{u}_{j,k}^{(\ell-1)},\overline{u}_{j,k}^{(\ell-1)}\right)\right],\\
&=\T_{j,k,s}^{-1}\left\lbrace\widetilde{\U}_{j,k}^{(\ell-1)} - \lambda_{j,k}^{(\ell-1)} \left[\widetilde{\mbH} \left(\omega_{j,k,s}^{(\ell-1)}, \widetilde{\U}_{j,k}^{(\ell-1)},\widetilde{\U}_{\widehat{j},\widehat{k}}^{(\ell-1)}\right) - \widetilde{\mbH}\left(\omega_{j,k,s}^{(\ell-1)},\widetilde{\U}_{j,k}^{(\ell-1)},\widetilde{\U}_{j,k}^{(\ell-1)}\right)\right]
\right\rbrace := \T_{j,k,s}^{-1}\widetilde{\U}_{j,k,s}^{(\ell-1),\rm L},
\end{aligned}
\end{equation*}
where the rotated flux satisfies 
\begin{equation*}
\widetilde{\mbH} \left(\omega_{j,k,s}^{(\ell-1)}, \widetilde{\U}_{j,k}^{(\ell-1)},\widetilde{\U}_{\widehat{j},\widehat{k}}^{(\ell-1)}\right) 
= \begin{cases}
\F(\widetilde{\U}_{j,k}^{(\ell-1)}) - \omega_{j,k,s}^{(\ell-1)}\widetilde{\U}_{j,k}^{(\ell-1)}, 
&\omega_{j,k,s}^{(\ell-1)} < \widetilde{a}_{j,k,s}^{(\ell-1),-}, \vspace{+3pt}\\

\widetilde{\F}_{j,k,s}^{(\ell-1),*} - \omega_{j,k,s}^{(\ell-1)}\widetilde{\U}_{j,k,s}^{(\ell-1),*}, 
&\widetilde{a}_{j,k,s}^{(\ell-1),-} \leq \omega_{j,k,s}^{(\ell-1)} \leq \widetilde{a}_{j,k,s}^{(\ell-1),+}, \vspace{+3pt}\\

\F(\widetilde{\U}_{\widehat{j},\widehat{k}}^{(\ell-1)}) - \omega_{j,k,s}^{(\ell-1)}\widetilde{\U}_{\widehat{j},\widehat{k}}^{(\ell-1)}, &\omega_{j,k,s}^{(\ell-1)} > \widetilde{a}_{j,k,s}^{(\ell-1),+},
\end{cases}
\end{equation*}
with the rotated flux $\widetilde{\F}_{j,k,s}^{(\ell-1),*}$ and intermediate state $\widetilde{\U}_{j,k,s}^{(\ell-1),*}$ being expressed by
\begin{equation}
\begin{aligned}
\widetilde{\F}_{j,k,s}^{(\ell-1),*}&=\dfrac{\widetilde{a}_{j,k,s}^{(\ell-1),+} \F(\widetilde{\U}_{j,k}^{(\ell-1)}) - \widetilde{a}_{j,k,s}^{(\ell-1),-} \F(\widetilde{\U}_{\widehat{j},\widehat{k}}^{(\ell-1)})}{\widetilde{a}_{j,k,s}^{(\ell-1),+} - \widetilde{a}_{j,k,s}^{(\ell-1),-}} + \dfrac{\widetilde{a}_{j,k,s}^{(\ell-1),-} \widetilde{a}_{j,k,s}^{(\ell-1),+}}{\widetilde{a}_{j,k,s}^{(\ell-1),+} - \widetilde{a}_{j,k,s}^{(\ell-1),-}}(\widetilde{\U}_{\widehat{j},\widehat{k}}^{(\ell-1)} - \widetilde{\U}_{j,k}^{(\ell-1)}), \\
\widetilde{\U}_{j,k,s}^{(\ell-1),*} &= \dfrac{\widetilde{a}_{j,k,s}^+ \widetilde{\U}_{\widehat{j},\widehat{k}}^{(\ell-1)} - \widetilde{a}_{j,k,s}^-\widetilde{\U}_{j,k}^{(\ell-1)} + \F(\widetilde{\U}_{j,k}^{(\ell-1)}) - \F(\widetilde{\U}_{\widehat{j},\widehat{k}}^{(\ell-1)})}{\widetilde{a}_{j,k,s}^+ - \widetilde{a}_{j,k,s}^-},
\end{aligned}
\end{equation}
and the auxiliary flux is defined as
\begin{equation*}
\widetilde{\mbH}\left(\omega_{j,k,s}^{(\ell-1)},\widetilde{\U}_{j,k}^{(\ell-1)},\widetilde{\U}_{j,k}^{(\ell-1)}\right) = \F(\widetilde{\U}_{j,k}^{(\ell-1)}) - \omega_{j,k,s}^{(\ell-1)}\widetilde{\U}_{j,k}^{(\ell-1)}.
\end{equation*}

It suffices to prove that $\widetilde{\U}_{j,k,s}^{(\ell-1),\rm L}\in\mG$ under conditions \eqref{EQ:BP-CFL-Scalar-area} and \eqref{EQ:BP-CFL-Scalar-dt}, which is verified according to the following three cases.

\noindent\textbf{Case I}: $\omega_{j,k,s}^{(\ell-1)} < \widetilde{a}_{j,k,s}^{(\ell-1),-}$. In this case, the rotated sub-cell scheme reduces to
\begin{equation}
\widetilde{\U}_{j,k,s}^{(\ell-1),\rm L} = 
\widetilde{\U}_{j,k}^{(\ell-1)} - \lambda_{j,k}^{(\ell-1)}\left[
\F(\widetilde{\U}_{j,k}^{(\ell-1)}) - \omega_{j,k,s}^{(\ell-1)}\widetilde{\U}_{j,k}^{(\ell-1)} - \left( \F(\widetilde{\U}_{j,k}^{(\ell-1)}) - \omega_{j,k,s}^{(\ell-1)}\widetilde{\U}_{j,k}^{(\ell-1)} \right)
\right] = \widetilde{\U}_{j,k}^{(\ell-1)}.
\end{equation}
The BP property $\widetilde{\U}_{j,k,s}^{(\ell-1),\rm L}\in\mG$ holds trivially, and the CFL constraint \eqref{EQ:BP-CFL-Scalar-dt} is automatically satisfied.

\noindent\textbf{Case II}: $\omega_{j,k,s}^{(\ell-1)} > \widetilde{a}_{j,k,s}^{(\ell-1),+}$. In this case, The rotated sub-cell scheme $\widetilde{\U}_{j,k,s}^{(\ell-1),\rm L}$ admits the decomposition
\begin{equation*}
\begin{aligned}
\widetilde{\U}_{j,k,s}^{(\ell-1),\rm L} &= \widetilde{\U}_{j,k}^{(\ell-1)} - \lambda_{j,k}^{(\ell-1)}\left[
\F(\widetilde{\U}_{\widehat{j},\widehat{k}}^{(\ell-1)}) - \omega_{j,k,s}^{(\ell-1)}\widetilde{\U}_{\widehat{j},\widehat{k}}^{(\ell-1)} - \left( \F(\widetilde{\U}_{j,k}^{(\ell-1)}) - \omega_{j,k,s}^{(\ell-1)}\widetilde{\U}_{j,k}^{(\ell-1)} \right)
\right]\\
&= \left[ 1- \lambda_{j,k}^{(\ell-1)}\left(\omega_{j,k,s}^{(\ell-1)} - \widetilde{a}_{j,k,s}^- \right) \right]\widetilde{\U}_{j,k}^{(\ell-1)} + \lambda_{j,k}^{(\ell-1)}\left(\omega_{j,k,s}^{(\ell-1)} - \widetilde{a}_{j,k,s}^+ \right)\widetilde{\U}_{\widehat{j},\widehat{k}}^{(\ell-1)} + \lambda_{j,k}^{(\ell-1)}\left(\widetilde{a}_{j,k,s}^+ -\widetilde{a}_{j,k,s}^-\right)\widetilde{\U}_{j,k,s}^{(\ell-1),*}.
\end{aligned}
\end{equation*}
Under the assumptions $\overline{\U}_{j,k}^{(\ell-1)}\in\mG$ and $\overline{\U}_{\widehat{j},\widehat{k}}^{(\ell-1)}\in\mG$, the rotated states satisfy $\widetilde{\U}_{j,k}^{(\ell-1)}\in\mG$ and $\widetilde{\U}_{\widehat{j},\widehat{k}}^{(\ell-1)}\in\mG$. By virtue of the geometric quasilinearization (GQL) framework proposed in \cite{Cui2024Bound,Wu2023Geometric}, the intermediate state $\widetilde{\U}_{j,k,s}^{(\ell-1),*}$ remains within the invariant domain $\mG$. When the constraints \eqref{EQ:BP-CFL-Scalar-area} and \eqref{EQ:BP-CFL-Scalar-dt} are satisfied, the first-order sub-cell scheme is a convex combination of $\widetilde{\U}_{j,k}^{(\ell-1)}$, $\widetilde{\U}_{\widehat{j},\widehat{k}}^{(\ell-1)}$, and $\widetilde{\U}_{j,k,s}^{(\ell-1),*}$, which guarantees $\widetilde{\U}_{j,k,s}^{(\ell-1),\rm L}\in\mG$, and hence $\U_{j,k,s}^{(\ell-1),\rm L}\in\mG$.

\noindent\textbf{Case III}: $\widetilde{a}_{j,k,s}^{(\ell-1),-} \leq \omega_{j,k,s}^{(\ell-1)} \leq \widetilde{a}_{j,k,s}^{(\ell-1),+}$. In this case, the sub-cell scheme $\U_{j,k,s}^{(\ell-1),\rm L}$ can be decomposed as
\begin{equation*}
\widetilde{\U}_{j,k,s}^{(\ell-1),\rm L} = \left[ 1- \lambda_{j,k}^{(\ell-1)}\left(\omega_{j,k,s}^{(\ell-1)} - \widetilde{a}_{j,k,s}^{(\ell-1),-} \right) \right]\widetilde{\U}_{j,k}^{(\ell-1)} + \lambda_{j,k}^{(\ell-1)}\left(\omega_{j,k,s}^{(\ell-1)} -\widetilde{a}_{j,k,s}^{(\ell-1),-} \right)\widetilde{\U}_{j,k,s}^{(\ell-1),*}.
\end{equation*}
Under the conditions \eqref{EQ:BP-CFL-Scalar-area} and \eqref{EQ:BP-CFL-Scalar-dt}, the first-order sub-cell scheme $\widetilde{\U}_{j,k,s}^{(\ell-1),\rm L}$ is
 a convex combination of $\U_{j,k}^{(\ell-1)}$ and $\widetilde{\U}_{j,k,s}^{(\ell-1),*}$. Therefore, $\widetilde{\U}_{j,k,s}^{(\ell-1),\rm L}\in\mG$, which implies $\U_{j,k,s}^{(\ell-1),\rm L}\in\mG$.
 
This completes the proof of Theorem \ref{Th:BP-CFL-Euler-l}.
\end{proof}

It is remarked that the D-GCL is preserved by the BP corrections, which can be proved the same way as in Theorem \ref{Th:GCL}. Therefore, we omit the proof for the sake of brevity. 

Taking all three stages of SSP-RK3 into account, we establish the global BP property for the fully discrete SSP-RK3 AMM-FVM, as stated in the following theorem.
\begin{theorem}\label{Th:BP-CFL-Euler}
Suppose that the initial cell average satisfies $\overline{\U}_{j,k}^{n}\in\mG$ for all $j,k$. The fully discrete scheme satisfies $\overline{\U}_{j,k}^{n+1}\in\mG$ for all grid indices, provided that the conditions \eqref{EQ:BP-CFL-Scalar-area} and \eqref{EQ:BP-CFL-Scalar-dt} hold for all indices $j$, $k$, $s$, and $\ell$.
\end{theorem}

\begin{proof}
The global BP property follows directly from mathematical induction based on Theorem \ref{Th:BP-CFL-Euler-l}. The detailed derivation is omitted for the sake of brevity.
\end{proof}

%\begin{remark}
%Characteristic projections are applied to reduce numerical oscillations. 
%\end{remark}

\section{Numerical examples}\label{Sec:Num}
In this section, we first test the accuracy of the BP AMM-FVM in Example \ref{Ex:Ad-Accuracy} and Example \ref{Ex:Burgers}, and then verify its robustness through several numerical examples. We note that the BP property of the proposed scheme relies on the D-GCL of the algorithm, which will be examined in Example \ref{Ex:GCL}. Unless stated otherwise, the CFL number defined in \eqref{EQ:Conventional-dt} is set to $\CFL=1$ for problems involving shock waves. 
In addition, the gradient of the conservative variable $\dVert{D\U}_{j,k}=\dvert{\nabla u}_{j,k}$ will be taken for scalar problems, and the gradient of the density $\dVert{D\U}_{j,k}=\dvert{\nabla \rho}_{j,k}$ will be employed for Euler equations in all numerical experiments, unless specified otherwise.
The computational efficiency is also investigated. We demonstrate that, for a fixed number of cells, the AMM-FVM captures sharper wave structures than the standard UM-FVM; see Example \ref{Ex:TP} and Example \ref{Ex:BW}. We also show that the AMM-FVM achieves superior efficiency: for comparable CPU times, it requires fewer cells while yielding sharper wave structures compared to the UM-FVM; see Example \ref{Ex:Burgers} and Example \ref{Ex:SVI}.
For the sake of brevity, the finite volume methods with and without BP limiters are denoted by $\textbf{BP}$ and $\textbf{NBP}$, respectively. Besides, the results obtained on adaptive moving mesh and uniform mesh are labeled $\textbf{MM}$ and $\textbf{UM}$, respectively.

\begin{example}\label{Ex:Ad-Accuracy} 
\em advection equation.
\end{example}
In the first example, the accuracy is examined by solving the linear advection equation:
\begin{equation*}
\left\lbrace
\begin{array}{l}
u_t + u_x + u_y = 0, \quad (x,y)\in \Omega = [0,1]^2,\\
u(x,y,0) = u_0(x,y) = \sin^2(2\pi x) \sin^2(2\pi y),
\end{array}
\right.
\end{equation*}
subject to periodic boundary conditions on all boundaries of the physical domain. The exact solution is given by $u(x,y,t) = u_0(x-t,y-t)$. 

The invariant domain for this test case is $\mG = \{u \mid 0\leq u\leq 1\}$. Without BP limiters, the third-order AMM-FVM fails to preserve the lower bound. We therefore solve the advection equation using the BP AMM-FVM proposed in this paper. In the simulation, we set $\beta=0.6$. The CFL number is taken as $1$. 
Table \ref{Tab:Ad-Accuracy} reports the $L^1$ errors $\varepsilon=\dVert{u-1}_1$, the maximum of areas of the control volumes $\dvert{C_{j,k}}_{\max}$, the corresponding convergence rates, and the minimum of the numerical solution $u_{\min}$ at the output time $t=1$. Computations start from a mesh with $N=N_x=N_y$, and the convergence rate is computed as
\begin{equation*}
\mathcal{R}^{2N} = \dfrac{\log_{10} \left( \varepsilon^{N}/\varepsilon^{2N} \right)}{\log_{10} \sqrt{\dvert{C_{j,k}}_{\max}^{N}/\dvert{C_{j,k}}_{\max}^{2N}}}.
\end{equation*}
It can be observed that the AMM-FVM scheme achieves the designed third-order of accuracy. Furthermore, the proposed BP AMM-FVM scheme preserves the invariant domain, whereas the standard AMM-FVM without BP limiters does not.
\begin{table}[htb!]
\caption{\sf Example \ref{Ex:Ad-Accuracy}: Convergence of AMM-FVM schemes with/without BP limiters, together with minimum of solution. The CFL number is taken as $1$, and $\beta=0.6$.}
\label{Tab:Ad-Accuracy}
\begin{center}
\begin{tabular}{||c|cccc|cccc||}
\hline
\multirow{2}{*}{$N$}
&\multicolumn{4}{c}{Without BP Limiters, $\CFL=1$, $\beta=0.6$}
&\multicolumn{4}{c||}{With BP Limiters, $\CFL=1$, $\beta=0.6$} \\ \cline{2-9}
&$L^1$-error &$\dvert{C_{j,k}}_{\max}$ &Rate $\mathcal{R}$ &$u_{\min}$
&$L^1$-error &$\dvert{C_{j,k}}_{\max}$ &Rate $\mathcal{R}$ &$u_{\min}$\\ \hline
$ 40$ &1.85E-02 &9.29E-04 &--   &-3.46E-02
      &1.65E-02 &1.02E-03 &--   & 3.22E-07 \\

$ 80$ &2.55E-03 &2.31E-04 &2.85 &-6.99E-03
      &2.43E-03 &2.40E-04 &2.65 & 2.51E-10 \\
      
$160$ &3.24E-04 &5.95E-05 &3.04 &-6.93E-04
      &3.14E-04 &5.93E-05 &2.93 & 1.27E-13 \\

$320$ &4.11E-05 &1.49E-05 &2.98 &-7.06E-05
      &4.09E-05 &1.49E-05 &2.95 & 1.03E-16 \\
\hline
\end{tabular}
\end{center}
\end{table}

\begin{example}\label{Ex:GCL} 
\em Discrete geometric conservation law.
\end{example}
In the second set of numerical tests, the D-GCL is verified by solving the same linear advection equation as in Example \ref{Ex:Ad-Accuracy}, but with the uniform initial condition $u(x,y,0) \equiv 1$, $(x,y)\in [0,1]^2$. The boundary conditions and final time are identical to those in Example \ref{Ex:Ad-Accuracy}. The physical domain is initially partitioned into $N\times N$ uniform cells, and the grid points are then redistributed dynamically in time.

\begin{figure}[ht!]
\centering
\includegraphics[width=0.31\textwidth]{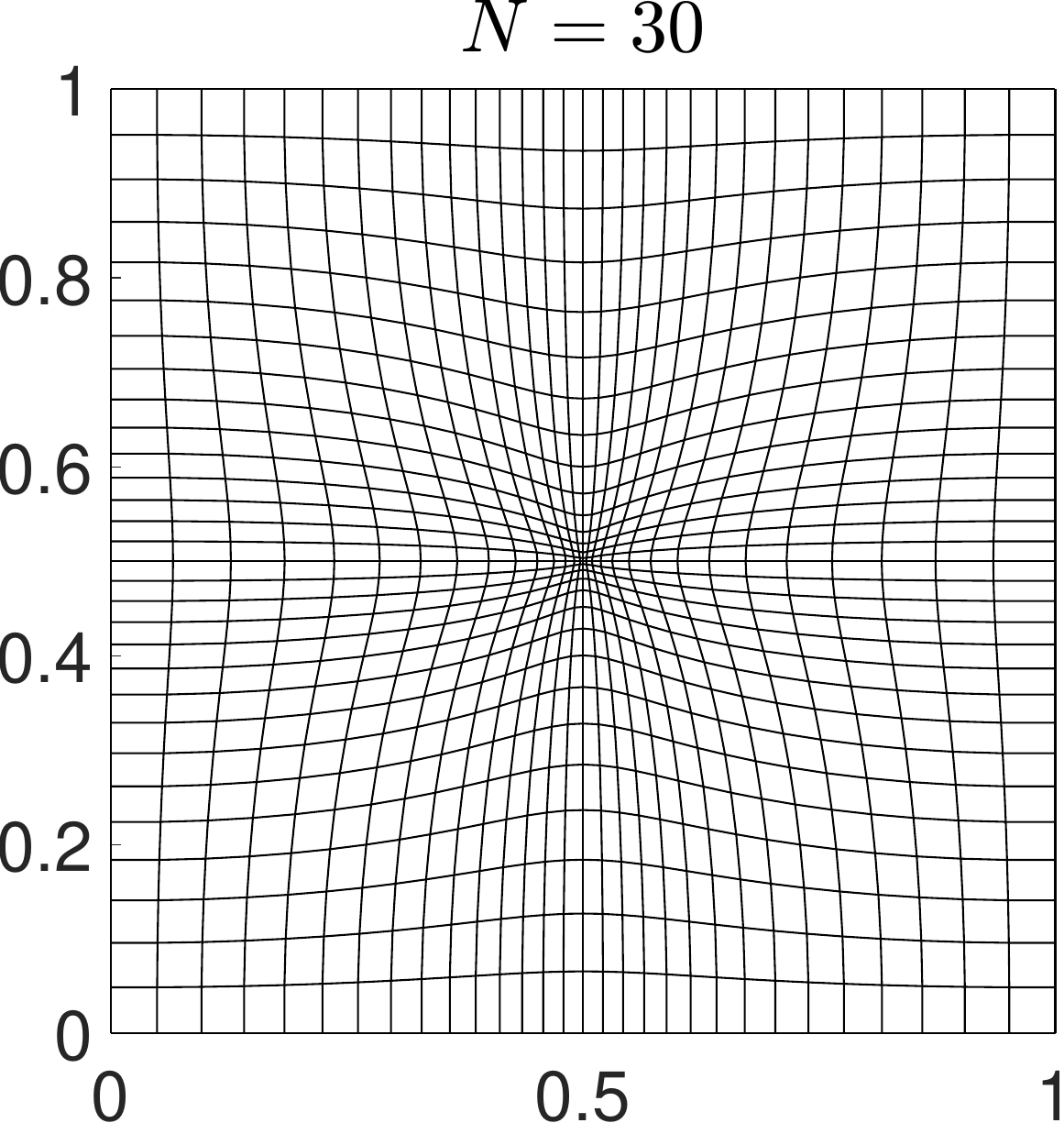}\qquad\qquad
\includegraphics[width=0.31\textwidth]{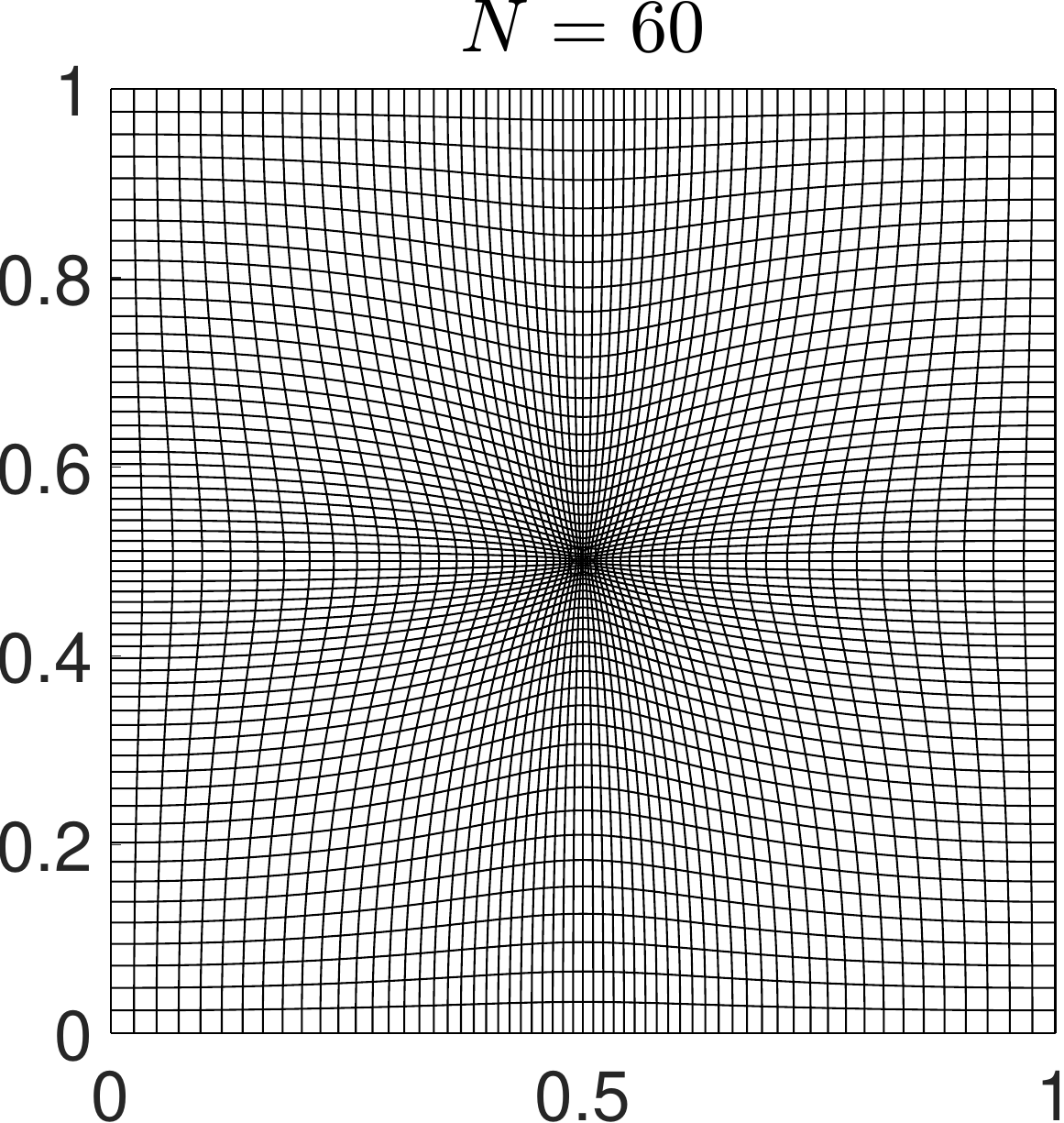}
\caption{\sf Example \ref{Ex:GCL}: Meshes for $N=N_x=N_y=30$ (left) and $N=N_x=N_y=60$ (right) at the output time $t=1$ in the D-GCL test.
\label{Fig:GCL-Mesh}}
\end{figure}

\begin{figure}[ht!]
\centering
\includegraphics[width=0.5\textwidth]{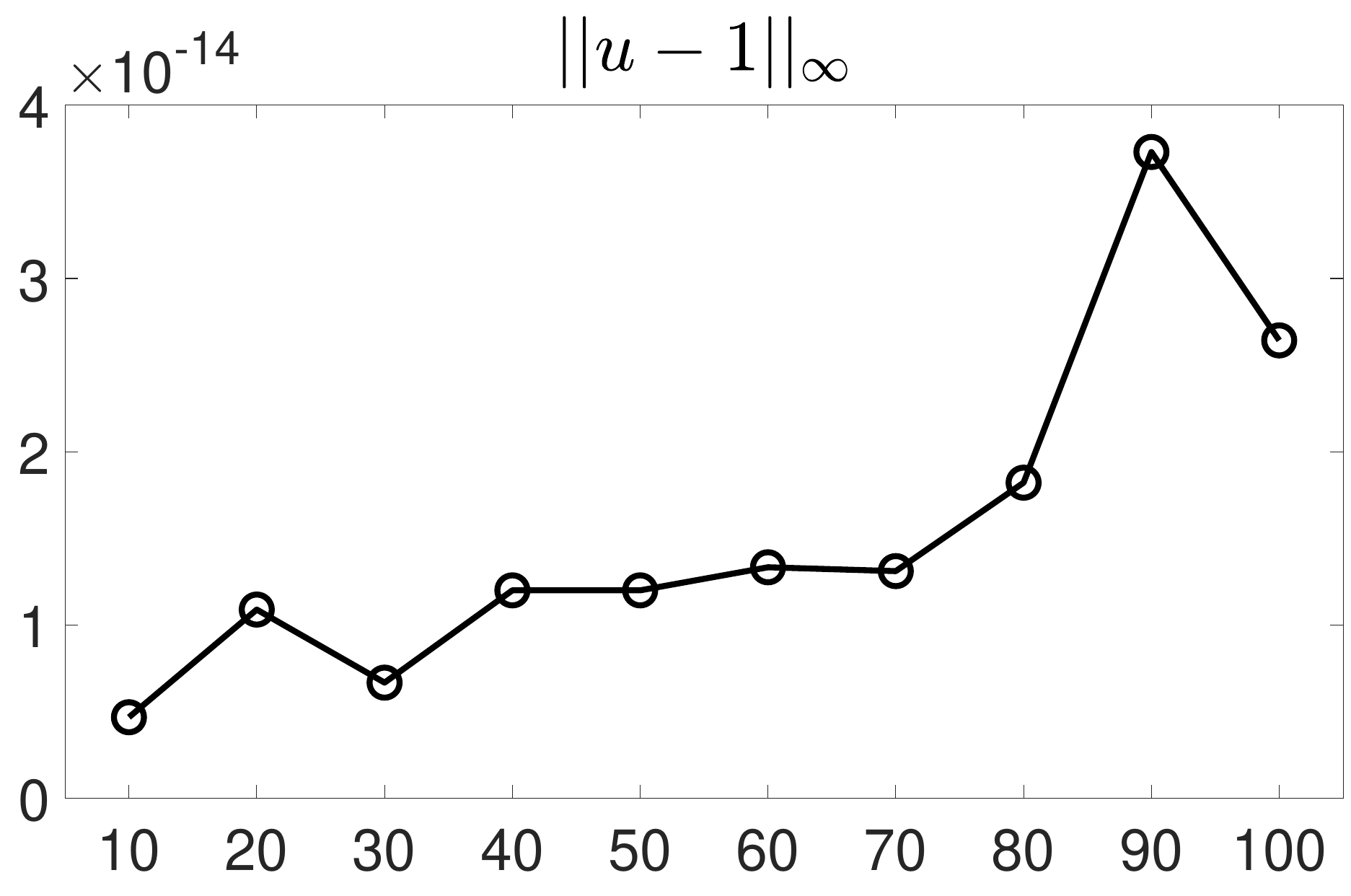}
\caption{\sf Example \ref{Ex:GCL}: The $L^\infty$ errors for different $N$ in the D-GCL test.
\label{Fig:GCL}}
\end{figure}

We note that the choice $Du=\dvert{\nabla u}$ is not suitable for this test case. The underlying reason is that the mesh remains stationary under this configuration, since the variation of $u$ is identically zero throughout the simulation.
To generate moving meshes, we set $Du = \nicefrac{1}{r}$, with $r=\sqrt{(x-0.5)^2+(y-0.5)^2}$; see for instance Figure \ref{Fig:GCL-Mesh}, which displays the resulting meshes at the final time $t=1$ for $N=N_x=N_y=30$ and $N=N_x=N_y=60$. The $L^\infty$ errors for various mesh resolutions are displayed in Figure \ref{Fig:GCL}. The D-GCL is clearly validated, as the errors are close to machine zero.

\begin{example}\label{Ex:Burgers} 
\em Burgers equation.
\end{example}
In the third example, the BP property of the proposed BP AMM-FVM is verified via the two-dimensional Burgers equation
\begin{equation}
u_t + \left(\frac{1}{2}u^2\right)_x + \left(\frac{1}{2}u^2\right)_y = 0,
\end{equation}
which is solved over the domain $\Omega=[-1,1]^2$. The setup for this problem is taken from \cite{ZHANG20103091}. The initial condition is given by
\begin{equation}
u(x,y,0) = u_0(x,y) = \sin(\pi (x+y)) + 0.5,
\end{equation}
with periodic boundary conditions imposed on the physical boundary. The invariant domain for this problem is $\mG = \{u \mid -0.5\leq u\leq 1.5\}$.

We first compute the Burgers equation up to $t=0.05$ when the solution remains smooth. The exact solution is obtained by Newton iterations. In the simulation, we use $\beta=0.3$ and $\beta=0.6$ to assess the accuracy. The CFL number is set to $1$. 
The $L^1$ errors, the maximum of areas of the control volumes $\dvert{C_{j,k}}_{\max}$, the corresponding convergence rates, and the discrepancy between the exact upper bound $u_M$ and the maximum of numerical solution $u_{\max}$ are reported in Table \ref{Tab:Burgers-Accuracy-1} for $\beta=0.3$ and in Table \ref{Tab:Burgers-Accuracy-2} for $\beta=0.6$. All computations start from a uniform mesh with $N=N_x=N_y$. It can be observed that the AMM-FVM scheme achieves third-order accuracy, but fails to preserve the upper bound. When BP limiters are applied, the bounds are preserved, while the accuracy is slightly affected. Such slight accuracy loss can be recovered by reducing the CFL number. In Table \ref{Tab:Burgers-Accuracy-0p4-1}, we set $\CFL=0.3$ and employ BP limiters to enforce the constraint $u\in \{u \mid -0.5+\varepsilon^\mathrm{BP}\leq u\leq 1.5-\varepsilon^\mathrm{BP}\}$. Clearly, the scheme maintains third-order of accuracy even when the BP limiter is activated.

\begin{table}[htb!]
\caption{\sf Example \ref{Ex:Burgers}: Convergence of AMM-FVM schemes with/without BP limiters, together with the errors $u_M-u_{\max}$. The output time is $t=0.05$. The CFL number is taken as $1$, and $\beta$ is taken as $0.3$.}
\label{Tab:Burgers-Accuracy-1}
\begin{center}
\begin{tabular}{||c|cccc|cccc||}
\hline
\multirow{2}{*}{$N$}
&\multicolumn{4}{c}{Without BP Limiters, $\CFL=1$, $\beta=0.3$}
&\multicolumn{4}{c||}{With BP Limiters, $\CFL=1$, $\beta=0.3$} \\ \cline{2-9}
&$L^1$-error &$\dvert{C_{j,k}}_{\max}$ &Rate $\mathcal{R}$ &$u_M-u_{\max}$
&$L^1$-error &$\dvert{C_{j,k}}_{\max}$ &Rate $\mathcal{R}$ &$u_M-u_{\max}$\\ \hline
$ 20$ &6.65E-03 &1.04E-02 &--   &-3.09E-02
      &6.41E-03 &1.04E-02 &--   & 4.44E-03 \\
      
$ 40$ &1.03E-03 &2.83E-03 &2.87 &-8.88E-03
      &9.57E-04 &2.79E-03 &2.89 & 7.22E-04 \\

$ 80$ &1.40E-04 &7.16E-04 &2.90 &-2.39E-03
      &1.51E-04 &7.15E-04 &2.71 & 1.45E-04 \\
      
$160$ &1.82E-05 &1.76E-04 &2.91 &-6.20E-04
      &2.55E-05 &1.75E-04 &2.53 & 3.34E-05 \\
%
%$320$ &2.37E-06 &4.28E-05 &2.88 &-1.58E-04
%      &3.92E-06 &5.10E-05 &2.51 & 8.74E-06 \\
\hline
\end{tabular}
\end{center}
\end{table}

\begin{table}[htb!]
\caption{\sf Example \ref{Ex:Burgers}: Convergence of AMM-FVM schemes with/without BP limiters, together with the errors $u_M-u_{\max}$. The output time is $t=0.05$. The CFL number is taken as $1$, and $\beta$ is taken as $0.6$.}
\label{Tab:Burgers-Accuracy-2}
\begin{center}
\begin{tabular}{||c|cccc|cccc||}
\hline
\multirow{2}{*}{$N$}
&\multicolumn{4}{c}{Without BP Limiters, $\CFL=1$, $\beta=0.6$}
&\multicolumn{4}{c||}{With BP Limiters, $\CFL=1$, $\beta=0.6$} \\ \cline{2-9}
&$L^1$-error &$\dvert{C_{j,k}}_{\max}$ &Rate $\mathcal{R}$ &$u_M-u_{\max}$
&$L^1$-error &$\dvert{C_{j,k}}_{\max}$ &Rate $\mathcal{R}$ &$u_M-u_{\max}$\\ \hline
$ 20$ &6.62E-03 &1.08E-02 &--   &-2.86E-02
      &6.25E-03 &1.08E-02 &--   & 4.74E-03 \\
      
$ 40$ &1.08E-03 &3.18E-03 &2.96 &-7.79E-03
      &9.69E-04 &3.10E-03 &2.99 & 8.30E-04 \\

$ 80$ &1.54E-04 &8.44E-04 &2.94 &-2.15E-03
      &1.48E-04 &8.43E-04 &2.89 & 1.69E-04 \\
      
$160$ &2.27E-05 &2.13E-04 &2.78 &-5.78E-04
      &2.36E-05 &2.13E-04 &2.67 & 3.76E-05 \\
%
%$320$ &2.37E-06 &4.28E-05 &2.88 &-1.58E-04
%      &3.92E-06 &5.10E-05 &2.51 & 8.74E-06 \\
\hline
\end{tabular}
\end{center}
\end{table}

\begin{table}[htb!]
\caption{\sf Example \ref{Ex:Burgers}: Convergence of AMM-FVM schemes with/without BP limiters, together with the errors $u_M-u_{\max}$. The output time is $t=0.05$. The CFL number is taken as $0.3$, and $\beta$ is taken as $0.3$ and $0.6$.}
\label{Tab:Burgers-Accuracy-0p4-1}
\begin{center}
\begin{tabular}{||c|cccc|cccc||}
\hline
\multirow{2}{*}{$N$}
&\multicolumn{4}{c||}{With BP Limiters, $\CFL=0.3$, $\beta=0.3$}
&\multicolumn{4}{c||}{With BP Limiters, $\CFL=0.3$, $\beta=0.6$} \\ \cline{2-9}
&$L^1$-error &$\dvert{C_{j,k}}_{\max}$ &Rate $\mathcal{R}$ &$u_M-u_{\max}$
&$L^1$-error &$\dvert{C_{j,k}}_{\max}$ &Rate $\mathcal{R}$ &$u_M-u_{\max}$\\ \hline
$ 20$ &5.98E-03 &1.06E-02 &--   & 7.53E-03 
      &5.93E-03 &1.12E-02 &--   & 7.99E-03 \\
      
$ 40$ &8.78E-04 &3.01E-03 &3.05 & 1.68E-03 
      &8.97E-04 &3.51E-03 &3.26 & 2.03E-03 \\

$ 80$ &1.22E-04 &7.98E-04 &2.97 & 3.61E-04 
      &1.48E-04 &1.04E-04 &2.96 & 4.57E-04 \\
      
$160$ &1.64E-05 &1.98E-04 &2.89 & 8.17E-05
      &2.18E-05 &2.76E-04 &2.89 & 1.02E-05 \\
%
%$320$ &2.37E-06 &4.28E-05 &2.88 &-1.58E-04
%      &4.11E-06 &4.13E-05 &2.61 & 7.94E-06 \\
\hline
\end{tabular}
\end{center}
\end{table}

%\begin{table}[htb!]
%\caption{\sf Example \ref{Ex:Burgers}: Convergence of AMM-FVM schemes with/without BP limiters, together with the errors $u_M-u_{\max}$. The output time is $t=0.05$. The CFL number is taken as $0.6$.}
%\label{Tab:Burgers-Accuracy-0p4-2}
%\begin{center}
%\begin{tabular}{||c|cccc||}
%\hline
%\multirow{2}{*}{$N$}
%&\multicolumn{4}{c||}{With BP Limiters, $\CFL=0.4$} \\ \cline{2-5}
%&$L^1$-error &$\dvert{C_{j,k}}_{\max}$ &Rate $\mathcal{R}$ &$u_M-u_{\max}$\\ \hline
%$ 20$ &5.93E-03 &1.12E-02 &--   & 7.99E-03 \\
%      
%$ 40$ &8.97E-04 &3.51E-03 &3.26 & 2.03E-03 \\
%
%$ 80$ &1.48E-04 &1.04E-04 &2.96 & 4.57E-04 \\
%      
%$160$ &2.18E-05 &2.76E-04 &2.89 & 1.02E-05 \\
%%
%%$320$ &2.37E-06 &4.28E-05 &2.88 &-1.58E-04
%%      &4.11E-06 &4.13E-05 &2.61 & 7.94E-06 \\
%\hline
%\end{tabular}
%\end{center}
%\end{table}

\begin{figure}[ht!]
\centering
\includegraphics[width=0.31\textwidth]{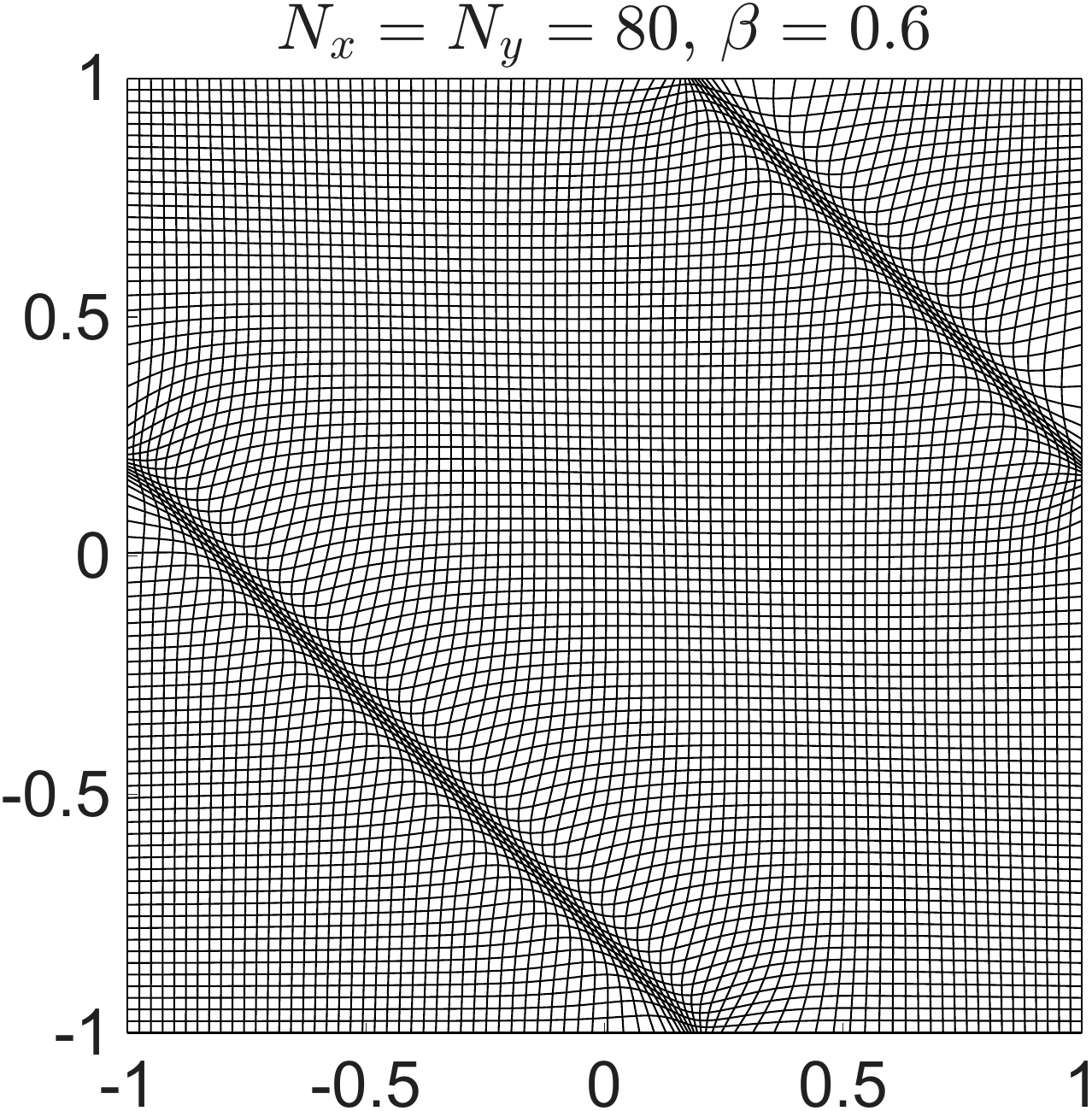}\quad
\includegraphics[width=0.31\textwidth]{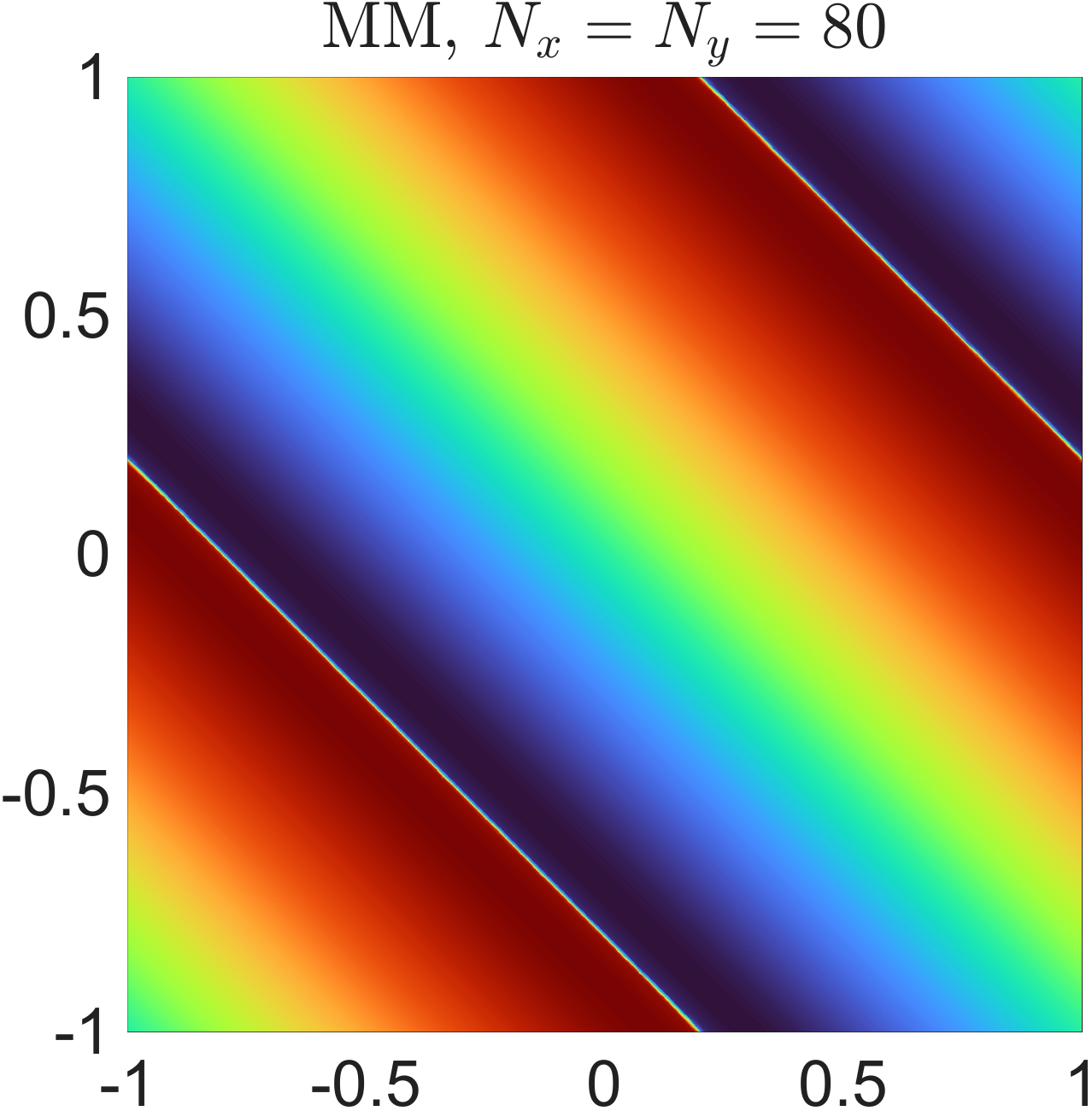}\quad
\includegraphics[width=0.31\textwidth]{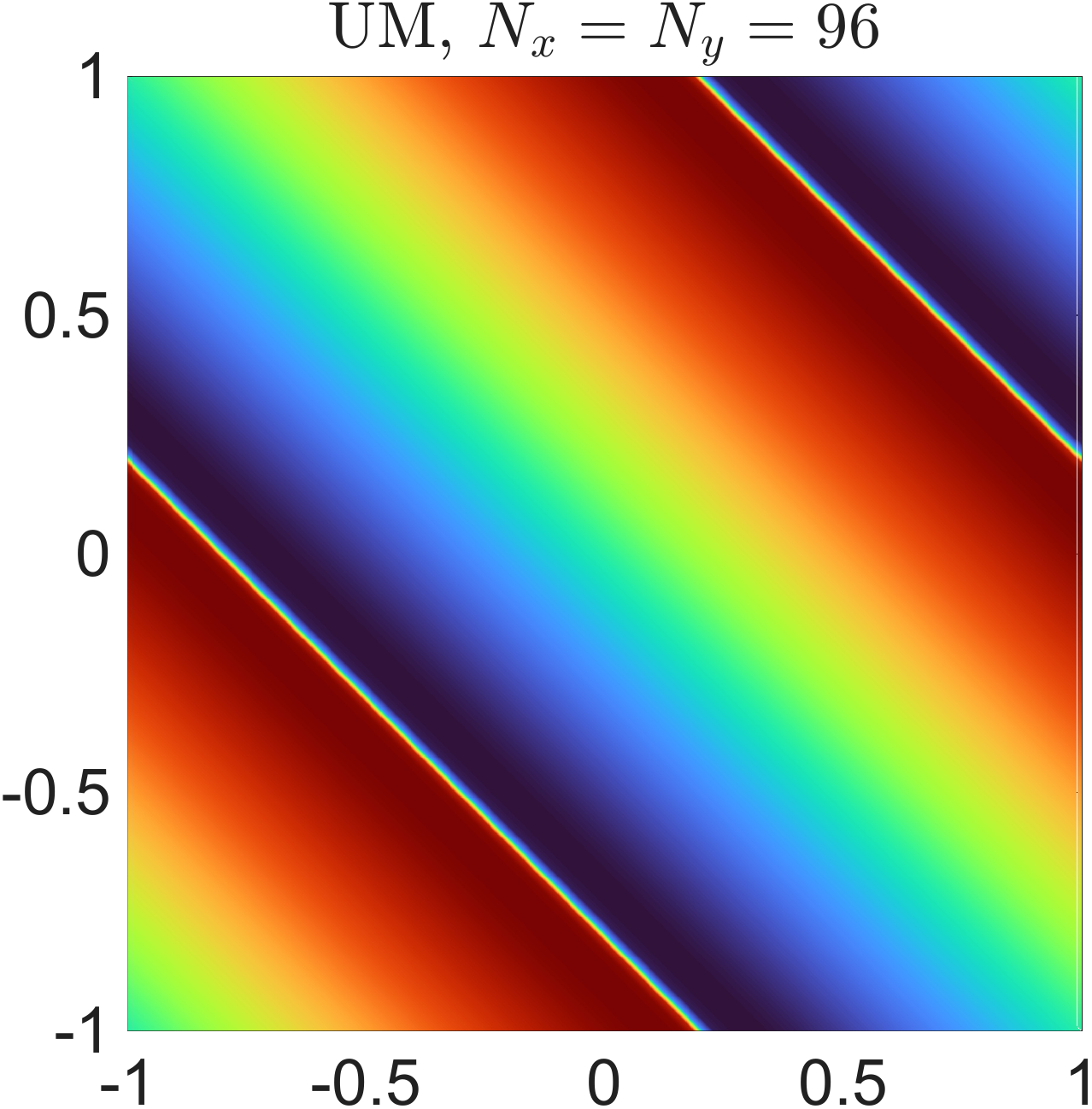}
\caption{\sf Example \ref{Ex:Burgers}. Left: computational mesh. Middle: numerical solution on MM using mesh resolution $N_x=N_y=80$. Right: numerical solution on UM using mesh resolution $N_x=N_y=96$.
\label{Fig:Burgers-Mesh-Solu}}
\end{figure}

Next, we simulate the problem up to $t=0.2$, at which two shock waves have formed. In this regime, the solution is discontinuous, so that the BP property becomes the primary concern rather than formal accuracy. We therefore set $\beta=0.6$ to concentrate grid points adaptively and use $\CFL=1$ to enhance computational efficiency. Computations are performed on an adaptive moving mesh containing $N_x\times N_y = 80\times 80$ cells. The adaptive mesh and corresponding numerical solution are displayed in the left two panels of Figure \ref{Fig:Burgers-Mesh-Solu}. It can be seen that the grid points are adaptively concentrated near the shock waves. Consequently, both shocks are well resolved without spurious oscillations. 
To further demonstrate the efficiency of the adaptive moving mesh method, we compare the solution produced by BP AMM-FVM with that produced by BP UM-FVM. Here, we take a $96\times 96$ uniform mesh, and the corresponding solution is exhibited in the right panel of Figure \ref{Fig:Burgers-Mesh-Solu}. The CPU time for AMM-FVM with $N_x\times N_y = 80\times 80$ is $1.68$s, and the CPU time for UM-FVM with $N_x\times N_y = 96\times 96$ is $1.72$s.
As can be seen from Figure \ref{Fig:Burgers-Mesh-Solu}, when comparable CPU times are consumed, the AMM-FVM produces solution with sharper shock waves than UM-FVM, but requires fewer grid points.
In other words, if the resolution of solution produced by FVM on uniform mesh is comparable to that on adaptive moving mesh, more gird points need to be used for UM-FVM, and more CPU time will be consumed.

\begin{example}\label{Ex:TP} 
\em Triple-point problem of Euler equations.
\end{example}

In this example, the triple-point problem is used to validate performance of the BP AMM-FVM upon the occurrence of large mesh deformation. We compute the Euler equations on $[0,9]\times [0,3]$ with initial data
\begin{equation*}
(\rho, u, v, p) = \begin{cases}
(1, 0, 0, 1), & [0,1]\times [0,3],\\
(0.125, 0, 0, 0.1), & [1,9]\times [0,1.5],\\
(1, 0, 0, 0.1), & [1,9]\times [1.5,3],
\end{cases}
\end{equation*}
subject to the free boundary conditions. The parameter $\beta$ in the monitor function is taken as $0.5$.
 
In the simulations, the BP limiters are activated only for AMM-FVM. Accordingly, we employ a fixed mesh resolution $N_x\times N_y=450\times 150$ and compare the numerical solutions obtained on MM and UM in Figure \ref{Ex:TP}, at different output times $t=1$, $t=2.5$, and $t=4$. For clarity, the mesh is plotted using every other grid point in each coordinate direction, while all computations are performed on the full mesh. Benefiting from mesh adaption, Figure \ref{Ex:TP} clearly illustrates that the AMM-FVM captures significantly sharper shock waves than the UM-FVM with the same number of cells. Furthermore, the FVM on adaptive moving mesh exhibits higher resolution than that on uniform mesh, enabling the AMM-FVM to resolve Kelvin-Helmholtz instabilities; see the rollup structures.

\begin{figure}[ht!]
\centering
\includegraphics[width=0.32\textwidth]{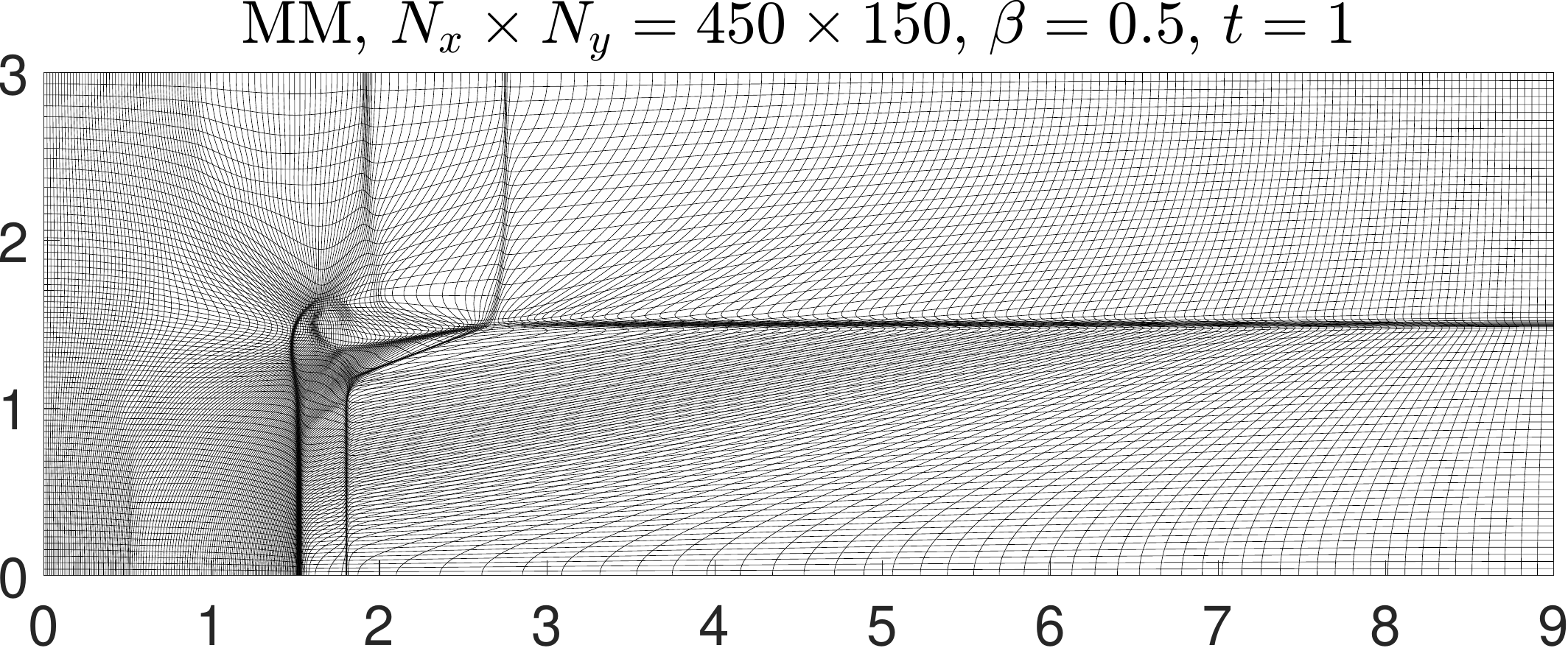}\quad
\includegraphics[width=0.32\textwidth]{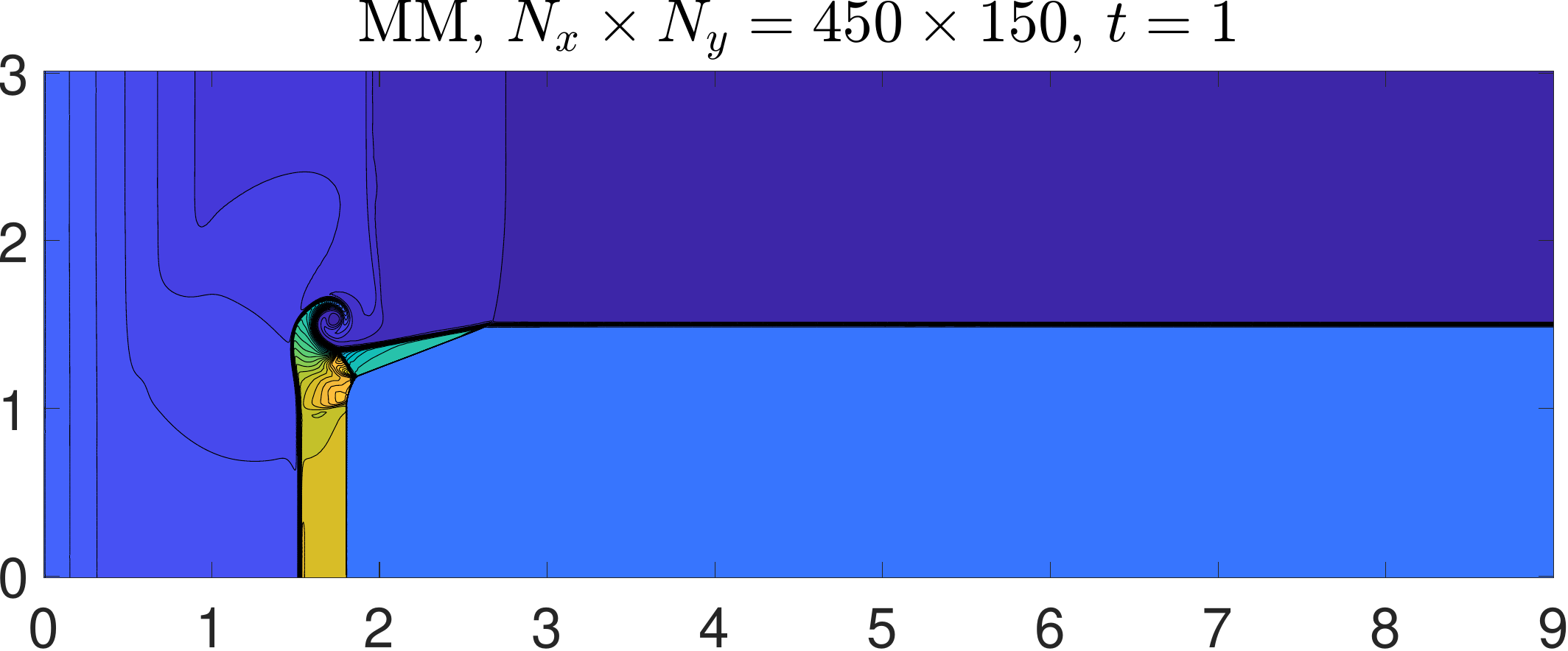}\quad
\includegraphics[width=0.32\textwidth]{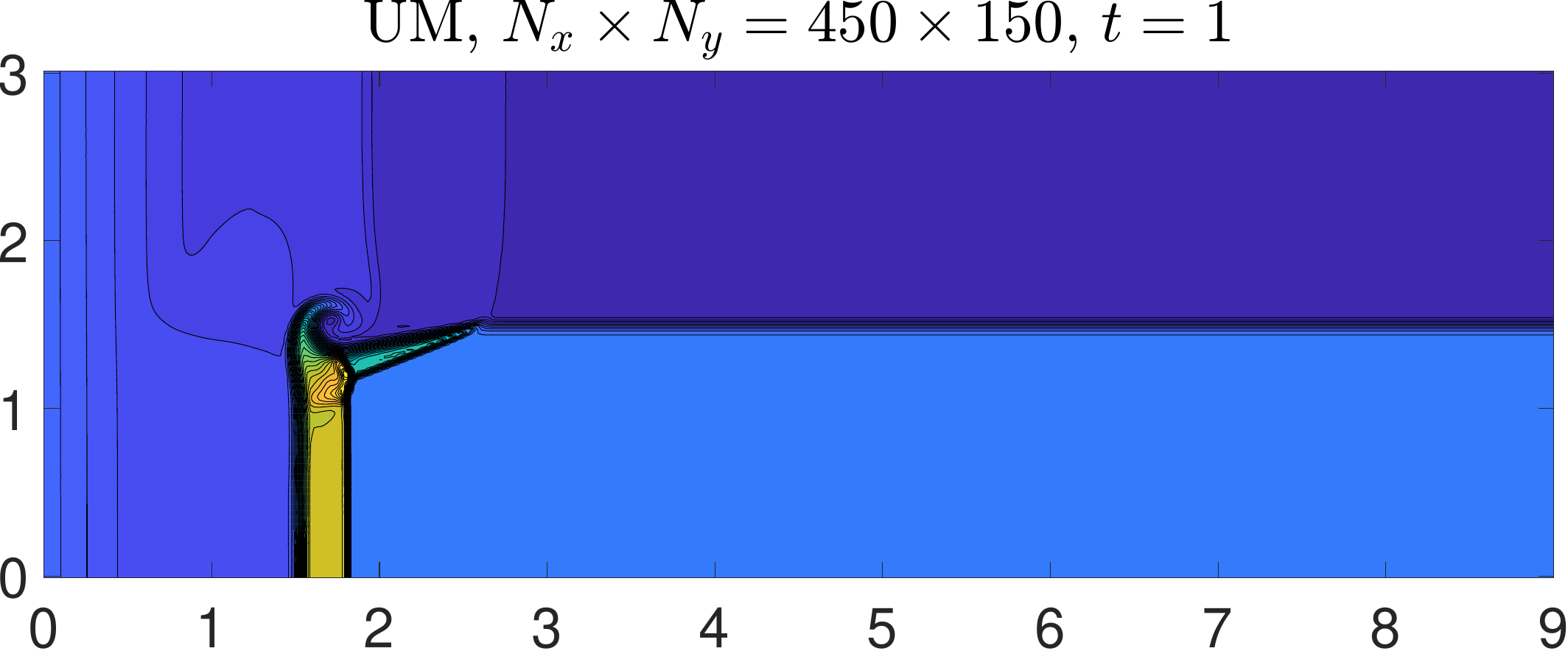}\vspace{+10pt}\\

\includegraphics[width=0.32\textwidth]{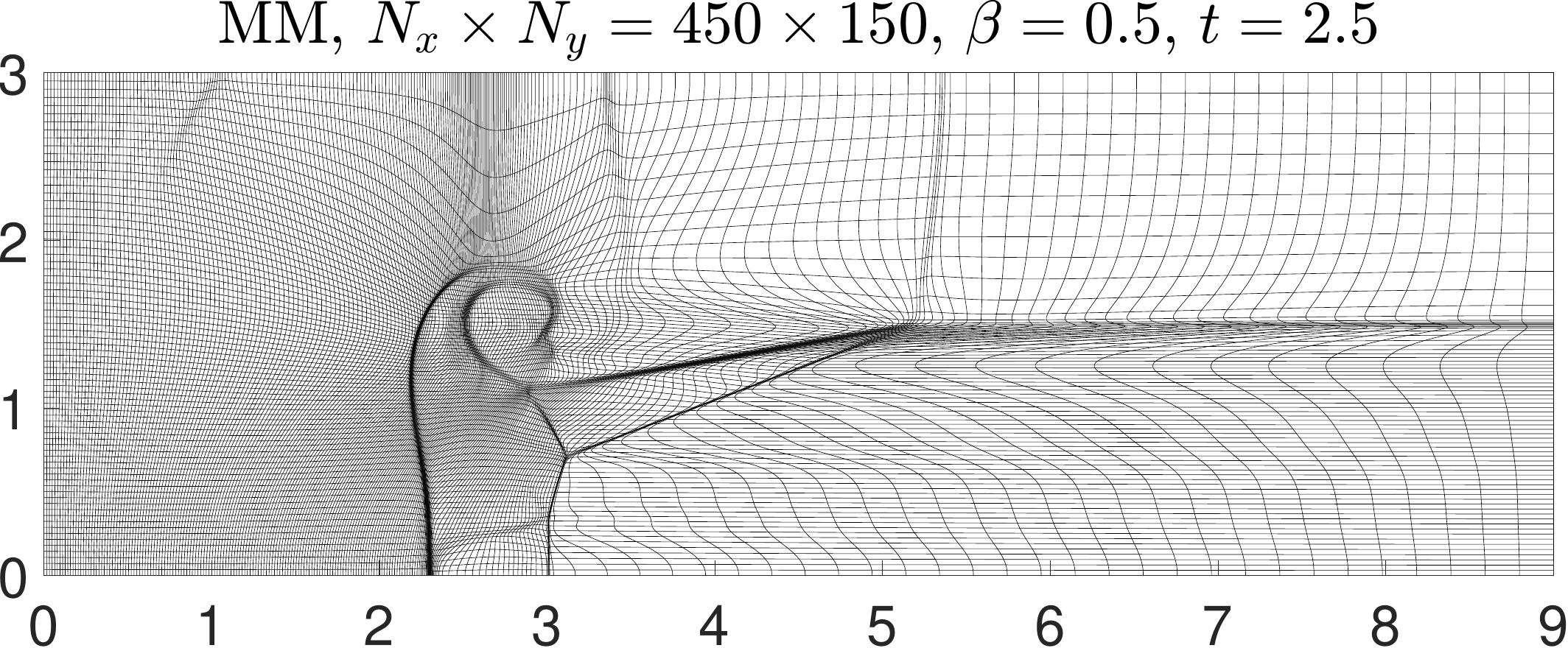}\quad
\includegraphics[width=0.32\textwidth]{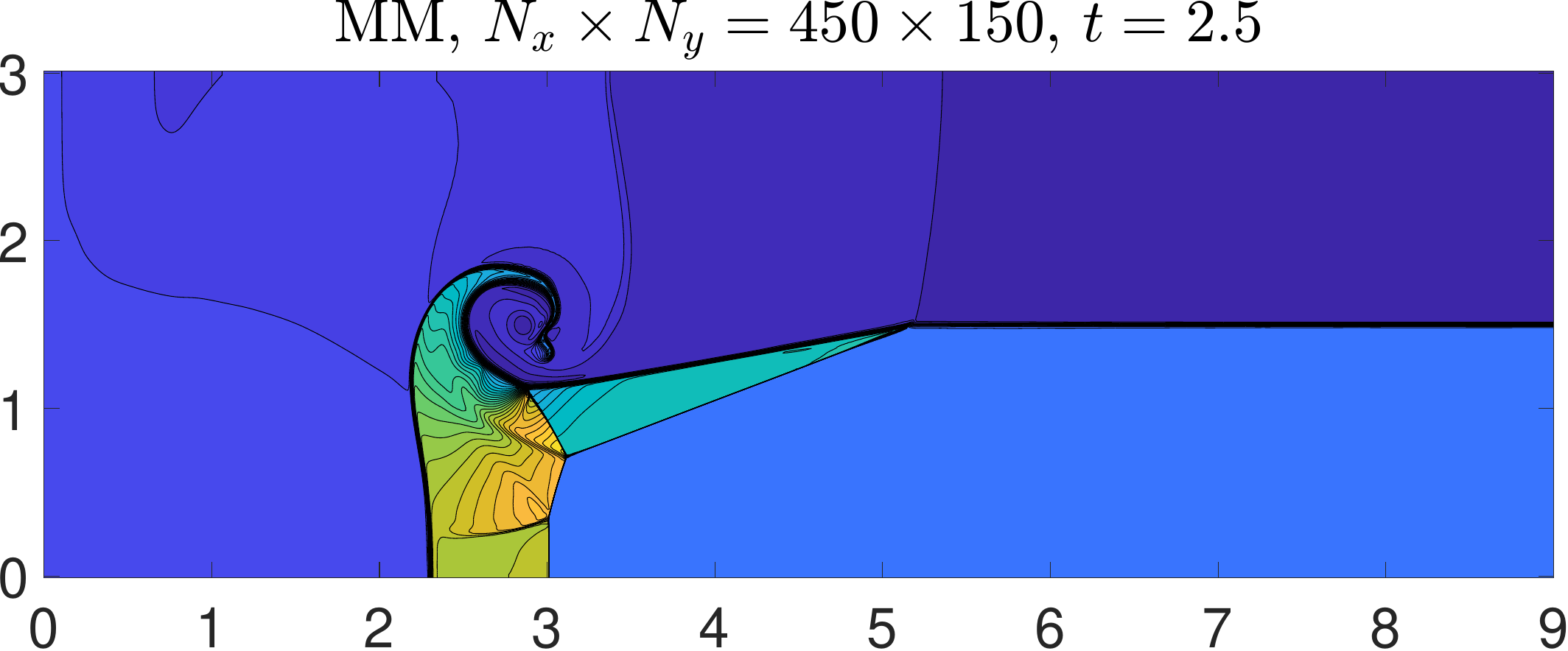}\quad
\includegraphics[width=0.32\textwidth]{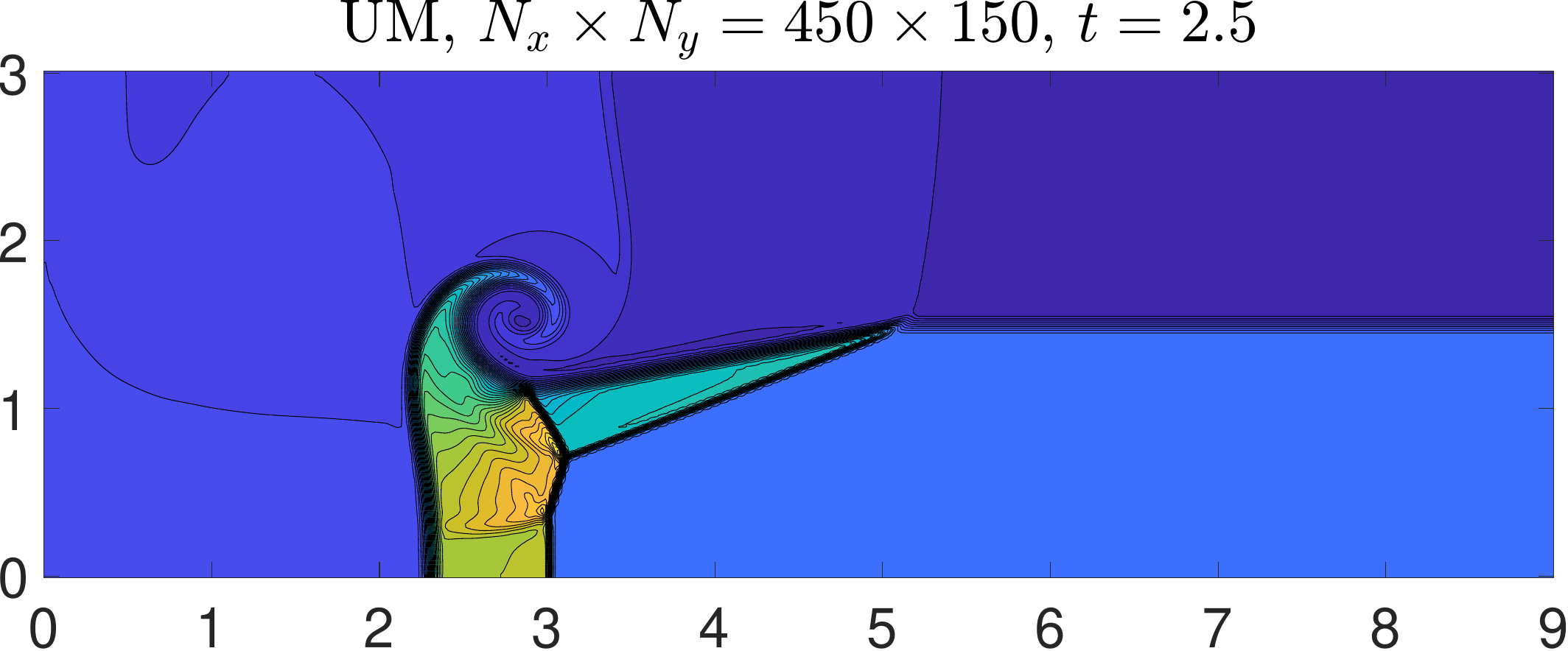}\vspace{+10pt}\\

\includegraphics[width=0.32\textwidth]{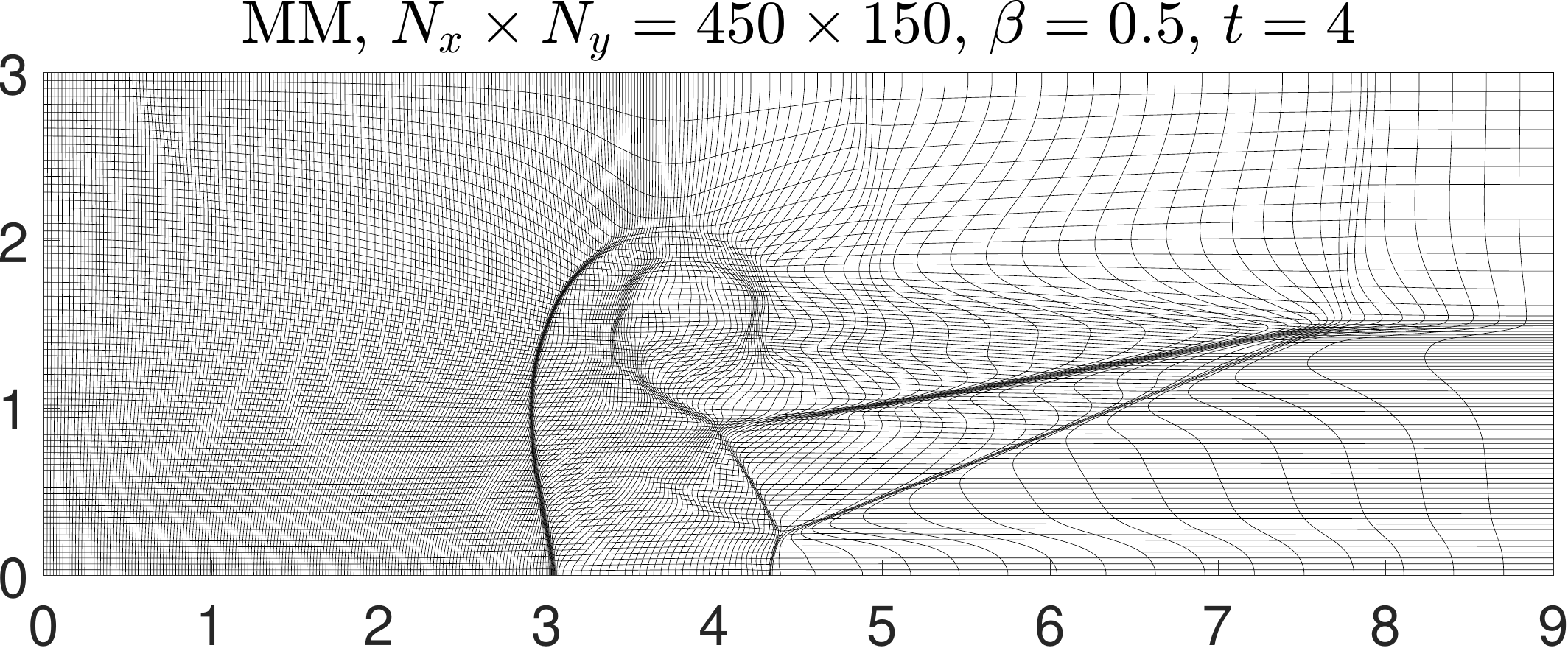}\quad
\includegraphics[width=0.32\textwidth]{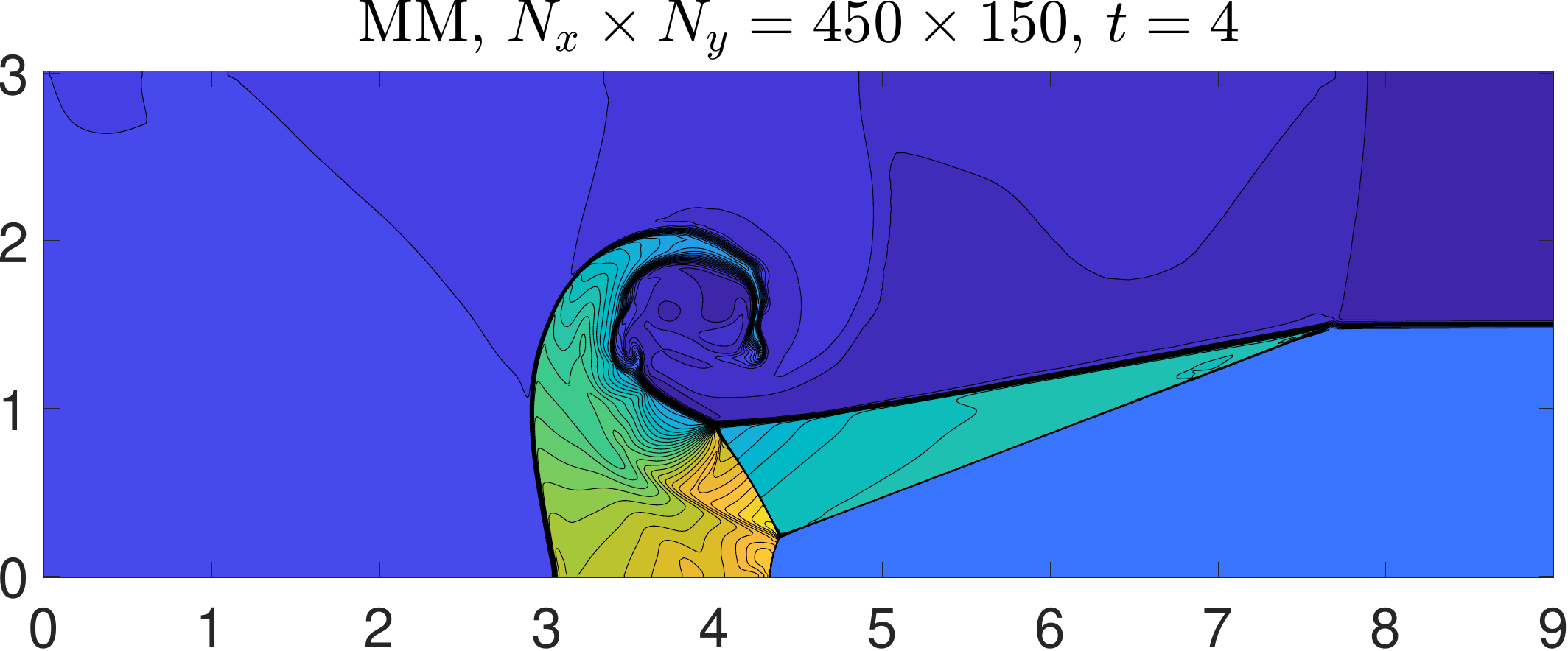}\quad
\includegraphics[width=0.32\textwidth]{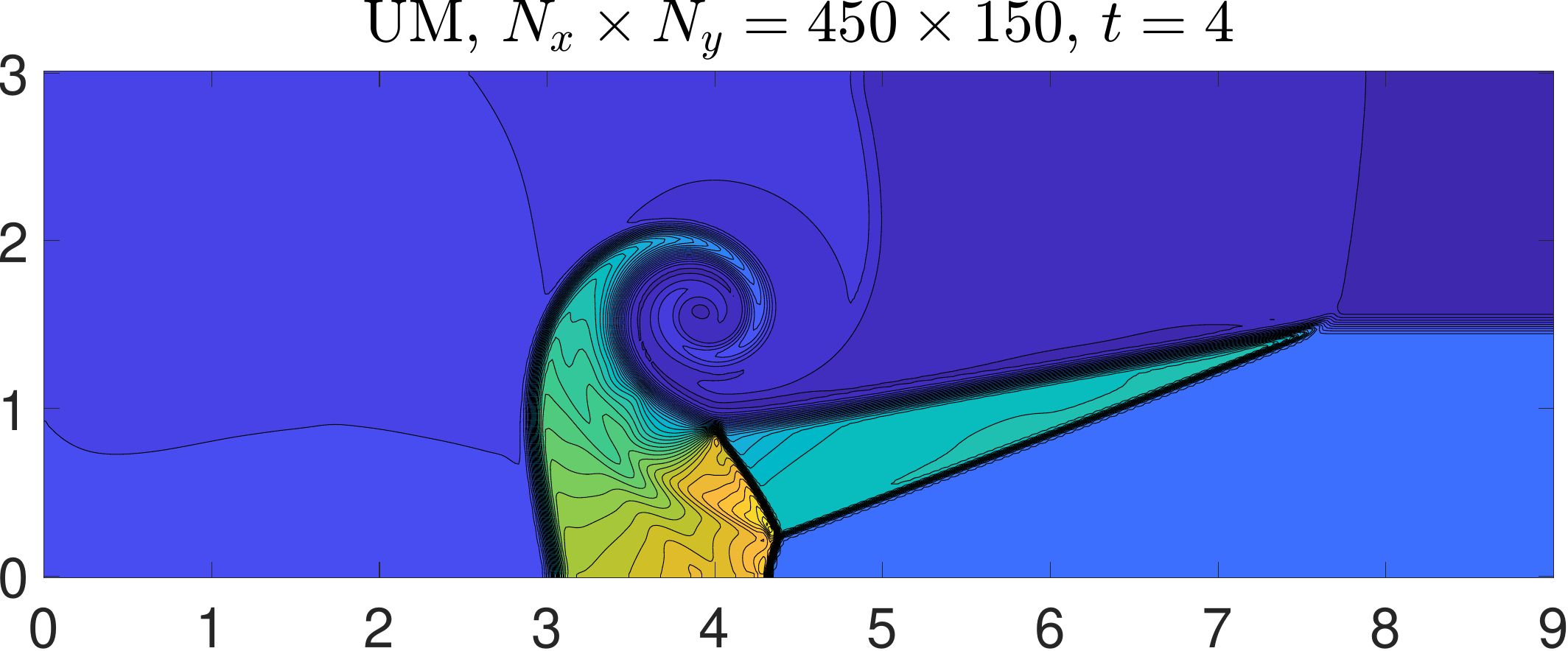}
\caption{\sf Example \ref{Ex:TP}: Mesh distribution (first column) and solution (second column) of the BP AMM-FVM, and solution of FVM on UM without BP limiters (third column) at $t=1$ (first row), $t=2.5$ (second row), and $t=4$ (third row).
\label{Fig:SVI}}
\end{figure}

\begin{example}\label{Ex:BW} 
\em Blast wave problem of Euler equations.
\end{example}
In this example, the robustness and efficiency of the BP AMM-FVM are validated through the blast wave problem of Euler equations taken from \cite{Zhang2012Maximum}. The computational domain is $[0,1.2]^2$. Initially, the density is $1$, the velocities $u$ and $v$ are zero, and the total energy is $10^{-12}$ everywhere except at the bottom left cell, where it is set to $\nicefrac{0.244816}{|C^\textrm{UM}|}$, with $|C^\textrm{UM}|$ denoting the cell area for UM. Reflective boundary conditions are applied on the left and bottom boundaries, while free boundary conditions are used elsewhere. The simulation is run up to the final time $t=1$.

\begin{figure}[ht!]
\centering
\includegraphics[width=0.31\textwidth]{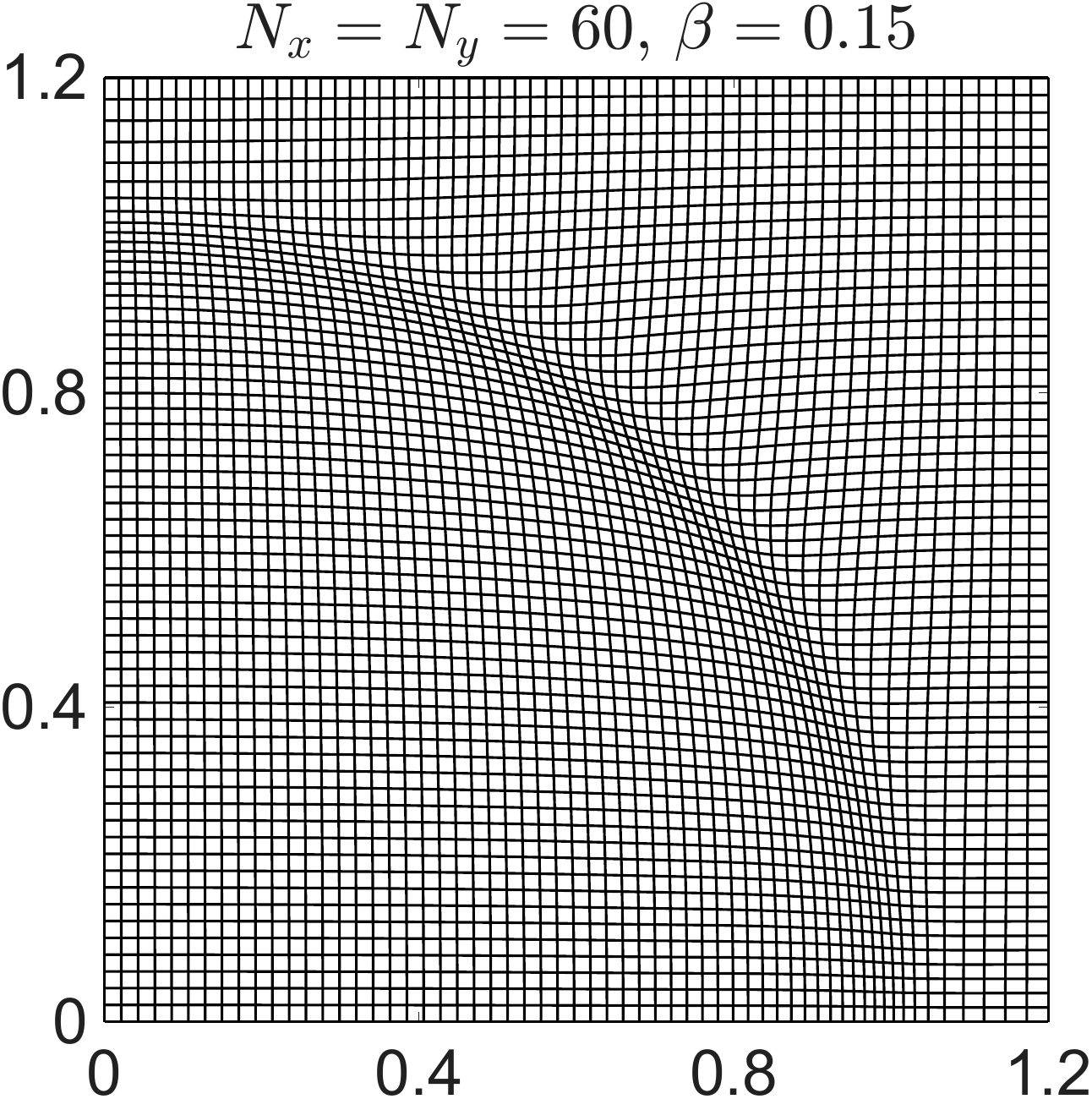}\quad
\includegraphics[width=0.31\textwidth]{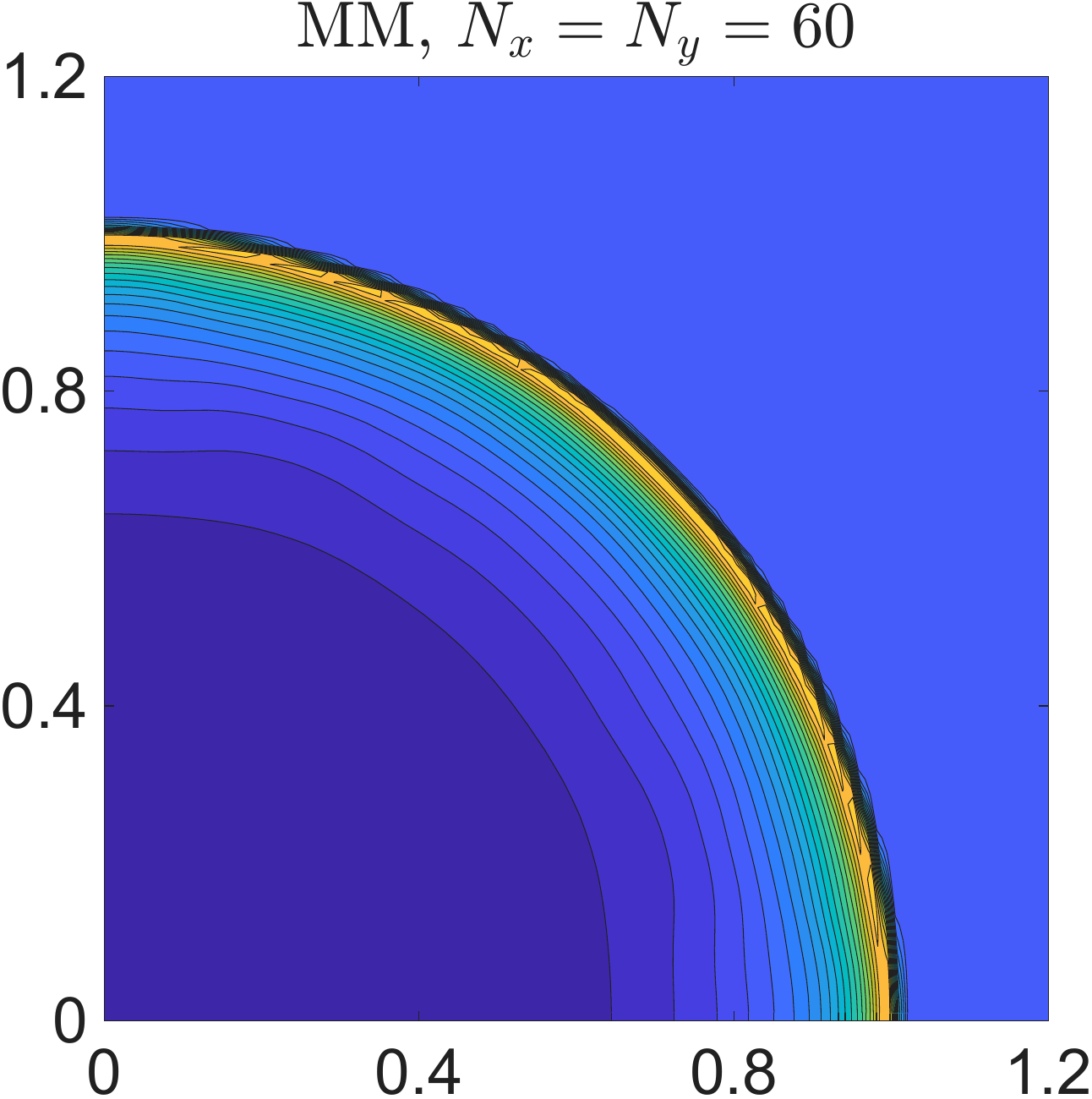}\quad
\includegraphics[width=0.31\textwidth]{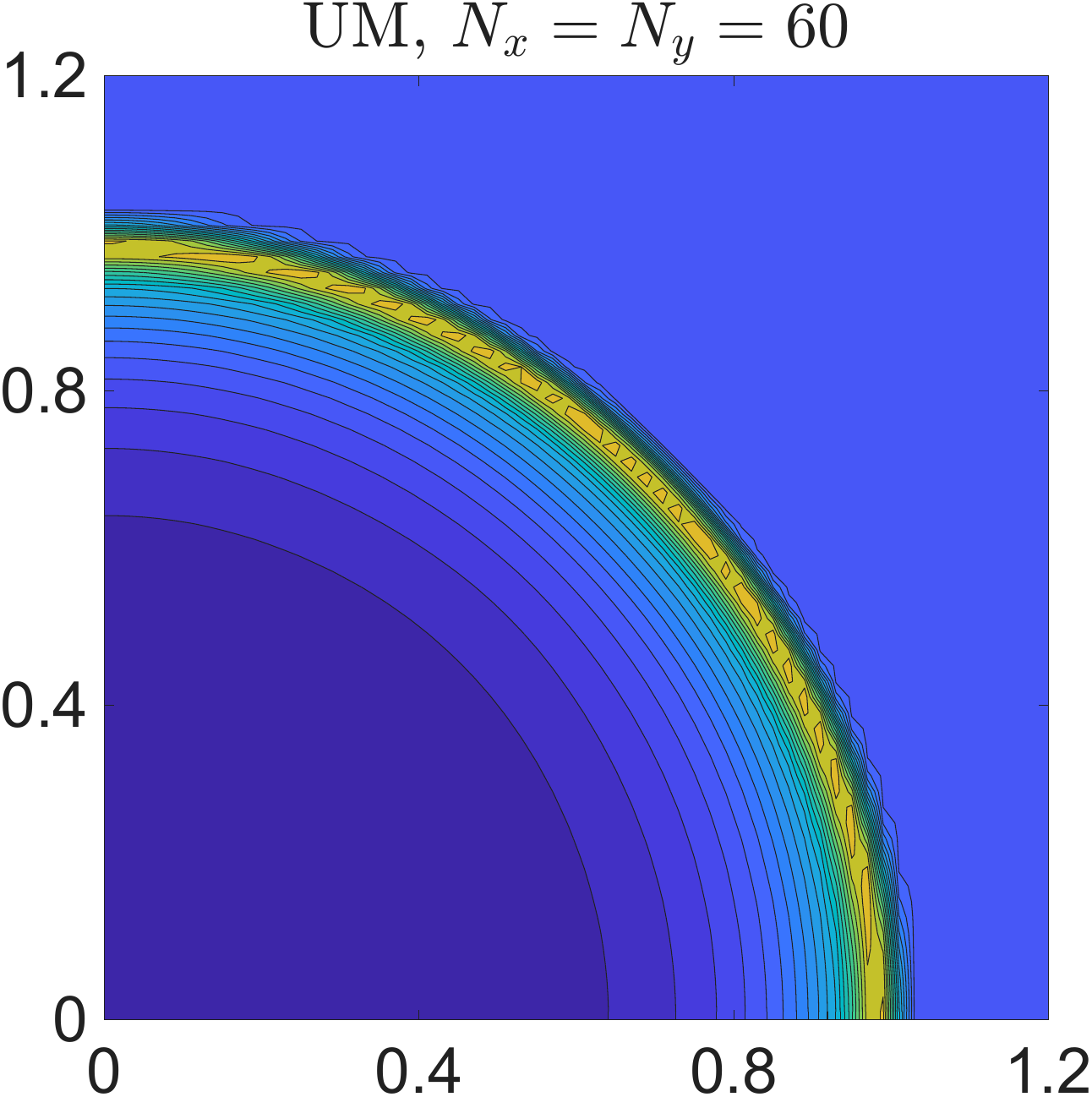}
\caption{\sf Example \ref{Ex:BW}: Mesh distribution (first column) and solutions of the BP FVM on MM (second column) and UM (third column) at $t=1$, Mesh resolution is $N_x\times N_y = 60\times 60$.
\label{Fig:BW}}
\end{figure}

In the computations, we set $\beta=0.15$. BP limiters are employed for FVMs on both MM and UM.
Since the initial condition depends on the cell area, a fixed mesh resolution $N_x\times N_y=60\times 60$ is used to compare numerical results obtained on both MM and UM. 
%The CPU time is $18$h $17$min for MM and $17$h $53$min for UM. The slightly larger CPU time for MM arises from mesh adaptivity. 
As observed in Figure \ref{Fig:BW}, the grid points concentrate adaptively to the shock wave. Consequently, the shock wave obtained on MM is remarkably sharper than that obtained on UM, which verifies the efficiency of the moving mesh methods. 

\begin{example}\label{Ex:SVI} 
\em Shock-vortex interaction problem of Euler equations.
\end{example}

In this example, the robustness and efficiency of the BP AMM-FVM is demonstrated by testing the shock-vortex interaction problem of Euler equations, taken from \cite{Cui2024On}. The computational domain is $[0,2]\times [0,1]$. Initially, a steady shock of Mach number $M_s = 1.1$ is positioned at $x=0.5$, with the upstream state on the left given by $(\rho_l,u_l,v_l,p_l) = (1, M_s\gamma, 0, 1)$, and the downstream state $(\rho_r,u_r,v_r,p_r)$ on the right given by 
\begin{equation}
\frac{\rho_r}{\rho_l} = \frac{(\gamma+1)M_s^2}{2+(\gamma-1)M_s^2}, \qquad
\frac{u_r}{u_l} = \frac{2+(\gamma-1)M_s^2}{(\gamma+1)M_s^2}, \qquad
v_r = 0, \qquad
\frac{p_r}{p_l} = 1 + \frac{2\gamma}{\gamma+1}\left(M_s^2 - 1\right).
\end{equation}
Then, an isentropic vortex is imposed to the left state. The perturbations to the velocities $(\delta u, \delta v)$, temperature $T_l = p_l/\rho_l$, and entropy $s_l = \ln(p_l/\rho_l^\gamma)$ associated with the vortex are given by
\begin{equation}
(\delta u, \delta v) = \frac{\varepsilon}{r_c} e^{\alpha(1-\tau^2)} (\bar{y}, -\bar{x}), \qquad
\delta T = -\frac{(\gamma-1)\varepsilon^2}{4\alpha\gamma} e^{2\alpha(1-\tau^2)}, \qquad
\delta s = 0,
\end{equation}
with $\tau = \frac{1}{r_c}\sqrt{(x-x_c)^2+(y-y_c)^2}$, where $x_c = 0.25$ and $y_c=0.5$. The rest of the parameters are set as follows. The strength of the vortex $\varepsilon$ is set as $1.378106$, the decay rate of the vortex $\alpha$ is set as $0.204$, and the critical radius $r_c=0.4$. Eventually, the initial state to the left of the shock positioned at $x=0.5$ is 
\begin{equation}
\rho = \left(\frac{T_l+\delta T}{s_l+\delta s}\right)^{\frac{1}{\gamma-1}}, \qquad
u = u_l + \delta u, \qquad
v = v_l + \delta v, \qquad
p = \rho(T_l+\delta T).
\end{equation}

In this example, BP limiters are required for FVMs on both UM and MM. Results obtained on MM are compared with those on UM in Figure \ref{Fig:SVI}. The mesh resolutions for MM and UM are adjusted such that the CPU times of the BP FVMs are comparable (approximately $46$min). It can be observed from Figure \ref{Fig:SVI} that the BP AMM-FVM captures significantly sharper shock waves than the BP FVM on the uniform mesh, while using far fewer grid points. This behavior is not supervising, as it stems from the efficient adaptive movement of the grid points towards shock waves.

\begin{figure}[ht!]
\centering
\includegraphics[width=0.3\textwidth]{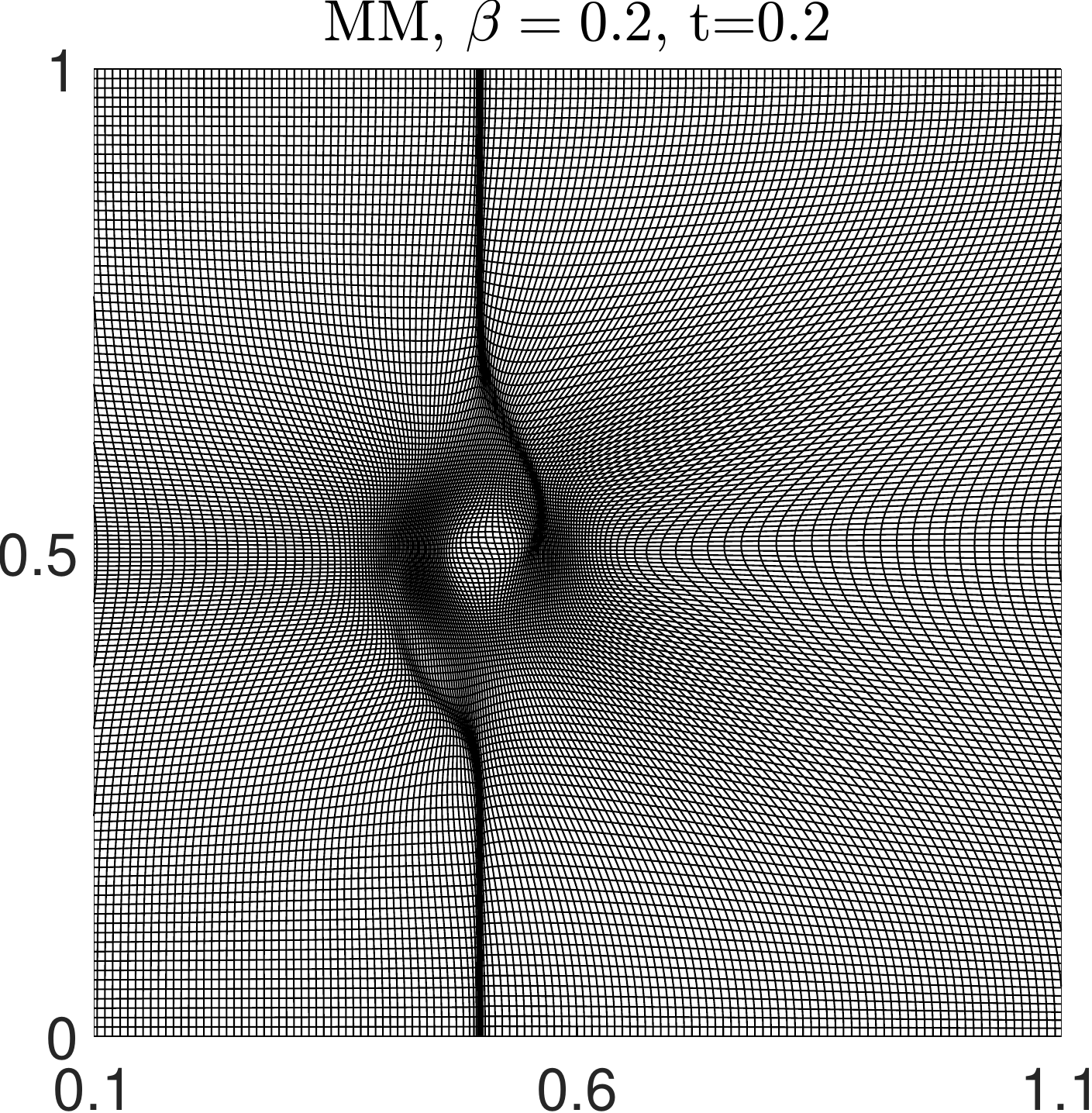}\quad
\includegraphics[width=0.3\textwidth]{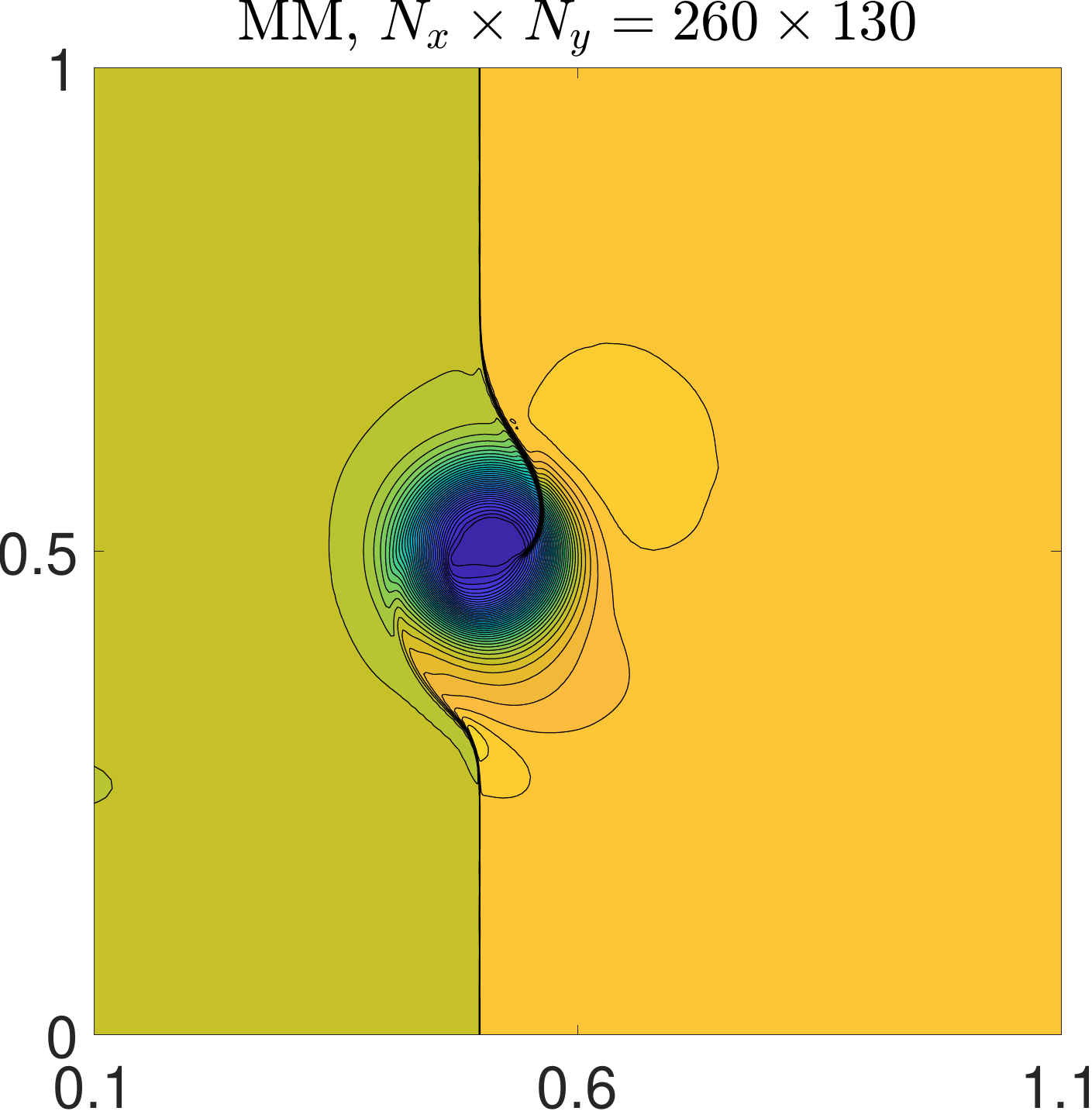}\quad
\includegraphics[width=0.3\textwidth]{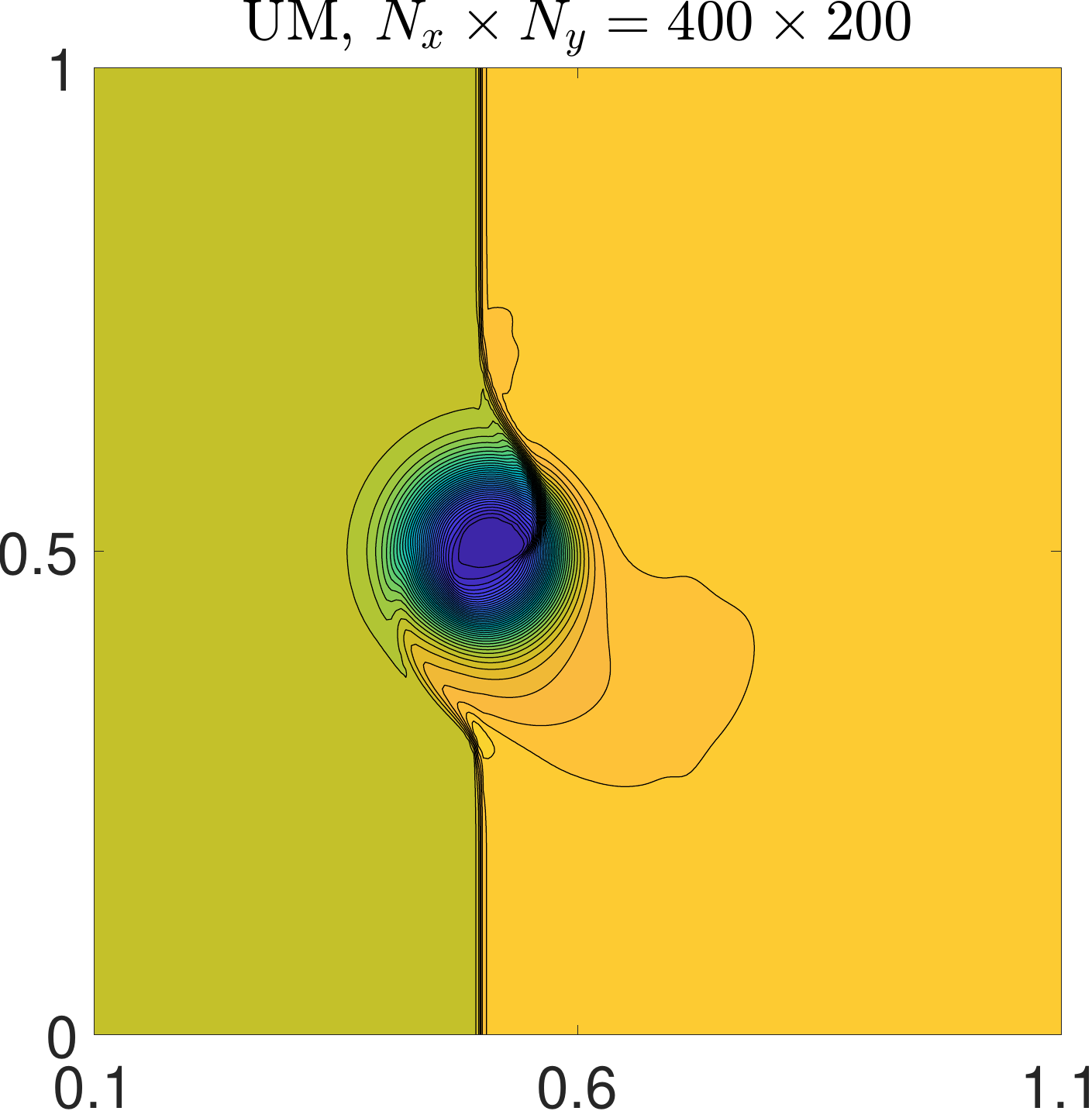}\vspace{+10pt}\\

\includegraphics[width=0.3\textwidth]{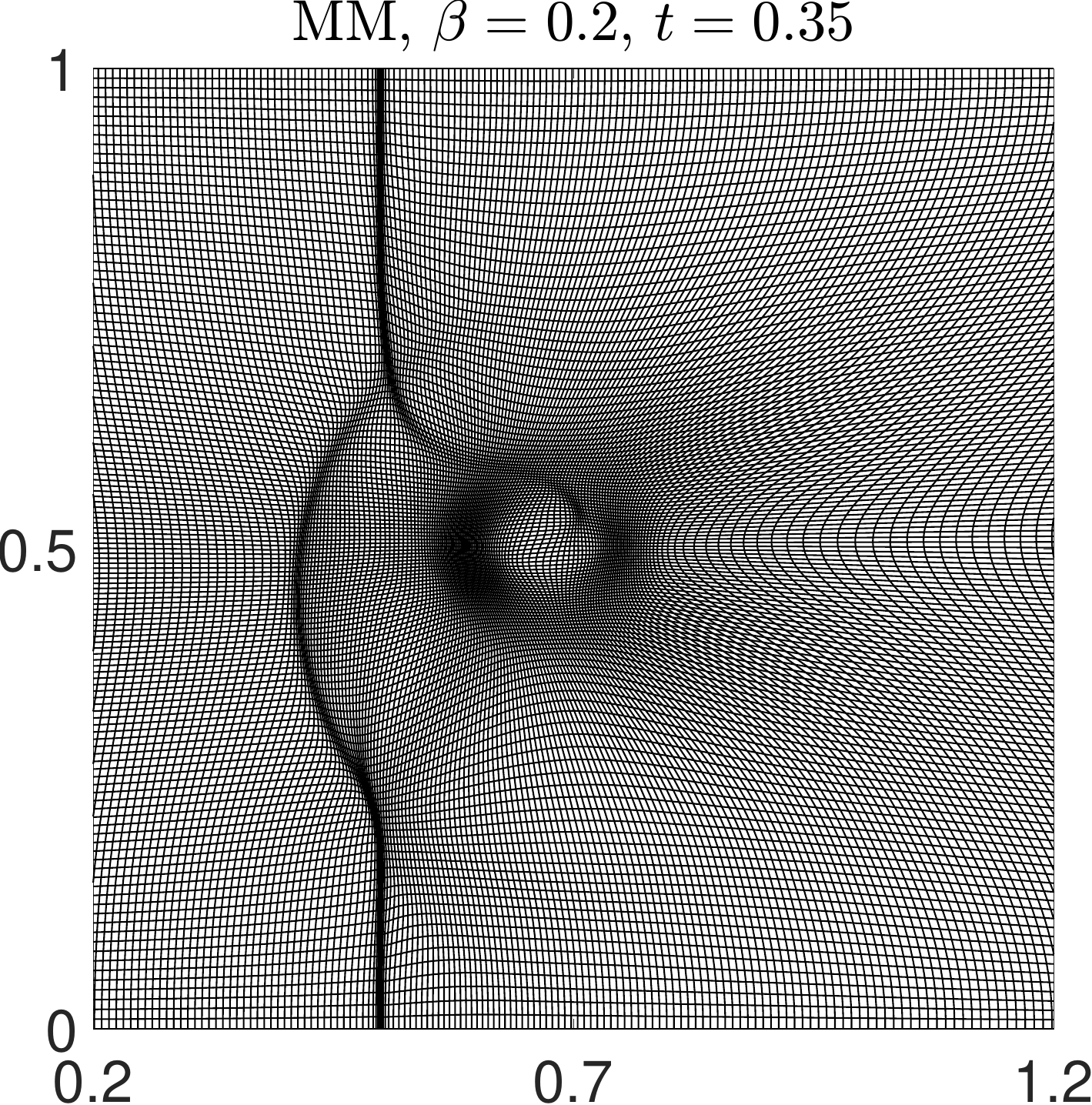}\quad
\includegraphics[width=0.3\textwidth]{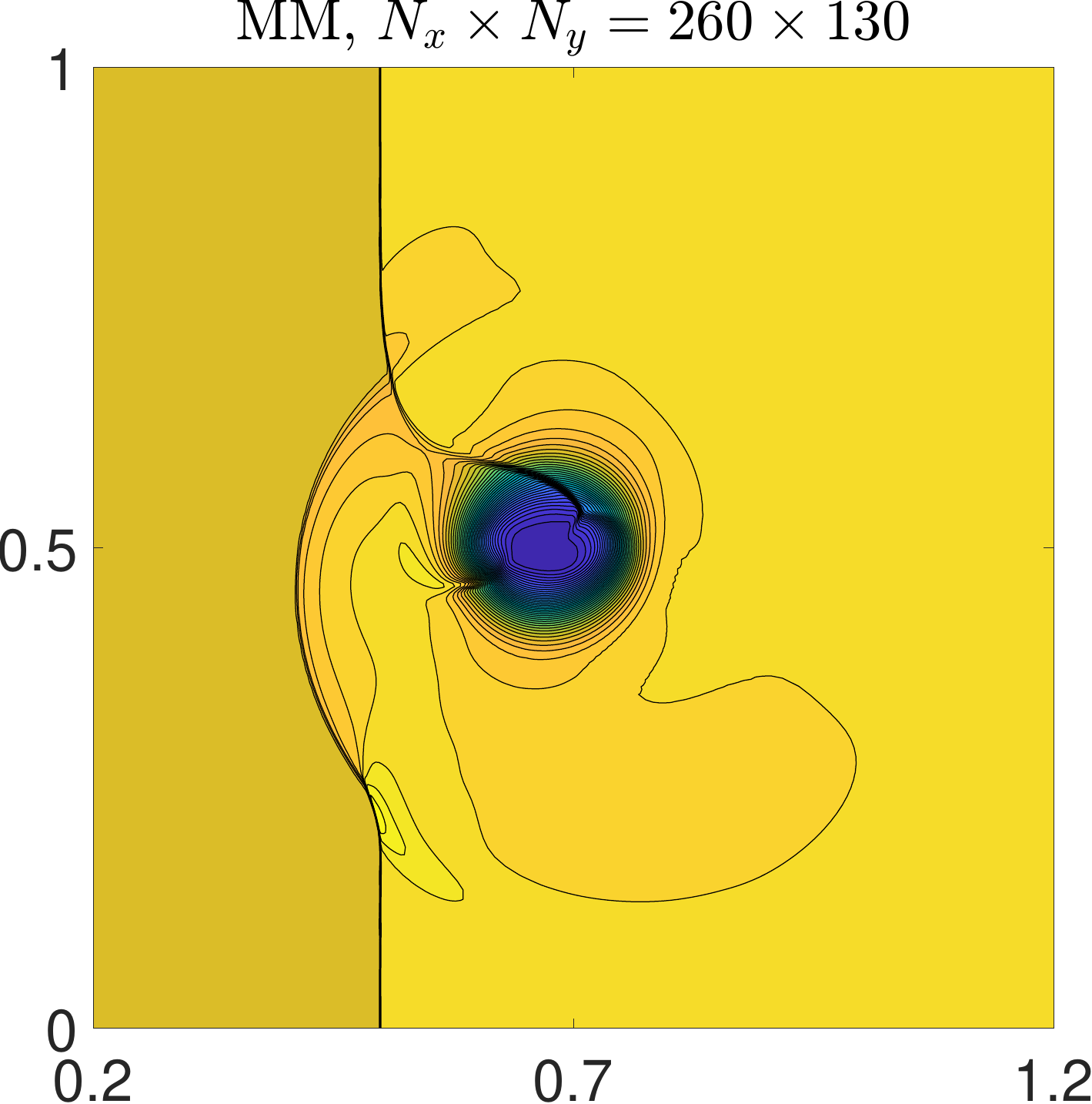}\quad
\includegraphics[width=0.3\textwidth]{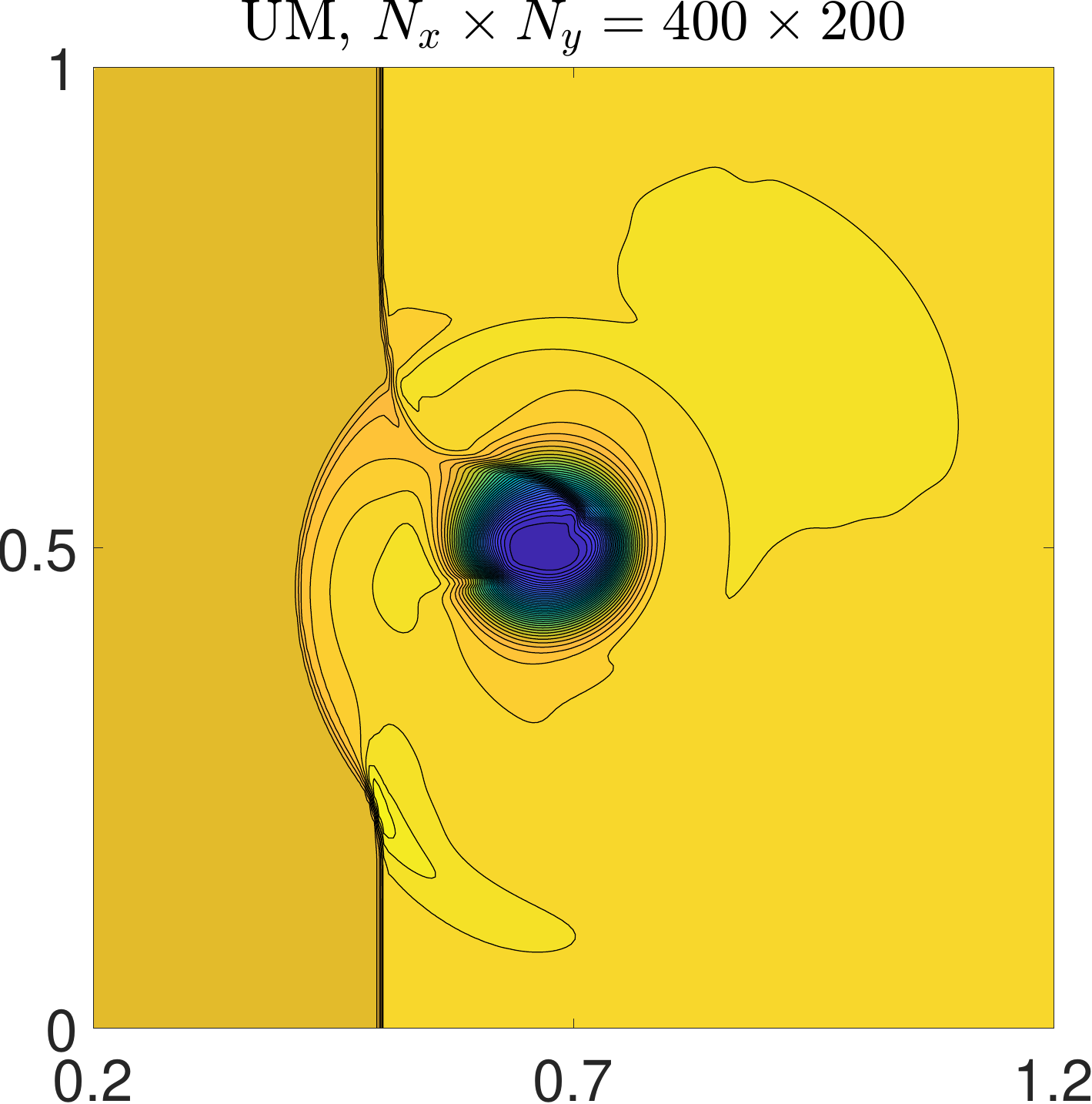}\vspace{+10pt}\\

\includegraphics[width=0.3\textwidth]{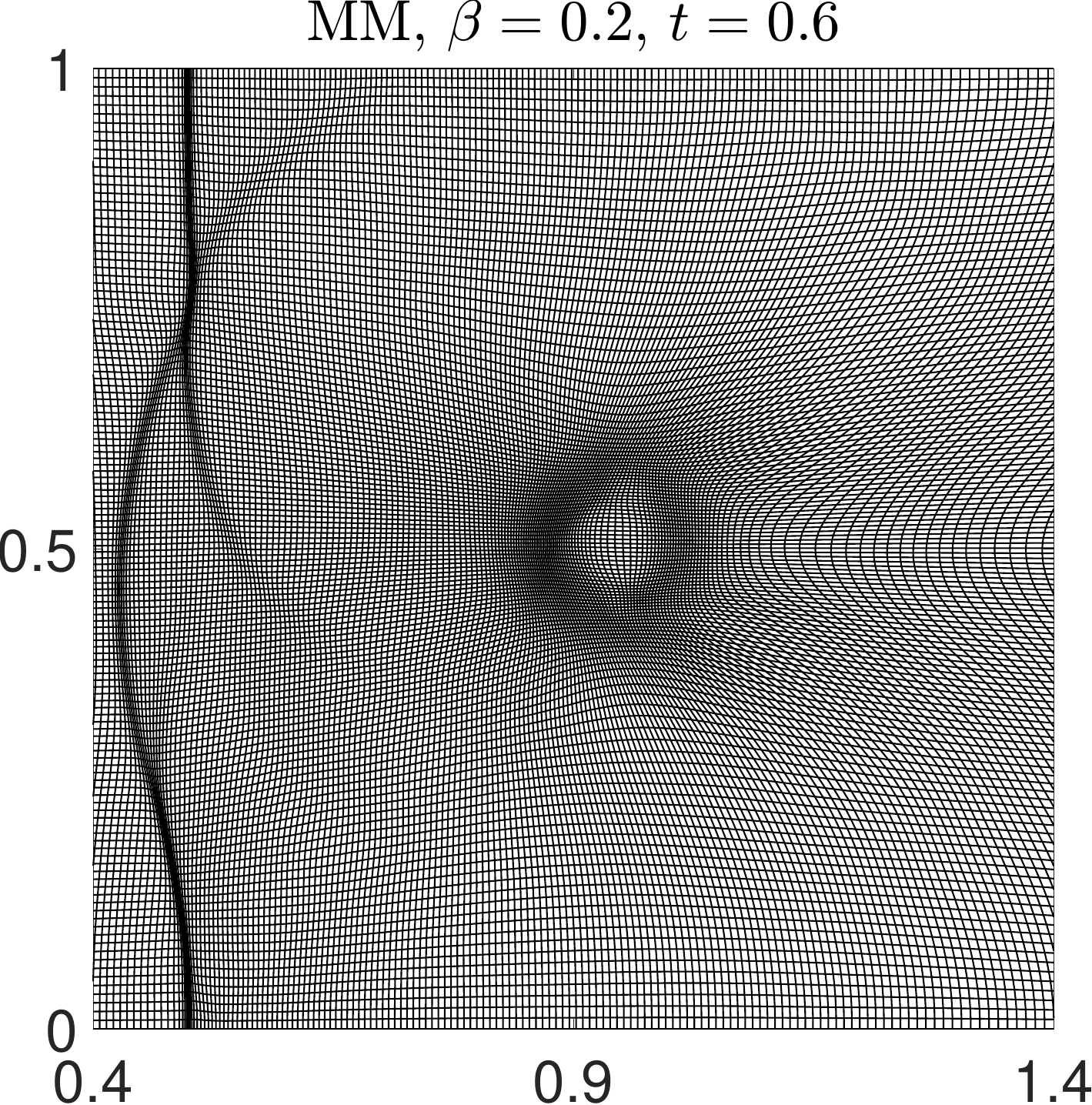}\quad
\includegraphics[width=0.3\textwidth]{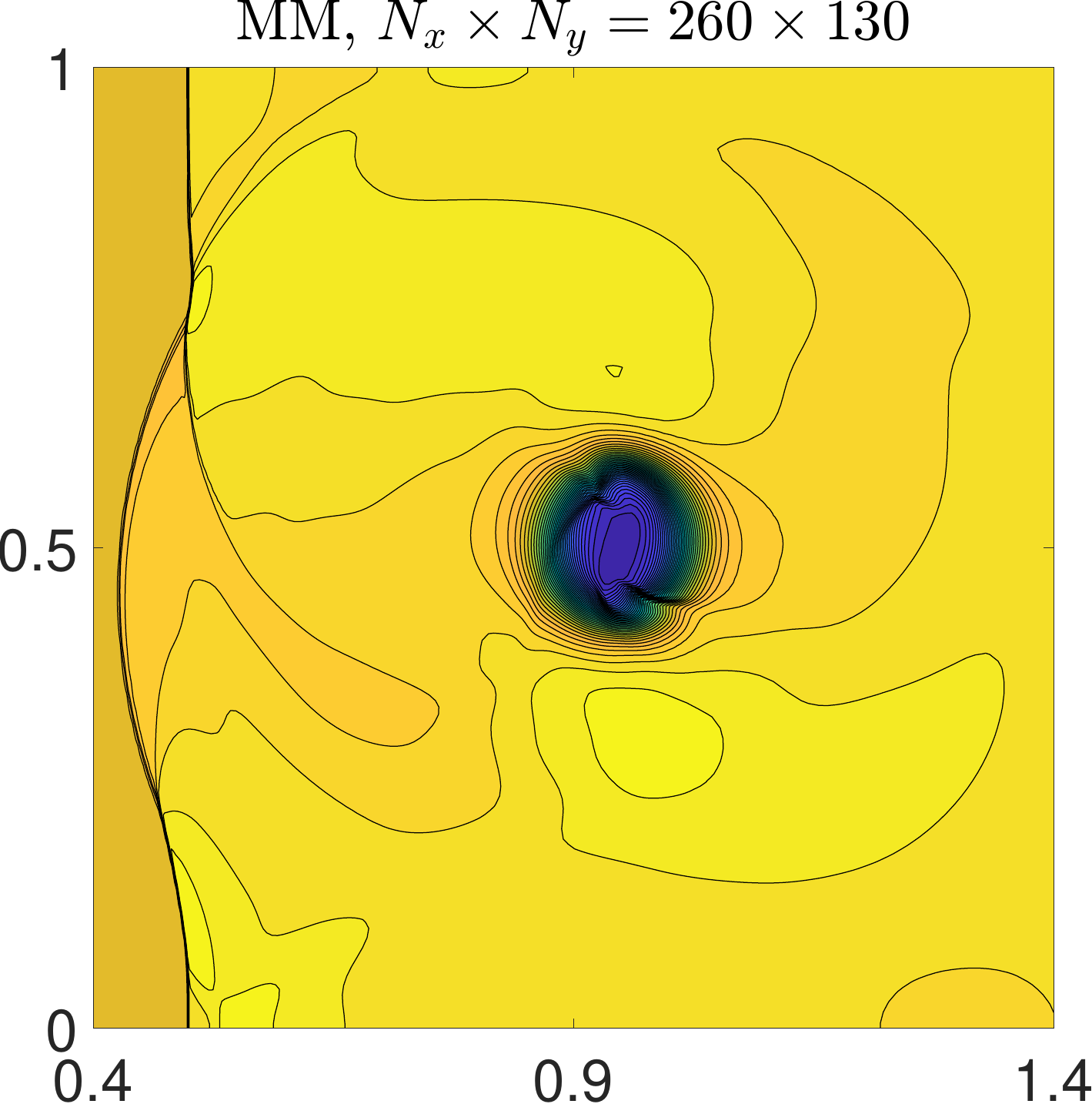}\quad
\includegraphics[width=0.3\textwidth]{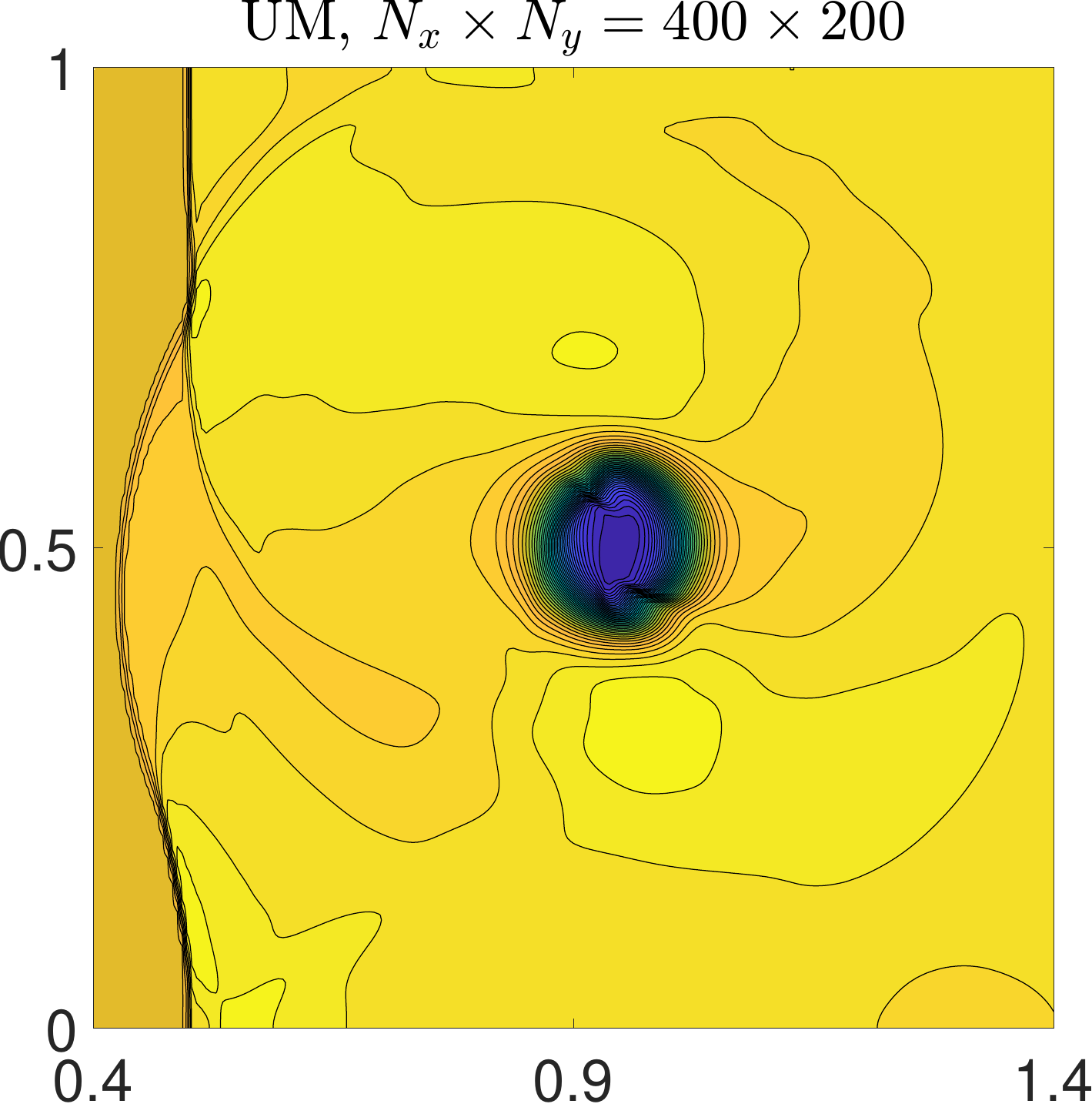}
\caption{\sf Example \ref{Ex:SVI}: Mesh distribution (first column) and solution (second column) of the BP AMM-FVM, and solution of BP FVM on UM (third column) at $t=0.2$ (first row), $t=0.35$ (second row), and $t=0.6$ (third row).
\label{Fig:SVI}}
\end{figure}

\section{Conclusion}
In this paper, we have developed a new framework of bound-preserving finite volume method on adaptive moving meshes for two-dimensional hyperbolic conservation laws. In the discretization part, the solution is directly evolved from $t^n$ to the subsequent time level $t^{n+1}$ via the standard finite volume method on adaptive moving meshes. To achieve the BP property, the high-order AMM-FVM is first reformulated as a convex combination of several high-order sub-cell schemes, and then the BP property of each sub-cell scheme is enforced by blending it with its first-order counterpart. It has been demonstrated that this procedure essentially introduces a BP flux limiter, which allows each high-order numerical flux on the cell boundary to be modified by blending with its first-order counterpart so as to enforce the BP property. A variety of two-dimensional numerical examples have been presented to illustrate and validate the accuracy, robustness, and effectiveness of the proposed BP AMM-FVM.

It is worth noting that the BP framework proposed in this paper is suitable for all hyperbolic equations that can be written in fully conservative forms. Extensions of the present BP framework to hyperbolic systems with non-conservative products and to three-dimensional problems will be investigated in our future work.
 
%\newpage
\section*{Acknowledgement} 
The authors would like to thank Prof. Shumo Cui and Prof. Shengrong Ding from Sun Yat-Sen University for insightful discussions regarding Example \ref{Ex:SVI} of this paper.
 
%\section*{CRediT authorship contribution statement}
%\textbf{Yaguang Gu}: 
%\textbf{Guanghui Hu}: 

\bibliographystyle{siam}
%\bibliography{Reference}

\end{document}